\documentclass[10pt]{article}
\usepackage{amsmath,amssymb,amsfonts,amsthm,amscd,graphicx,psfrag,epsfig}
\usepackage{color}
\definecolor{blue}{rgb}{0,0,0.7}
\definecolor{red}{rgb}{0.75, 0, 0}
\usepackage{stmaryrd}
\usepackage[titletoc,toc]{appendix}
\usepackage{hyperref}
\hypersetup{colorlinks, linkcolor=blue, citecolor=cyan}
\usepackage[all]{xy}           
\usepackage{marginnote}  
\usepackage{float}
\usepackage{soul}

\newtheorem{theorem}{Theorem}[section]
\newtheorem{lemma}[theorem]{Lemma}
\newtheorem{proposition}[theorem]{Proposition}
\newtheorem{corollary}[theorem]{Corollary}
\newtheorem{conjecture}[theorem]{Conjecture}
\newtheorem{definition}[theorem]{Definition}

\newcommand{\bpf}{\begin{proof}}
\newcommand{\epf}{\end{proof}}
\newcommand{\bcr}{\begin{color}{red}}
\newcommand{\ecr}{\end{color}}
\newcommand{\bs}{\begin{split}}
\newcommand{\es}{\end{split}}
\newcommand{\be}{\begin{equation}}
\newcommand{\ee}{\end{equation}}
\newcommand{\bt}{\begin{theorem}}
\newcommand{\et}{\end{theorem}}
\newcommand{\bd}{\begin{definition}}
\newcommand{\ed}{\end{definition}}
\newcommand{\bp}{\begin{proposition}}
\newcommand{\ep}{\end{proposition}}
\newcommand{\bl}{\begin{lemma}}
\newcommand{\el}{\end{lemma}}
\newcommand{\bc}{\begin{corollary}}
\newcommand{\ec}{\end{corollary}}
\newcommand{\bcon}{\begin{conjecture}}
\newcommand{\econ}{\end{conjecture}}
\newcommand{\la}{\label}
\newcommand{\B}{{\rm B}}
\newcommand{\A}{{\rm A}}

\newcommand{\Z}{{\mathbb Z}}
\newcommand{\R}{{\mathbb R}}

\newcommand{\C}{{\mathbb C}}

\newcommand{\hra}{\hookrightarrow}
\newcommand{\lra}{\longrightarrow}
\newcommand{\lms}{\longmapsto}
\newcommand{\bS}{{\Bbb S}}

\usepackage{yfonts}
 
\begin{document}

 \setcounter{tocdepth}{2}

\date{}
 
 \title{Non-commutative Cluster Lagrangians}

 \author{Alexander B. Goncharov, Maxim Kontsevich}
\maketitle

 \tableofcontents

 \begin{abstract}

   The space  ${\rm Loc}_m(S)$ of rank $m$ flat bundles on a closed surface $S$ is $K_2-$symplectic. 
   A threefold $M$ bounding $S$ gives rise a $K_2-$Lagrangian ${\cal L}_M$ in ${\rm Loc}_m(S)$ given by the   flat bundles on $S$ extending to  $M$. 
   We generalize this, replacing the zero section in $T^*M$ by certain  singular Lagrangians in $T^*M$. 
  \vskip 1mm

    First,  we  introduce {\it ${\cal Q}-$diagrams} in threefolds. 
    They are collections  ${\cal Q}$ of smooth cooriented surfaces $\{S_i\}$, intersecting transversally everywhere 
    but in a finite set of quadruple crossing points. We require that  shifting any surface $S_i$  from such a point    in  the direction of its coorientation creates a simplex with the  cooriented out  faces.     
   The ${\cal Q}-$diagrams are 3d analogs of bipartite ribbon graphs.     \vskip 1mm
    
    Let ${\Bbb L}$ be the Lagrangian in $T^*M$ given by the union of the zero section and the conormal bundles to the cooriented surfaces $S_i$ of ${\cal Q}$.  Let ${\cal X}_{\Bbb L}$ be  the 
    stack of      admissible dg-sheaves ${\cal F}$ on $M$ with the microlocal support in  ${\Bbb L}$, whose microlocalization  at $T^*_{S_i}M-0$ is   a rank one  local system.
    \vskip 1mm
    
    We introduce the boundary  $\partial {\Bbb L}$ of  ${\Bbb L}$. It is a singular Lagrangian in a symplectic space, providing a symplectic  stack  ${\cal X}_{\partial \Bbb L}$, 
    and  a restriction functor ${\cal X}_{\Bbb L} \lra {\cal X}_{\partial \Bbb L}$. The image of the latter  is  Lagrangian. 
    We show that, under mild conditions on ${\cal Q}\cap \partial M$, this Lagrangian  has a  cluster 
    description, and so  it is a $K_2-$Lagrangian.  It  also has a simple description in the   non-commutative setting.     
           
   
\end{abstract}

    \section{Introduction}

    \subsection{${\cal Q}-$diagrams of cooriented hypersurfaces}  \la{sect3aa} 
    
     All Lemmas in Section \ref{sect3aa} are  evident.
    
    \subsubsection{${\cal Q}-$diagrams of hypersurfaces} \la{SEC3b}

         Let $X$ be a manifold. A {\it coorientation} of a hypersurface $H$ in $X$ is  a choice of a connected component 
 of the  conormal bundle to $H$ minus the zero section. 
  This component is called the  conormal bundle to  $H$. If $X$ is oriented,  a coorientation of a hypersurface is equivalent to its orientation.

\bd  \la{L1a1}
A collection of smooth  cooriented hypersurfaces in an $m-$dimensional  manifold    is   a \underline{${\cal Q}-$diagram} if the intersection points where $k<m$ hypersurfaces meet are transversal intersections, and 
\begin{enumerate} 

\item The remaining isolated intersection points  $q$ are   
intersections of $m+1$ hypersurfaces.  

\item 

Shifting   any of the   hypersurfaces intersecting at $q$   in the direction of its coorientation  we create an $m-$dimensional  simplex with  the cooriented out  faces.

\end{enumerate}

\ed

 \noindent 
 
Condition (2) just means that   the  convex hull of the endpoints of  oriented conormals $h_0, ..., h_m$ 
 to the hypersurfaces intersecting at $q$ contains zero, i.e.  they  satisfy a single up to a constant linear relation 
\be \la{convc}
\alpha_0 h_0 +   \ldots  + \alpha_mh_{m}=0, \ \ \ \ \alpha_0, \ldots, \alpha_{m} >0. 
\ee

      \begin{figure}[ht]
\centerline{\epsfbox{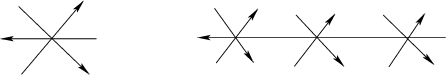}}
\caption{${\cal Q}-$diagrams  of oriented curves on a surface are triple crossing diagrams.} 
\label{ncl2}
\end{figure}  
 
\bl  A  ${\cal Q}-$diagram  
in a manifold  cuts a ${\cal Q}-$diagram on each   hypersurface   of the collection.
\el

 A ${\cal Q}-$diagram ${\cal Q}$ in $X$ provides a decomposition of $X-{\cal Q}$ into connected  domains ${\cal C}$. 
  A vertex $q\in {\cal C}$ is called {\it consistently cooriented} if the all faces  
  of ${\cal C}$ meeting at $q$ are either cooriented in, or cooriented out of ${\cal C}$.   
 A domain  is  {\it colored} if all its vertices  are consistently oriented, and  {\it mixed} if none   is consistently oriented. 
 A colored domain is a {\it $\circ-$domain} if its faces are cooriented in, and a {\it $\bullet-$domain} otherwise. 
 
 \bl \la{L2.3} Given a ${\cal Q}-$diagram ${\cal Q}$ in $X$, every component of $X-{\cal Q}$ is either colored, or mixed. 
  
A componsnt  in $X-{\cal Q}$ is colored if and only if all its boundary faces are colored by the same color.
\el

 \subsubsection{Alternating diagrams of hypersurfaces}  \la{2.0.2}  A  finite collection  of cooriented points on a line  is called an {alternating diagram} if traversing  the line  the coorientations of the points alternate.  
 
      \begin{figure}[ht]
\centerline{\epsfbox{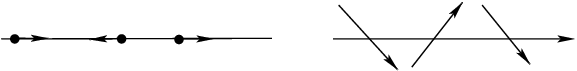}}
\caption{Alternating diagrams of points on a line, and of oriented curves on a surface.} 
\label{nc2025}
\end{figure} 
 
 \bd  \la{DEF2.1} An  \underline{alternating diagram} in a manifold $X$ is a normal crossing connected  
collection  ${\cal H}$ of smooth cooriented  hypersurfaces in  $X$ such that

 \begin{itemize}

 \item For any intersection line in ${\cal H}$ the coorientations of the hypersurfaces intersecting  
 it alternate. 
 \end{itemize}
 
 We require that an {alternating diagram}  intersects the boundary of the manifold $X$ transversally. 
  \ed

   An alternating diagram on $X$  does not necessarily induce an alternating diagram  on the boundary.

The following Lemma provides an equivalent   inductive definition.

 \bl  Under the same assumptions as in Definition \ref{DEF2.1}, 
  an ${\cal H}$  is   an  {\it alternating diagram}  if 
 intersecting the collection with any of its hypersurfaces $H$ 
we get an alternating diagram in $H$.
\el 
    
We define colored and mixed domains for alternating diagrams the same way as for ${\cal Q}-$diagrams. Lemma \ref{L2.3} holds for alternating diagrams, see Figure \ref{ncls104ax}:

\begin{figure}[ht]
\centerline{\epsfbox{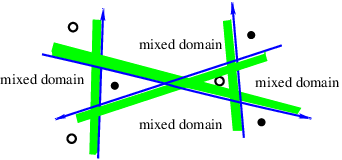}}
\caption{(Co)oriented lines and $\bullet$, $\circ$ and mixed domains for an alternating diagram.}
\label{ncls104ax}
\end{figure}

\subsubsection{${\cal Q}-$diagrams $\lra$  alternating diagrams.}  \la{2.0.3} A ${\cal Q}-$diagram gives rise to an alternating diagram: we  resolve each isolated intersection point $q$ 
  by moving any hypersurface containing $q$ in the direction 
  of its coorientation.    This way we get an equivalence
\be \la{1b}
\begin{split}
& \{\mbox{${\cal Q}-$diagrams}\}  \longleftrightarrow \{\mbox{Alternating diagrams  whose $\bullet-$domains are simplices}\}. \\
\end{split}
\ee
 Indeed,  shrinking the $\bullet-$domains to points we get a ${\cal Q}-$diagram.  
 
 \bl  
Let $q$ be an   isolated intersection point   of a 
${\cal Q}-$diagram   ${\cal Q}$. Then the 
  hypersurfaces of ${\cal Q}$  induce  an alternating collection on the boundary of a little sphere $S_q$ around $q$. 
\el

\subsubsection{${\cal Q}-$diagrams, alternating diagrams,   and bipartite  graphs on surfaces}
 
Here are the basic facts concerning this,  following \cite[Sections 2.2, 2.3, 2.6]{GKe}. 
\vskip 1mm
 
1. By  definition, ${\cal Q}-$diagrams on a surface $S$ are the same as   {\it   triple crossing diagrams} on $S$ \cite{Th}.\footnote{D. Thurston considered triple crossing diagrams only in a disc.}
\be \la{1aaa}
 \{\mbox{${\cal Q}-$diagrams on a surface $S$}\} =  \{\mbox{triple crossing  diagrams on a surface $S$}\}.
\ee 

\begin{figure}[ht]
\centerline{\epsfbox{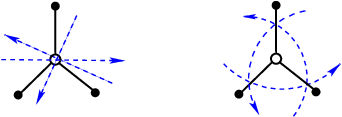}}
\caption{Resolving  triple crossings we get a bipartite graph.}
\label{ncls29}
\end{figure}

2.  A {\it zig-zag strand} on a bipartite surface graph   is a path on the graph  turning left at the $\bullet-$vertices, and right at the $\circ-$vertices. Assigning to a   bipartite  surface graph 
 the collection  of its zig-zag strands pushed out of the vertices   
  provides  an equivalence 
\be \la{1aa}
\begin{split}
 & \{\mbox{minimal bipartite  graphs on a surface $S$}\} =\\
  &  \{\mbox{minimal alternating diagrams on $S$ with contractible colored domains}\}.   \\
   \end{split}
\ee

3.   Any bipartite graph can be transformed to a one with 3-valent $\bullet-$vertices  by splitting the $\bullet-$vertices of valency $>3$ by   adding  edges, and 
 introducing 2-valent    $\circ-$vertices in the middle of these edges. 
 \vskip 1mm
 
 Given a collection  of  zig-zags  strands for a bipartite surface graph $\Gamma$ with 3-valent $\bullet-$vertices, draw zig-zag strands assigned to all edges of $\Gamma$, as on Figure \ref{ncls29}. Then each vertex of $\Gamma$ is surrounded by an oriented polygon of arcs.   Shrink the clockwise triangles near the  $\bullet-$vertices to points, getting a triple crossing diagram.   This way we get an equivalence 
\be \la{1aa}
\begin{split}
 & \{\mbox{minimal bipartite  graphs on a surface $S$ with 3-valent $\bullet-$vertices}\} =\\
  &  \{\mbox{minimal triple crossing diagrams on $S$ with contractible $\circ-$domains}\}.   \\
   \end{split}
\ee
  
  We conclude that 
  
  \begin{itemize}
  
  \item 
 {\it Alternating diagrams on surfaces are slighly more general geometric objects than bipartite  graphs and triple crossing diagrams}.


\item 

{\it ${\cal Q}-$diagrams in threefolds are 3d  analogs of bipartite  graphs / triple crossing diagrams on surfaces}. 

 \end{itemize}
\subsection{Lagrangians assigned to ${\cal Q}-$diagrams in surfaces and threefolds}

We start from  general definitions. 
 Recall that a coorientation of a hypersurface $H$ in a manifold $ X$ is a section of the   bundle of conormal rays $(T^*X- 0)/\R^\times_{>0}$. 
  This section is the Legendrian of the cooriented hypersurface. 
  The {\it conormal bundle} Lagrangian $T^*_HX$  is its preimage   in $T^*X- 0$.
  
    \bd Let ${\cal H}$ be a  finite collection $\{H_i\}$ of smooth cooriented hypersurfaces with disjunct  Legendrians in a manifold $X$. 
  It   gives rise to the following 
  Lagrangian   in $T^*X$:
   \be \la{LL}
{\Bbb L}= {\Bbb L}_{\cal H}:= \mbox{ the union of the zero section $X^\circ$ and the conormal bundles    $T^*_{H_i}X$ to  $H_i$}.
  \ee
  \ed
  
\subsubsection{Lagrangians and spectral surfaces assigned to alternating diagrams on surfaces}   \bd \la{DEF1.4} Let ${\cal T}$ be an alternating diagram on a surface.   Consider   the following   Lagrangian subset:
\be
{\Bbb L}^\circ = {\Bbb L}_{\cal T}^\circ:= {\Bbb L}_{\cal T}- \{\mbox{the union of  zero sections for mixed domains}\}.
  \ee
\ed

\bl \la{L1.5}    The Lagrangian ${\Bbb L}_{\cal T}^\circ$ is homeomorphic to a smooth surface $\Sigma_{\cal T}$, called the {\it spectral surface}:
\be \la{SigmaS}
{\Bbb L}_{\cal T}^\circ = \Sigma_{\cal T}. 
\ee    
The loops of the  diagram $\cal T$  are identified with the boundary circles  on the spectral surface. 
\el

  \begin{proof} Let us describe the spectral surface.   
 Take a disjoint collection of domains ${\cal D}'_*$ matching the colored domains ${\cal D}_*$ of $S -{\cal T}$. 
 For each intersection point $e$ on ${\cal T}$ there are two colored domains 
 ${\cal D}_\bullet$ and ${\cal D}_\circ$ sharing the  $e$. Denote by  ${\cal D}'_\bullet$ and ${\cal D}'_\circ$   disjoint copies of these domains. Glue them along a bridge  replacing the intersection point, with the twist by $180^\circ$,  getting a smooth surface  denoted $\Sigma_{\cal T}$. Then, just as, say,  in  \cite{GKo}, we get a homeomorphism (\ref{SigmaS}).
\end{proof}
   
The Lagrangian  ${\Bbb L}^\circ$ for  a triple crossing diagram   on a surface is Hamiltonian  isotopic to the one for the 
  associated  alternating diagram.   In particular, it is non-singular. 
  \vskip 2mm
  
 {\it Example.} Assume that  colored domains of $S-{\cal T}$ are contractible. 
  Let $\Gamma$ be the bipartite ribbon graph assigned to ${\cal T}$ in (\ref{1aa}).    
  Let  $\Gamma^*$ be the conjugate bipartite ribbon graph,  obtained by reversing the cyclic order of  edges at the $\bullet-$vertices of $\Gamma$ \cite{GKe}, and $S_{\Gamma^*}$ the associated surface glued from the ribbons of   $\Gamma^*$. 
   Then the spectral surface 
  $\Sigma_{\cal T}$ is  the surface $S_{\Gamma^*}$, and  
 \be \la{ESS}
  \{\mbox{zig-zag loops on $\Gamma$}\}=  \{\mbox{boundary loops on $\Gamma^*$}\} = \{\mbox{boundary circles on $\Sigma_{\cal T}$}\}. 
  \ee
   \vskip 1mm

    \subsubsection{Lagrangians and spectral threefolds assigned to ${\cal Q}-$diagrams in  threefolds}
   
          \bd \la{LAGR0} Given a ${\cal Q}-$diagram in a threefold $M$,   consider   the    Lagrangian subset
\be
{\Bbb L}_{\cal Q}^\circ:= {\Bbb L}_{\cal Q}- \{\mbox{the union of  zero sections for mixed domains}\}.
  \ee
\ed  

  The Lagrangian  ${\Bbb L}^\circ_{\cal Q}$ has  singularities   at the quadruple intersection points.       
       Denote by  ${\rm Sing}({\Bbb L}^\circ_{\cal Q})$ the set of  singularities.    The {\it link} of any singular point  of ${\Bbb L}^\circ_{\cal Q}$ is
      the intersection of ${\Bbb L}_{\cal Q}^\circ$ with a small sphere  around the point. 
    Each quadruple intersection point $q$ of ${\cal Q}$ gives rise to a point 
   $q_\circ$ of ${\Bbb L}_{\cal Q}^\circ$ given by the zero cotangent vector at $q$.   
    \bl \la{BST}
The projection $T^*M \lra M$ induces a bijection 
    \be \la{M*}
 \begin{split}
&{\rm Sing}({\Bbb L}^\circ_{\cal Q}) \stackrel{\sim}{\lra}  \{\mbox{the points $q_\circ$ for all quadruple intersection points of ${\cal Q}$}\}.\\  
\end{split}
  \ee 
  The link of a singular point  of ${\Bbb L}_{\cal Q}^\circ$     is homeomorphic to a torus,   equipped with an action of the group $\A_4$.  
    \el
    
   \begin{proof} The complement to  quadruple intersection points  is locally isomorphic  to a transversal intersection of two planes in $\R^3$. 
 So the first claim 
    follows from Lemma \ref{L1.5}.   The second is proved in Lemma \ref{L1.8}. \end{proof}
    
 Let us remove a small ball $\B_q$ around each 
 quadruple intersection points $q$ of ${\cal Q}$. 
 We get a  threefold  $M_\times$ with the induced collection  of surfaces ${\cal Q}_\times\subset M_\times$. 
 It gives rise to a Lagrangian ${\Bbb L}_\times\subset T^*M_\times$, as well as its closed subset ${\Bbb L}_\times^\circ$
 obtained by removing zero sections over all mixed domains in $M_\times-{\cal Q}_\times$.

 \bl \la{L1.12} The Lagrangian 
 ${\Bbb L}^\circ_\times$ is a smooth threefold with  boundary. 
 \el
 
 \begin{proof}  The diagram  ${\cal Q}_\times$ is locally a product of an interval and a cooriented coordinate cross.  So we can apply the 3d version of the construction from the proof of Lemma \ref{L1.5}.   
 Take a disjoint collection of the 3d domains ${\cal D}'_*$ matching the colored domains ${\cal D}_*$ of $M -{\cal Q}$. 
 For each intersection edge $e$ on ${\cal Q}$ there are two colored domains 
 ${\cal D}_\bullet$ and ${\cal D}_\circ$ sharing the edge  $e$. Denote by  ${\cal D}'_\bullet$ and ${\cal D}'_\circ$   their disjoint copies. Glue them with the twist by $180^\circ$,  getting a smooth threefold. 
\end{proof}
 
 \bd The threefold  with boundary ${\Bbb L}^\circ_\times$ is  the spectral threefold $\Sigma_{\cal Q}$ for the ${\cal Q}-$diagram ${\cal Q}$. 
 \ed
 
 The boundary $\partial {\Bbb L}_\times^\circ$  of the Lagrangian threefold ${\Bbb L}_\times^\circ$ 
 assigned to a ${\cal Q}-$diagram of cooriented surfaces $\{S_i\}$  in a threefold $M$ with its own boundary  $\partial M$ plays the crucial role in the story. It is again a Lagrangian in a germ of a symplectic manifold. 
 We give a general definition of the boundary of a Lagrangian with mild singularities in Section \ref{SECTION9}. The boundary  $\partial {\Bbb L}_\times^\circ$  is described in Section \ref{Sect1.3.5}. 
 \vskip 1mm
 
Let us now explain how we associate various stacks to the Lagrangians assigned to ${\cal Q}-$diagrams. 

     \subsection{${\cal Q}-$diagrams  and non-commutative cluster varieties}

A manifold $X$ with boundary $\partial X$ gives rise to the dg-stack of local systems ${\rm Loc}(X)$ and the restriction functor ${\rm Loc}(X) \lra {\rm Loc}(\partial X)$, whose image is a shifted derived Lagrangian. 
One can replace manifolds by singular Lagrangian subsets  with "tame" singularities in symplectic manifolds, and get a similar picture.    
We start from the general set-up, and then proceed to Lagrangians corresponding to ${\cal Q}-$diagrams, which in the case when $X$ is a surface or a threefold provide  
the non-commutative cluster variant of the story. 
    
\subsubsection{Categories and stacks assigned to singular Lagrangians.} \la{1.1}
A closed Lagrangian subset ${\rm L}$ with "tame" singularities in a symplectic manifold  
gives rise to a sheaf of categories $\mathfrak{C}_{{\rm L}}$  on ${\rm  L}$, depending only on the germ of ambient symplectic manifold at ${\rm L}$, see  \cite{K}, \cite{N}, \cite{GPS1}-\cite{GPS3}. 
The section of the sheaf $\mathfrak{C}_{{\rm L}}$  over an open subset $\rm U$ of ${\rm L}$ which can be smoothified to a manifold is equivalent to 
the category of complexes of local systems on ${\rm U}$. 
 Another example: let   ${\cal E}_{n}$ be a single vertex graph with $n$ cyclically ordered edges. Then the section of the sheaf $\mathfrak{C}_{{\rm L}}$ over an open 
subset isomorphic to a product of a disc and the graph ${\cal E}_{n}$   is  the derived category of representations of the quiver of  type $\A_{n-1}$.


There is a Lagrangian $\partial {\rm L}$ of dimension one less than ${\rm L}$, called the {\it boundary of infinity} of ${\rm L}$, or just the {\it boundary of} ${\rm L}$, living in  a germ of a symplectic manifold ${\cal S}$, see Section \ref{SECTION9}, 
  and  the restriction functor:
\be \la{RFC}
{\rm Res}: \mathfrak{C}_{{\rm L}}\lra \mathfrak{C}_{\partial {\rm L}}.
\ee

From now on, ${\rm L}$ is  a closed conical Lagrangian subset  in the cotangent bundle $T^*X$ to a manifold $X$. 
 Then the category $\mathfrak{C}_{{\rm L}}$ can be described in terms of the  
 category $D_{\rm L}^b(X)$ of complexes of constructible sheaves on $X$  with the microlocal support in ${\rm L}$, twisted by a $\Z/2\Z-$gerb with the second Stiefel-Whitney class ${w}_2[T^*X] \in H^2(X, \Z/2\Z)$\footnote{We can twist $\mathfrak{C}_{\rm L}$ by any ${\Bbb G}_m-$gerb. If we use for the twist the pull back of ${w}_2(X)$ by the natural map 
${\rm L}\lra X$, we get the usual category of constructible sheaves.} \cite{KS}, \cite{GKS}.  
This description uses   the microlocalisation functor \cite{KS} 
$$
\mu: D_{\rm L}^b(X) \lra D^b({\rm L}).
$$ 
 Below, given a skew field  $R$,  we consider as the base category the category of complexes of constructible $R-$sheaves on $X$, referred below just as sheaves. 
  We denote by  $\mathfrak{C}_{\rm L}({\rm L})$  the category  of global sections of the sheaf of categories $\mathfrak{C}_{\rm L}$  on ${\rm L}$.   
\vskip 1mm

 Let ${\cal H}$ be a collection  of smooth cooriented hypersurfaces $\{H_i\}$  in $X$. The assigned in (\ref{LL})   Lagrangian ${\Bbb L}= {\Bbb L}_{\cal H}$  in $T^*X$ gives rise to  the 
 subcategory $\mathfrak{C}^{(1)}_{{\Bbb L}}({\Bbb L})\subset \mathfrak{C}_{{\Bbb L}}({\Bbb L})$,    
   given by    the 
 sheaves ${\cal F}$ on $X$ such that: 
\be \la{ASP}
 \begin{split}
& \mbox {\it   
   For each cooriented hypersurface $H_i$ of  ${\cal H}$, the microlocalization  of the sheaf ${\cal F}$}\\
   &\mbox{\it at the punctured conormal bundle $T^*_{H_i }X-\{\mbox{zero section}\}$ is   an $R-$rank one  local system}.\\
   \end{split}
   \ee
   
 \vskip 1mm
 The category $\mathfrak{C}^{(1)}_{{\Bbb L}}({\Bbb L})$ contains the subcategory ${\cal C}_{{\Bbb L}}$  of {\it ${\cal H}-$admissible sheaves} on $X$ \cite{GKo}, see  Section \ref{SEC7.1}. 
 Similar categories for surfaces,  in the commutative set up,   were studied in \cite{STZ}, \cite{STWZ}, \cite{STW}. \vskip 2mm  
  
An  {\it admissible deformation} of the collection ${\cal H}$ is a deformation keeping  the hypersurfaces smooth and their Legendrians disjunct, while  the intersection pattern  can change.   
The category of  ${\cal H}-$admissible sheaves  is invariant under  {admissible deformations}:    this was proved when $X$ is a surface in     \cite{GKo}, and the proof is easily extended to the crucial for us case when $X$ is a threefold.    

\bd Denote by    ${[\cal H]}$ the admissible deformation class of ${\cal H}$. 
We denote by
$
{\cal X}_{\Bbb L} = {\cal X}_{[\cal H]}
$
the derived stack of  objects in ${\cal C}_{{\Bbb L}}$, and by 
${\cal X}_{\partial {\Bbb L}}$ the derived stack of objects in ${\cal C}_{\partial {\Bbb L}}$. 
 \ed
   The restriction functor (\ref{RFC}) induces a map of derived stacks
\be \la{RF}
{\rm Res}: {\cal X}_{\Bbb L}\lra {\cal X}_{\partial  {\Bbb L}}.
\ee

 \subsubsection{Main results.}

 When the manifolds $X$ are surfaces and threefolds with boundary,  the stacks ${\cal X}_{[\cal H]}$ together with  the restriction functors (\ref{RF}) give rise to    
 non-commutative cluster structures of various  types. 
  Their clusters are provided by the ${\cal Q}-$diagrams  in the same admissible deformation class. Precisely:

\begin{itemize}

\item Let $S$ be a  smooth surface with boundary, and 
 ${\cal T}$ an alternating diagram  on $S$. Then the stack ${\cal X}_{[\cal T]}$ carries a non-commutative cluster Poisson structure \cite{GKo}, generalising cluster Poisson varieties assigned to ${\cal T}$ in the commutative setting \cite{GKe}.  The restriction functor (\ref{RF}) describes its symplectic leaves. The clusters are given by the non-commutative tori 
  of flat line bundles on  {spectral surfaces} $\Sigma_{\cal T}$. If the boundary of the surface $S$ is empty, we get a $K_2-$symplectic stack ${\cal X}_{[\cal T]}$. 

\item 
For  threefolds $M$ with boundary, 
   the image of the restriction functor is a 
    non-commutative  cluster Lagrangian ${\cal L}_{[{\cal Q}]}$. 
 The adjective {\it cluster} means that we describe it  in an explicit and uniform way. 
   In the commutative case we get cluster $K_2-$Lagrangians.  
   The clusters are given by the Lagrangians ${\cal L}_{\cal Q}$  assigned to  ${\cal Q}-$diagrams ${\cal Q}$  in the given admissible deformation class $[{\cal Q}]$, which  intersect the boundary 
   by an alternating diagram  ${\cal T}$.  
   Precisely, we define a smooth closed surface $\Upsilon_{\cal Q}$ by  gluing the surfaces of the collection ${\cal Q}$ with the spectral surface 
   $\Sigma_{\cal T}$ of  ${\cal T}$. Then ${\cal L}_{\cal Q}$ is an explicitly described  Lagrangian in the non-commutative torus of flat $R-$line bundles on the surface $\Upsilon_{\cal Q}$:
   \be
   {\cal L}_{\cal Q}\subset   {\rm Loc}_1(\Upsilon_{\cal Q}).
    \ee 
   If all surfaces of the collection ${\cal Q}$ are discs, then $\Upsilon_{\cal Q}$ is the compactified spectral surface ${\bf \Sigma}_{\cal T}$. 
   
   In the commutative case the $K_2-$symplectic structure of the cluster variety ${\cal X}_{[\cal T]}$   provides, by applying the Beilinson-Deligne {\it symbole mod\'er\'e},  a line bundle with connection 
   on the cluster variety, whose  restriction to the cluster Lagrangian ${\cal L}_{\cal Q}$ is canonically trivialised.

Here is an example. Take a threefold $M$ and a set of points $\{p_1, ..., p_n\}$ on its boundary $\partial M$. We construct ${\cal Q}-$diagrams of discs in $M$ describing  the cluster nature of the 
stack of flat $m-$dimensional $R-$vector bundles  on the punctured surface $\partial M- \{p_1, ..., p_n\}$ with the following properties:

 {\it They have  unipotent monodromies around the punctures $p_i$, carry monodromy invariant flags near the punctures,  and  can be extended to $M$}.

\end{itemize}  

  \subsubsection{Why ${\cal Q}-$diagrams, rather than alternating diagrams?} \la{QAL} 

Let  ${\cal H}$ be an alternating diagram of cooriented surfaces obtained by resolution of the quadruple intersection points of a ${\cal Q}-$diagram. Then  we have canonical equivalences of stacks
\be
{\cal X}_{[\cal H]} = {\cal X}_{[\cal Q]}.
\ee 

  The singular set  of the Lagrangian subset ${\Bbb L}_{\cal H}^\circ$ 
consists of the points $t_\circ$ assigned to the triple intersection points $t$ of ${\cal H}$. However the link of  a singular point $t_\circ$ is  a  {\it surface with boundary} - a pair of pants,  
rather than a torus from Lemma \ref{BST}.  So, although we can define the stacks ${\cal X}^\circ_{\cal H}$ for alternating diagrams of surfaces, they are meaningless:  the 
 admissible sheaves vanishing on mixed domains are forced to have the zero microlocalization on the punctured cotangent bundles  to the surfaces.\footnote{The category assigned to a closed Lagrangian subset 
 in the vicinity of  a (single vertex ribbon graph with $n$ edges) $\times \R^m$ is equivalent to the derived category of representations of the quiver $\A_n$, see \cite{K}, \cite{GPS1}-\cite{GPS3}. So near the boundary of a surface we get the category 
 of representations of $\A_0-$quiver, which consists of the single object: the zero.}  
 
 Therefore we do not get a cluster cover of the Lagrangian ${\Bbb L}_{\cal H}$ directly from alternating diagrams.

This is  why we  
deal with ${\cal Q}-$ diagrams of surfaces rather than with the related alternating diagrams.

      \subsubsection{Non-commutative cluster Poisson varieties from alternate diagrams on surfaces.}

 Let ${\cal T}$ be an alternating diagram on a surface $S$.  Recall the Lagrangian ${\Bbb L}_{\cal T}$ assigned to  ${\cal T}$ in 
 (\ref{LL}). The stack ${\cal X}_{[{\cal T}]}$  parametrises  
    admissible dg-sheaves ${\cal F}$ on $S$ with the microlocal support in  ${\Bbb L}_{\cal T}$ satisfying  condition (\ref{ASP}).   
 It carries a non-commutative cluster Poisson structure. 
  Restriction functor (\ref{RF}) describes its symplectic leaves. Let us describe the clusters. 

   \bd \la{DEF1.4a}    Denote by   ${\cal X}^\circ_{\cal T}$  the substack of  stack\ ${\cal X}_{[\cal T]}$ parametrising    objects vanishing on  mixed domains for the alternating diagram   ${\cal T}$. 
So their  microlocal supports are contained in the Lagrangian ${\Bbb L}_{\cal T}^\circ$. 
\ed

Recall the spectral surface $\Sigma_{\cal T}$ for the alternating diagram ${\cal T}$ from Lemma \ref{L1.5}. 
  Let ${\rm Loc}_1(\Sigma_{\cal T})$ be the non-commutative tori  parametrising the $R-$rank one local systems on the  {surface} $\Sigma_{\cal T}$. 
 Thanks to (\ref{SigmaS}),  
\be \la{St1}
  {\cal X}_{{\cal T}}^\circ =   {\rm Loc}_1(\Sigma_{\cal T}).
  \ee
The intersection form on the spectral surface gives rise to a non-commutative Poisson structure on   ${\cal X}_{{\cal T}}^\circ $.   
   The isotopy classes of alternating diagrams ${\cal T}$ in a given admissible deformation class  provide 
  cluster Poisson tori (\ref{St1}),  defining the  non-commutative cluster Poisson structure on the stack ${\cal X}_{[\cal T]}$. 
Modifications   ${\cal T} \lra {\cal T}'$, called  {\it two by two moves},  provide   birational isomorphisms of  non-commutative Poisson tori
 \be \la{2b2Ix}
    {\rm Loc}_1(\Sigma_{\cal T}) \lra   {\rm Loc}_1({\cal T}').
\ee
We conclude that:
\vskip 1mm
 {\it The stack $  {\cal X}_{[{\cal T}]}$ assigned to an alternating diagram ${\cal T}$ on a surface carries  a  structure of a non-commutative cluster Poisson variety. Its cluster tori  are the tori (\ref{St1}) provided by  the alternate diagrams admissibly equivalent to  $\cal T$. The restriction functor (\ref{RF}) describes its symplectic leaves.  The birational maps  (\ref{2b2Ix}) generate the cluster Poisson transformations}.

    \subsubsection{Cluster Lagrangians from ${\cal Q}-$diagrams in threefolds.}   \la{Sect1.3.5}  

  Let $M$ be a threefold, possibly with  boundary, and   ${\cal Q} = \{S_i\}$ a  ${\cal Q}-$diagram   of smooth cooriented surfaces  in $M$  transversal to the boundary $\partial M$. 
  Recall the Lagrangian ${\Bbb L}= {\Bbb L}_{\cal Q}$   assigned   to   ${\cal Q}$ in (\ref{LL}). In Section \ref{SECTION9} we define 
its  {boundary at infinity}
 $$
  \partial {\Bbb L}  \subset {\cal S}.
    $$ 
    It is a singular Lagrangian subset 
$
\partial {\Bbb L} 
$ of the dimension one less than ${\Bbb L}$ 
  in a certain symplectic space ${\cal S}$. 
      It gives rise to a derived symplectic  stack ${\cal X}_{\partial {\Bbb L}}$ 
  and  the {restriction map}
 \be \la{DRM}
 {\rm Res}: {\cal X}_{{\Bbb L}} \lra {\cal X}_{\partial {\Bbb L}}.
  \ee
  The 
    image of the map (\ref{DRM}) is a derived Lagrangian.
  \vskip 1mm
  
 Recall the Lagrangian subset ${\Bbb L}_{\cal Q}^\circ\subset T^*M$ from Definition \ref{LAGR0}. 
  
       \bd \la{LAGR0} Given a ${\cal Q}-$diagram in a threefold $M$, 
  let   ${\cal X}^\circ_{\cal Q}$  the substack of   ${\cal X}_{\Bbb L}$ parametrising   the objects  
with the microlocal supports   in the Lagrangian ${\Bbb L}_{\cal Q}^\circ$. 
\ed  

We denote by $\partial {\Bbb L}^\circ_{\cal Q}$ the boundary of the Lagrangian ${\Bbb L}^\circ_{\cal Q}$, and by ${\cal X}^\circ_{\partial {\Bbb L}}$ the stack of objects corresponding to the Lagrangian 
$\partial {\Bbb L}^\circ_{\cal Q}$. 
  Our next goal is to describe the  Lagrangian $\partial {\Bbb L}^\circ_{\cal Q}$ and the  stack  ${\cal X}^\circ_{\partial {\Bbb L}}$.  
  \vskip 1mm


  Let  us assume from now on  that the  ${\cal Q}-$diagram ${\cal Q}= \{S_i\}$ 
   intersects the boundary $\partial M$ transversally by an alternating  collection   ${\cal T}=\{\alpha_{ij}\}$ of loops in $\partial M$.   
   Then there is  the 
 spectral surface  $\Sigma_{\cal T}$, see (\ref{SigmaS}). It   is a smooth surface, whose 
 boundary is identified with the disjoint union of  loops $\alpha_{ij}$. Indeed,  the loops  of the alternating diagram $\partial {\cal Q}$  become the  boundary circles  on the spectral surface.\footnote{See 
 Lemma \ref{L1.5} and  (\ref{ESS}), describing the case when $\partial {\cal Q}$ comes from a  bipartite graph.}  
  The boundary of each surface $S_i$ is  a disjoint union of the loops $\alpha_{i\ast}$. So: 
 \be \la {uc}
  \partial \Sigma_{\cal T} = \coprod_{j}\alpha_{ij}= \coprod_i  \partial S_i.
 \ee
 
\begin{figure}[ht]
\centerline{\epsfbox{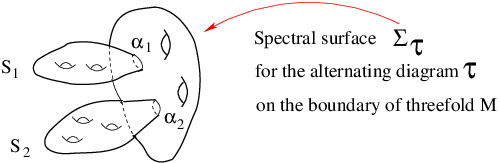}}
\caption{The surface $\Upsilon$ is obtained by gluing surfaces $S_i$  to the spectral surface $\Sigma_{\cal T}$ along  the loops $\alpha_{i*}$.}
\label{ncls26}
\end{figure}

 \bd \la{DEF1.12}The surface $\Upsilon =  \Upsilon_{\cal Q}$ is the smooth closed surface obtained by gluing 
 the disjoint union of  the surfaces $S_i$  with the spectral surface $\Sigma_{\cal T}$ over the matching boundary components, see Figure \ref{ncls26}:
 \be \la{SSS*}
 \Upsilon_{\cal Q}= \Upsilon:= \Bigl(\coprod_iS_i \Bigr)\ast_{\{\alpha_{ij}\}}\Sigma_{\cal T}.
 \ee
\ed
 \bp The   Lagrangian  $\partial {\Bbb L}_{\cal Q}^\circ \subset {\cal S}_{\cal Q}^\circ$  is isomorphic to the zero section Lagrangian $\Upsilon\subset T^*\Upsilon$: 
 $$
 \partial {\Bbb L}_{\cal Q}^\circ  = \Upsilon; \ \ \ \  {\cal S}_{\cal Q}^\circ = T^*\Upsilon.
 $$
 The symplectic  stack $ {\cal X}^\circ_{\partial {\Bbb L}} $ is the non-commutative symplectic  torus of flat line bundles  on the surface $\Upsilon$:
 \be \la{DL0}
 {\cal X}^\circ_{\partial {\Bbb L}} = {\rm Loc}_1(\Upsilon).
  \ee
  \ep
  
    The category ${\cal C}_{\partial {\Bbb L}}$ is the fibered product of the categories $\prod_i {\rm Loc}_1(S_i)$ and  ${\cal C}_{{\cal T}}$ 
   over 
   ${\rm Loc}_1(\{\alpha_{ij}\})$: 
 \be
 {\cal C}_{\partial  {\Bbb  L}} =  \prod_i {\rm Loc}_1(S_i)\times_{{\rm Loc}_1(\{\alpha_{ij}\})} {\cal C}_{{\cal T}}. 
 \ee
 Let us describe  the restriction functor
\be \la{RESF1i}
\begin{split}
&{\rm Res}_{\partial {\Bbb L}}:   {\cal X}_{\Bbb L} \lra    \prod_i {\rm Loc}_1(S_i) \times_{{\rm Loc}_1(\{\alpha_{ij}\})} {\cal X}_{[{\cal T}]}. \\
\end{split}
\ee
  Assigning to an admissible sheaf  the flat $R-$line bundle on $S_i$  given by  the restriction  of its   microlocalisation to $T^*_{S_i}M-\{\mbox{zero section}\}$ we get a functor
$
{\rm Res}_{\infty}:  
 {\cal X}_{\Bbb L} \lra    \prod_i {\rm Loc}_1(S_i). 
$ 
 The restriction to the boundary provides a  functor 
$
 {\rm Res}_{\partial M}:  
 {\cal X}_{\Bbb L} \lra    {\cal X}_{[{\cal T}]}.
$ 
Combining the two functors, we get the  functor (\ref{RESF1i}).
\vskip 2mm

Restricting the functor ${\rm Res}_{\partial M}$ to the  sheaves vanishing on mixed domains we get  a functor 
 \be \la{RESFa} 
\begin{split}
& {\rm Res}^\circ_{\partial M}:  
 {\cal X}^\circ_{{\cal Q}} \lra  {\cal X}^\circ_{{\cal T}} \stackrel{(\ref{St1})}{=} {\rm Loc}_1(\Sigma_{\cal T}).\\
\end{split}
\ee  
So we arrive at the restriction functor 
\be
  \la{RESF1ii}
\begin{split}
{\rm Res}_{\partial {\Bbb L}}^\circ:   {\cal X}^\circ_{\Bbb L} \lra  {\cal X}^\circ_{\partial {\Bbb L}} &
\stackrel{(\ref{DL0})}{=}    {\rm Loc}_1(\Upsilon) \  \\ & \stackrel{(\ref{SSS*})}{=}  \  \prod_i {\rm Loc}_1(S_i) \times_{{\rm Loc}_1(\{\alpha_{ij}\})} {\rm Loc}_1(\Sigma_{\cal T}). \\
\end{split}
\ee
The images  of  stacks  ${\cal X}^\circ_{{\Bbb L} }$ and ${\cal X}_{{\Bbb L} }$ under  restriction functors (\ref{RESF1ii}) and (\ref{RESF1i}) are Lagrangians, denoted by:
 \be \la{14} \begin{split}
& {\cal L}^\circ_{{\cal Q}} := ~{\rm Res}_{\partial {\Bbb L}}^\circ ({\cal X}^\circ_{{\Bbb L}})~\subset   ~{\rm Loc}_1({\Upsilon}).  \\
& {\cal L}_{[{\cal Q}]} :=  ~{\rm Res} _{\partial {\Bbb L}} ({\cal X}_{{\Bbb L}})~\subset  ~\prod_i {\rm Loc}_1(S_i)  \times_{{\rm Loc}_1(\{\alpha_{ij}\})} {\cal X}_{[{\cal T}]}.    \\
\end{split}
\ee
 
  \vskip 1mm
  
    The stacks  ${\cal X}_{{\Bbb L}}$,  $ {\cal X}_{\partial {\Bbb L}}$  and the restriction map do not change under   {admissible deformations}.  However      admissible deformations may not preserve intersection patterns of  surfaces, and hence  ${\cal Q}-$diagrams.  
The  Lagrangians ${\Bbb L}_{\cal Q}^\circ$ and $\partial {\Bbb L}_{\cal Q}^\circ$ and the related stacks ${\cal X}^\circ_{{\Bbb L}}$ and ${\cal X}^\circ_{\partial {\Bbb L}}$
    depend on the choice of a ${\cal Q}-$diagram ${\cal Q}$. 
   Admissible deformations of  ${\cal Q}$ may lead to different isotopy classes of the  alternating diagram ${\cal T}$ on the boundary, altering the spectral surface associated with ${\cal T}$. 
      
   \vskip 1mm
    
    Therefore the ${\cal Q}-$diagrams in $M$  in a given admissible deformation class, considered up to isotopy, give rise to    a cover of the  space ${\cal X}_{\partial {\Bbb L}}$   
  by  Zariski open split tori ${\cal X}^\circ_{\partial {\Bbb L}}$. They are the cluster symplectic tori,  
  defining  a cluster symplectic structure on ${\cal X}_{\partial {\Bbb L}}$.\footnote{Note that since $\partial M$ has no boundary,  the surface $\partial {\Bbb L}^\circ$ are closed, and the related cluster Poisson stacks are  symplectic.} The subvarieties ${\rm Res}_{\partial {\Bbb L}}^\circ({\cal X}^\circ_{{\Bbb L}})$ 
  form a Zariski open cover of the Lagrangian ${\rm Res}_{\partial {\Bbb L}}({\cal X}_{{\Bbb L}})$ in ${\cal X}_{\partial {\Bbb L}}$. 
  We call  the obtained structure  a {\it cluster Lagrangian}.  
   \vskip 1mm
    
\bt \la{MTHGK*} Let ${\cal Q}$ be a ${\cal Q}-$diagram in a threefold $M$ intersecting the boundary transversally by an alternating diagram ${\cal T}$. Then  ${\cal L}_{[{\cal Q}]}$ is a
  non-commutative  cluster Lagrangian. It contains a collection of open Lagrangians 
  ${\cal L}^\circ_{{\cal Q}}$ - the clusters -  parametrised by the ${\cal Q}-$diagrams ${\cal Q}$ in the same admissible deformation class. They are described by systems of equations, depending  uniformly on the cluster. \et

Theorem \ref{MTHGK*} is deduced from the crucial Theorem \ref{T1.14} below. 
   We conclude that: 
   
  \vskip 1mm
{\it ${\cal Q}-$diagrams in threefolds are a geometric source of non-commutative cluster Lagrangians}. 
 
{\it In the commutative setting, they are a geometric source of $K_2-$Lagrangians}.

 \vskip 2mm

  Let us now describe the restriction map (\ref{RESF1ii}) and its image. \vskip 1mm
  
 1. We start with the simplest, and at the same time most basic example, when $M$ is a cube ${\rm C}$. The 1-skeleton of the cube boundary $\partial {\rm C}$  
 has a unique up to a cube rotation structure of a bipartite graph $\Gamma_{\rm C}$, with the  
 $\bullet$ and $\circ$ vertices.  We define the ${\cal Q}-$diagram ${\cal Q}_{\rm C}$ as the collection of four planes passing through the center of the cube, 
 perpendicular to the cube diagonals, and cooriented  $\circ \to \bullet$. 

     \begin{figure}[ht]
\centerline{\epsfbox{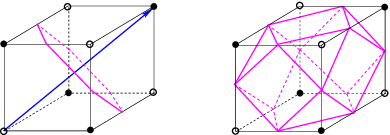}}
\caption{The diagram ${\cal Q}_{\rm C}$ is the union of  four central plane sections perpendicular to the  cube diagonals.} 
\label{AB8*}
\end{figure} 

We give an explicit description of the Lagrangian ${\rm Res}({\cal X}_{{\Bbb L}^\circ_{\rm C}} )$ in Theorem \ref{T1.14}, and an alternative description in Theorem \ref{S111}.
 \vskip 1mm
 
 2. In  general, the key fact is that intersecting  a ${\cal Q}-$diagram with a small ball $\B_q$ around any  
 quadruple intersection points $q$ we get a ${\cal Q}-$diagram isomorphic to  ${\cal Q}_{\rm C}$. 
 Removing all  balls $\B_q$ we get a  threefold  $M_\times$ with the induced collection  of cooriented surfaces ${\cal Q}_\times\subset M_\times$. 
 It gives rise to a Lagrangian ${\Bbb L}_\times\subset T^*M_\times$, as well as its closed subset 
 $
 {\Bbb L}_\times^\circ$ 
 obtained by removing zero sections over all mixed domains in $M_\times-{\cal Q}_\times$.

 \bl \la{L1.20} The Lagrangian 
 ${\Bbb L}^\circ_\times$ is a smooth threefold with  boundary. 
 Its boundary $\partial {\Bbb L}^\circ_\times$ is isomorphic to the disjoint union of the surface $\Upsilon$ and the cube surfaces $\partial {\rm C}_q$ for all 
 quadruple intersection points $q$:
 \be \la{BLY}
 \partial {\Bbb L}^\circ_\times = \Upsilon \cup    \partial {\rm C}_{q_1}\cup ... \cup \partial {\rm C}_{q_m}. \ee
 \el
 
 \begin{proof}  The diagram  ${\cal Q}_\times$ is locally a product of an interval and a cooriented coordinate cross.  So we can apply the 3d version of the construction from the proof of Lemma \ref{L1.5}.   
 Precisely, take a disjoint collection of the 3d domains ${\cal D}'_*$ matching the colored domains ${\cal D}_*$ of $M -{\cal Q}$. 
 For each intersection edge $e$ on ${\cal Q}$ there are two colored domains 
 ${\cal D}_\bullet$ and ${\cal D}_\circ$ sharing the edge  $e$. Denote by  ${\cal D}'_\bullet$ and ${\cal D}'_\circ$   their disjoint copies. Glue them along a bridge  replacing the intersection edge, with the twist by $180^\circ$,  getting a smooth threefold  denoted $\Sigma_{\cal Q}$. Then we get the first claim. The second claim is clear.
\end{proof}
 
 \bd The threefold  with boundary ${\Bbb L}^\circ_\times$ is  the spectral threefold $\Sigma_{\cal Q}$ for the ${\cal Q}-$diagram ${\cal Q}$. 
 \ed
 
Next,  let $N$ be any threefold with boundary $S$. Then, in  the commutative case, the intersection pairing on $S$ provides the canonical element 
 \be
 W_N \in {\cal O}^\times({\rm Loc}_1(S)) \wedge {\cal O}^\times({\rm Loc}_1(S)).  
 \ee
 
 \bl \la{K2}
 The restriction of the element $W_N$ to the subspace ${\rm Loc}_1(S, N)$ of ${\rm Loc}_1(S)$ given by the line bundles on $S$ which extend to $N$ is equal to zero. 
  \el 
  
  \begin{proof} Follows from the fact that the image of the restriction map 
  $
  H^1(N; \Z) \lra H^1(\partial N; \Z)
  $
  is Lagrangian for the $\cup-$pairing on $H^1(S; \Z)$. 
  \end{proof}

  Recall  the map 
  \be
  \begin{split}
  &d\log:  {\cal O}^\times({\rm Loc}_1(S)) \wedge {\cal O}^\times({\rm Loc}_1(S)) \lra \Omega^2_{\rm log}({\rm Loc}_1(S)),\\
   &d\log: f\wedge g  \lra d\log (f) \wedge d\log (g).\\
    \end{split}
\ee 
 The  element $W_N$ provides the $K_2-$symplectic structure $d\log (W_N)$ on ${\rm Loc}_1(S)$. By  Lemma \ref{K2}  the subspace ${\rm Loc}_1(S, N)$  is a $K_2-$Lagrangian.

\begin{figure}[ht]
\centerline{\epsfbox{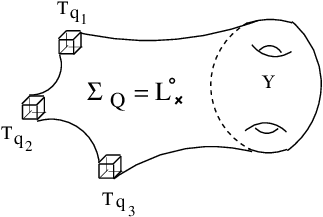}}
\caption{Spectral threefold $\Sigma_{\cal Q}= {\Bbb L}^\circ_\times$. Its boundary is  disjoint union of the closed surface $\Upsilon$, see Figure \ref{ncls26},  and  tori ${\rm T}_{q_i}$. The torus ${\rm T}_{q}$ 
  is the spectral surfaces for the bipartite graph on the cube surface ${\rm C}_{q}$.}
\label{ncls28}
\end{figure}

 Applying this to the spectral threefold $\Sigma_{\cal Q}={\Bbb L}^\circ_\times$, and using (\ref{BLY}),  we get a $K_2-$Lagrangian  
 $$
 \Lambda_{\cal Q} \subset {\rm Loc}_1(\Upsilon) \times \prod_q  {\rm Loc}_1(\partial {\rm C}_q).
  $$
 In Theorem \ref{S111} we describe   the basic $K_2-$Lagrangian ${\cal L}_{\rm C_q} \subset {\rm Loc}_1({\rm C}_q)$. Their product is a $K_2-$Lagrangian  
 $$
{\cal L}_{\rm C}:=  \prod_q  {\cal L}_{{\rm C}_q}\subset  \prod_q  {\rm Loc}_1({\rm C}_q).
   $$
 We consider $\Lambda_{\cal Q}$ as a $K_2-$Lagrangian  correspondence between the $K_2-$symplectic varieties ${\rm Loc}_1(\Upsilon)$ and $\prod_q  {\rm Loc}_1(\partial {\rm C}_q)$. 
 Applying it  to the $K_2-$Lagrangian  ${\cal L}_{\rm C}$ we get the $K_2-$Lagrangian  we were looking for:
 $$
 \Lambda_{\cal Q}\circ  {\cal L}_{\rm C} \subset  {\rm Loc}_1(\Upsilon).
 $$
 By the construction, it comes with an explicit defining system of equations, see Sections \ref{I} $\&$ \ref{II}. 
 \vskip 2mm
 
 The non-commutative version of Lemma \ref{K2} is more complicated.  We will discuss its explicit version elsewhere. It also follows from the general machinery developed in \cite{BrD1}, \cite{BrD2}.

   \subsubsection{The case when all surfaces $S_i$ of the collection ${\cal Q}$ are discs.} 
   
     Filling all holes on the surface $\Sigma_{\cal Q}$ by discs,  we get the compactified spectral surface ${\bf \Sigma}_{\cal T}$. 
     
 If all surfaces of the diagram ${\cal Q}$ are discs, the surface  $\Upsilon_{\cal Q}$ is the compactified  spectral surface ${\bf \Sigma}_{\cal T}$:
\be
\Upsilon = {\bf \Sigma}_{\cal T} \ \ \ \ \ \mbox{\it if all surfaces $S_i$ of the collection ${\cal Q}$ are discs.}\ee
   
  This case is  important for  applications to character varieties.
     Denote by ${\cal U}_{[{\cal T}]}\subset {\cal X}_{[{\cal T}]}$  the substack of    sheafs with trivial microlocalisation at the conormal bundles to  strands of ${\cal T}$. 
   Since any local system on a disc is trivial,     ${\rm Loc}_1(S_i)$ are points,  and 
     the image of the restriction map lands in  ${\cal U}_{[{\cal T}]}$. So 
\be \la{14} \begin{split}
& {\cal L}^\circ_{{\cal Q}}  \subset   {\rm Loc}_1(\bf \Sigma_{\cal T}).    \\
& {\cal L}_{[{\cal Q}]} \subset    {\cal U}_{[{\cal T}]}.    \\
\end{split}
\ee

\bt \la{MTHGK**} Let ${\cal Q}$ be a ${\cal Q}-$diagram of discs in a threefold  inducing   an alternating diagram ${\cal T}$ on the boundary. Then  
  ${\cal L}_{[{\cal Q}]} \subset    {\cal U}_{[{\cal T}]}$ is  a non-commutative cluster Lagrangian.  Its cluster tori are  symplectic tori ${\rm Loc}_1({\bf \Sigma}_{\cal T})$ for the alternating diagrams ${\cal T}$ in the same admissible deformation class. 
\et  
 Theorem \ref{MTHGK**} follows immediately from Theorem \ref{MTHGK*}. See  Section \ref{sec7.3}.

  \subsubsection{Applications to (non-commutative) character varieties.} 
  
   \vskip 1mm
 In Section \ref{Sec7.3} we describe a class of non-commutative cluster Lagrangians in the stack ${\cal U}_{m, S}$ parametrising  $m-$dimensional $R-$local systems on a surface $S$ 
 with punctures, with unipotent monodromies around the punctures. Namely, we show that any threefold $M$  bounding  the surface $S$ with the filled punctures 
gives rise to such  non-commutative cluster Lagrangian ${\cal L}_m(M)$, and  describe it  by a system of equations. 

  The crucial step  is done in Section \ref{sec5.3}, where 
 we  define a ${\cal Q}-$diagram   of cooriented discs in  $M$ starting from a decomposition of $M$ into tetrahedra, and an integer $m\geq 2$ -   the rank of local systems on $\partial M$.  
 Then we use Theorem \ref{MTHGK**}. We describe  the Lagrangian ${\cal L}_m(M)$  using the cluster ${\cal A}-$coordinates assigned in Section \ref{SECTAC} to any ${\cal Q}-$diagram of cooriented discs in a threefold. 
\vskip 1mm
 In the commutative setting the cluster $K_2-$symplectic structure on the space ${\cal U}_{m, S}$ gives rise, via the Beilinson-Deligne {\it symbole mod\'er\'e} \cite{B}, \cite{D}, to a canonical line bundle with connection 
 on the cluster part of ${\cal U}_{m, S}$. Results of Section \ref{Sec7.3} provide a trivialization of  its restriction to the Lagrangian ${\cal L}_m(M)$, depending on a choice of a ${\cal Q}-$diagram in $M$. 
 This generalizes the work \cite{DGG}.
   If $m=2$, this  is closely related to the work of Freed and Neitzke \cite{FN}.

 \subsubsection{Organization of the paper.}

In Section \ref{sectdva} we  describe the  non-commutative Lagrangian for the basic ${\cal Q}-$diagram in the cube. 

In Section \ref{Sec7.2} we prove one of the crucial results, Theorem \ref{T1.14}.

In Section \ref{SECT7} we give another description of the basic non-commutative cluster Lagrangian. We relate it to the ${\cal A}-$coordinate description of the two by two move  \cite{GKo}, uncovering the  surprising $\A_4-$symmetry of the latter.

\vskip 1mm
In Section \ref{SECTION9} we introduce the boundary at infinity of a conical Lagrangian. It  enters to the formulation of  Theorem \ref{MTHGK*}. 

 \vskip 1mm
In Section \ref{SECT3aa} we 
construct spectral covers associated with a certain class (called ideal) of alternating and   ${\cal Q}-$diagrams in threefolds. 
The construction  uses zig-zag surfaces  assigned to 
 these diagrams.

\vskip 1mm
 In Section \ref{SECT6} we discuss non-commutative cluster ${\cal A}-$varieties and their canonical quotients, which carry a cluster symplectic structure, arising from ${\cal Q}-$diagrams in threefolds. 

In Section \ref{SECTAC} we show that a ${\cal Q}-$diagram of cooriented discs in a threefold $M$ gives rise to a cluster Lagrangian in the cluster symplectic quotient of the ${\cal A}-$variety assigned to the boundary of $M$. 

   In Section \ref{sec5.3} we  define a ${\cal Q}-$diagram   of cooriented discs in a threefold $M$ starting from a decomposition of $M$ into tetrahedra, and an integer $m\geq 2$ -   the rank of local systems on $\partial M$.  

In Section \ref{sec7.3} we use  all this to describe the non-commutative cluster Lagrangian ${\cal L}_m(M) \subset {\cal U}_{m, S}$.

\paragraph{Acknowledgment.} The key result of the paper was obtained  during the summer of 2021, 
when the first author held the Gretchen and Barry Mazur Chair at IHES.  The paper  was finished during his  stay  in IHES in 2025. 
The first author was partially supported by the NSF grant DMS-2153059. We are grateful to IHES and NSF for hospitality and  support.

\section{The basic non-commutative cluster Lagrangian}

\subsection{Non-commutative cluster Lagrangian for the basic ${\cal Q}-$diagram} \la{sectdva}

\subsubsection{The basic ${\cal Q}-$diagram.} \la{S1.3.1}

    \begin{figure}[ht]
\centerline{\epsfbox{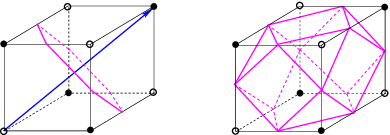}}
\caption{The central plane section perpendicular to a principal diagonal   of the cube cuts   a hexagon. The four planes cut  a polyhedron with 8 triangular and 6 square faces. They cut  the cube  into  four $\circ-$domains, four $\bullet-$domains, and six mixed domains. These domains are cones over the    triangles   assigned to the $\circ/\bullet-$vertices, 
  and the squares. } 
\label{AB8}
\end{figure}

\begin{figure}[ht]
\centerline{\epsfbox{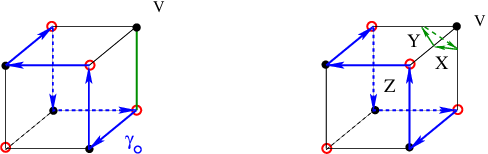}}
\caption{
Let  ${\rm L}_{v}$ be the fiber of a flat bundle ${\rm L}$ on the graph over a $\bullet-$vertex $v$. Then $X,Y,Z\in {\rm End}({\rm L}_{v})$ are  
 monodromies around  faces sharing  $v$  in the direction 
matching the cyclic order of the edges at  $v$. }
\label{3d9}
\end{figure}

Take a cube ${\rm C}$ as our threefold $M$. 
Its 1-skeleton is  a bipartite ribbon graph $\Gamma_{\rm cube}$,  well defined  up to switching $\bullet \leftrightarrow \circ$, see Figure \ref{AB8}.  Consider    
the four planes $H_1, ..., H_4$ passing through the center and perpendicular to the   diagonals of the cube. 
   We  orient  the diagonals by $\circ \to \bullet$. The planes are cooriented by the  diagonals. 
The  four planes form a ${\cal Q}-$diagram ${\cal Q}_{\rm cube}$ with a single quadruple intersection point $p$. It is called the {\it basic ${\cal Q}-$diagram}. 
Each  plane  intersects the cube  by a hexagon. These hexagons subdivide the cube  into six squares and eight triangles. 
  The hexagons can be deformed to a zig-zag loop $\gamma_i$ on the graph $\Gamma_{\rm cube}$, shown by arrows on the left of Figure \ref{3d9}.  There are four zig-zag loops on $\Gamma_{\rm cube}$, 
  which  match the $\circ-$vertices. 
 \vskip 1mm

Pick a  $\bullet-$vertex  $v$.   
 Denote by $X_1, X_2, X_3$ the  monodromies around the  faces sharing the vertex $v$,  
 see Figure \ref{3d9}, understood as endomorphisms of the fiber ${\rm L}_{v}$:
 \be \la{WVE}
 X_1, X_2, X_3\in {\rm End}({\rm L}_{v}).
  \ee 
 
Denote by ${\rm Loc}^{\rm tr}_1(\Gamma_{\rm cube})$ the substack of the non-commutative torus ${\rm Loc}^{\rm tr}(\Gamma_{\rm cube})$ parametrising flat 
rank one bundles on the bipartite ribbon graph $\Gamma_{\rm cube}$ with trivial monodromy around  zig-zag loops. 
 \vskip 1mm
 
The following Theorem \ref{T1.14} is one of the main results of the paper. It is proved in Section \ref{Sec7.2}.

\bt  \la{T1.14}  For the basic ${\cal Q}$-diagram ${\cal Q}_{\rm cube}$,  the restriction functor (\ref{RESF1i}) is  a fully faithful embedding
\be
{\rm Res}^\circ_{\partial {\rm C}}: {\cal X}_{{\cal Q}_{\rm cube}}^\circ \hra  {\rm Loc}^{\rm tr}_1(\Gamma_{\rm cube}).
\ee
Its image ${\cal L}_{{\cal Q}_{\rm cube}}^\circ$ 
 is described by the following equations  on  monodromies  $\{X^{}_1, X^{}_2, X^{}_3\}$ in (\ref{WVE}): 

\begin{itemize}

\item  For each $\bullet-$vertex $v$ of the cube ${\rm C}$ we have: 
\be \la{E2}
X_1 X_2 X_3    +   X_2  X_3  +  X_3 =0.
\ee
\be \la{E1}
 -X_1 X_2 X_3 =1.
\ee
\end{itemize}
The equations (\ref{E2}) at the $\bullet-$vertices  of the cube ${\rm C}$ are conjugate to each other, and thus are equivalent. 
\et

 The  monodromy ${\rm Mon}_{\gamma}$ along the zig-zag loop $\gamma$  opposite to the $\bullet-$vertex $v$ is  is given by:\footnote{The $-$ sign reflects the fact that we work with twisted sheaves.} 
\be \la{E1*}
{\rm Mon}_{\gamma} = -X_1 X_2 X_3 .
\ee
The  equation (\ref{E1}) just means that this monodromy is trivial.

 \bl \la{L1.8} The spectral surface for the bipartite surface graph $\Gamma_{\rm cube}$ is a torus ${\rm T}_4^*$ with four punctures.  The group $\A_4$ acts naturally on this torus. 
 The completion ${\rm T}$ of the torus ${\rm T}_4^*$ is naturally  isomorphic to the link of the Lagrangian ${\Bbb L}_{\rm cube}$ at the singular point. This isomorphism is $\A_4-$equivariant. 
 \el

\begin{proof} It is a connected  oriented surface homotopy equivalent to the graph $\Gamma_{\rm cube}$. Its punctures match the zig-zag loops on the graph $\Gamma_{\rm cube}$. So it must be a torus with four punctures. The group $\A_4$ is realized as the group of symmetries of the cube preserving the bipartite graph. 

The link of the Lagrangian ${\Bbb L}_{\rm cube}$ at the singular point is identified with the cone over 
the two dimensional Lagrangian associated with the zig-zag loops on the bipartite surface graph $\Gamma_{\rm cube}$. \end{proof} 

\vskip 1mm
 
 Let ${\cal F}$ be an admissible sheaf  inside of the cube ${\rm C}$ with the microlocal support at the  Lagrangian ${\Bbb L}_{{\rm cube}}$, which vanishes at the mixed domains. 
 Then its restriction to the complement ${\rm C}-q$ to the quadruple intersection point $q\in {\rm C}$ 
 is equivariant under dilations by the elements of $R^\times$. Therefore according to \cite{GKo} it is described by a flat line bundle on the graph $\Gamma_{\rm cube}$. 
 Precisely, there is canonical equivalence between the category of such 
 admissible sheaves in ${\rm C}-q$ vanishing on the mixed domains,  and the groupoid of flat line bundle on the graph $\Gamma_{\rm cube}$. The latter groupoid is identified by Lemma \ref{L1.8} 
 with the groupoid of flat line bundles on the punctured torus ${\rm T}_4^*$. As we prove in Section  \ref{Sec7.2}, the condition that the intersection of the 
 microlocal support of ${\cal F}$ with $T^*_p{\rm C}$ is supported on the union of the conormals to the planes $H_1, ..., H_4$ is {\it equivalent} to the equations (\ref{E2})-(\ref{E1}) for all $\bullet-$vertices of 
 $\Gamma_{\rm cube}$. In particular the flat line bundle  on the torus ${\rm T}_4^*$ assigned to ${\cal F}$ extends to the punctures. Indeed,   condition (\ref{E1}) just means that its monodromy around each of the punctures is trivial. 
 Then the most interesting condition (\ref{E2}) specifies a Lagrangian in ${\rm Loc}_1({\rm T})$. 
 \vskip 1mm
An important alternative description of this non-commutative Lagrangian is discussed in Section \ref{SECT7}.

 \subsubsection{The $K_2-$Lagrangian for  commutative $R$.} The moduli space of flat line bundles on the graph $\Gamma_{\rm cube}$ is 5-dimensional: the cube has 6 faces, but the product of all   face coordinates is equal to $1$. 
The subvariety    parametrising  flat line bundles  with trivial monodromies around  zig-zag loops can be described by imposing the four monomial equations (\ref{E1}) assigned to the $\bullet-$vertices of the cube. The product of relations (\ref{E1}) over all $\bullet-$vertices  is equal to $1$. 
 The resulting space is a two dimensional symplectic torus. It can be realized as the subtorus  in the torus $(\C^\times)^3$  with  coordinates $X_1, X_2, X_3$, subject to the relation $X_1X_2X_3=-1$. 
 Then the symplectic form is  
 $d\log(X_2) \wedge d\log(X_3)$. 
 Its subspace  ${\cal L}_{{\cal Q}_{\rm cube}}^\circ$ is one dimensional, and hence automatically Lagrangian. However it is in fact $K_2-$Lagrangian, which is a much stronger condition. It amounts to the fact 
 that the equation $X_2X_3+X_3 =1$ implies that\footnote{Recall that the Steinberg relation $\{x, 1-x\}$ in $K_2$ implies that $\{z, -z\}=0$, but $\{z, -1\}$ is only a 2-torsion.} 
 $$
 \{X_2, X_3\} = -\{X_2, 1+X_2\}  =0 \ \ \ \mbox{in $K_2$ mod 2-torsion}.
 $$

   \begin{figure}[ht]
\centerline{\epsfbox{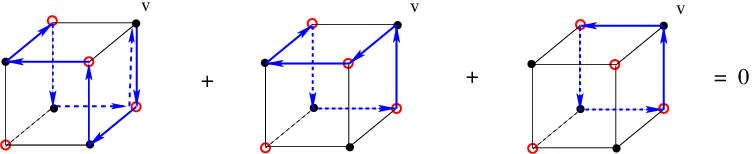}}
\caption{The relation $X_1X_2X_3+X_2X_3+X_3=0$. Together with the condition $X_1X_2X_3=-1$, it reads $X_2X_3 +X_3=1$. Here $X_i$ are the monodromies around the faces 
starting at the vertex $v$, as on Figure \ref{3d9}. }
\label{3d8}
\end{figure}

\subsubsection{ Describing the stack ${\cal X}_{{\cal Q}}^\circ$. } \la{I}
 Given any ${\cal Q}$-diagram ${\cal Q}$ in a threefold $M$,  let 
 \be \la{M*}
 \begin{split}
&M_\times= M - \{\mbox{small open balls in $M$ around the quadruple intersection points}\}.\\
&{\cal Q}_\times:= {\cal Q} \cap M_\times.\\  
\end{split}
  \ee

    \begin{figure}[ht]
\centerline{\epsfbox{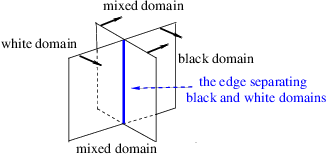}}
\caption{Intersection of two cooriented planes, and the edge separating the black and white domains.}
\label{3d11}
\end{figure}

 \bd  Given a ${\cal Q}$-diagram ${\cal Q}\subset M$, we define  a bipartite graph $\Gamma_{\cal Q}$. Its black (respectively white) vertices are  the $\bullet-$domains (respectively the $\circ-$domains) of ${\cal Q}$. The edges between the vertices   ${\cal D}_\bullet$ and   ${\cal D}_\circ$ 
   correspond to the  edges  of ${\cal Q}$ separating the two domains, see Figure \ref{3d11}. 
   \ed
   The edges of the graph $\Gamma_{\cal Q}$ which are  incident to the vertex assigned to a  colored polyhedron ${\cal D}$ match the edges of the polyhedron. 

          \bt The stack ${\cal X}_{{\cal Q}_\times}^\circ$ is equivalent to the stack of rank one local systems    on the graph $\Gamma_{\cal Q}$. 
   \et
   
   \begin{proof} The stack ${\cal X}_{{\cal Q}_\times}^\circ$ is described the same way as for   alternating diagrams on surfaces in \cite{GKo}.    
   \end{proof}
   
   So since $\pi_1(\Gamma_{\cal Q})$ is a free group, the stack ${\cal X}_{{\cal Q}_\times}$ is  a non-commutative torus.   
   \vskip 1mm
   
Let 
  \be \la{M**}
{\Bbb L}_\times:= {\Bbb L}- \{\mbox{small open balls in $T^*M$ around the singular points of ${\Bbb L}$}\}.
  \ee
  
  
  Let $q$ be a   quadruple intersection point of ${\cal Q}$, and  $S^2_q$ a little sphere  around  $q$. Then ${\cal Q}_\times\cap S^2_q$ is an alternating diagram on the sphere $S^2_q$. 
The related bipartite surface graph   is isomorphic to  $\Gamma_{\rm cube}$. 
The restriction to $S^2_q$ of any ${\cal F} \in {\cal X}_{{\cal Q}}^\circ $ 
is described by 
 {a flat line bundle ${\rm L}_{\cal F}$ on the graph $\Gamma_{{\rm cube}}$ for which  equations (\ref{E2}) from   Theorem \ref{T1.14} hold}.  
\vskip 1mm

The boundary  of the Lagrangian ${\Bbb L}_\times$ is the disjoint  union of the  surface $\Upsilon$ in (\ref{SSS*}),  
realised as the zero section in $T^*\Upsilon$, 
    and the tori ${\rm T}_q$ given by the links of the singular points $q$ of ${\Bbb L}$:
   \be
   \partial {\Bbb L}_\times =  \Upsilon ~\bigcup ~  \bigcup_{q\in {\rm Sing}({\Bbb L}) }  {\rm T}_q.
   \ee
 By Lemma \ref{L1.8},  the torus $\rm T_q$  is  the spectral surface of the bipartite surface graph $\Gamma_{\rm cube}$,  with filled punctures. 
 So  the space ${\rm Loc}_1({\rm T}_q)$ carries a natural symplectic structure.

   Therefore the functor of the restriction to the boundary  Lagrangian $\partial {\Bbb L}_\times$  provides a Lagrangian 
   \be \la{LLLL}
   \Lambda_{\cal Q}:= {\rm Res}_{\partial {\Bbb L}_\times} ({\cal X}_{{\cal Q}_\times}^\circ) \subset  {\rm Loc}_1(\Upsilon)  \times \prod_{q\in {\rm Sing}({\Bbb L}) }  {\rm Loc}_1({\rm T}_q). 
    \ee
    We view it as a Lagrangian correspondence. So it acts on the product of basic Lagrangians in ${\rm Loc}_1({\rm T}_q)$, and the image of this action is the Lagrangian $ {\cal L}_{\cal Q}$.  Let us elaborate this. 
   
  \subsubsection{Describing the Lagrangian $ {\cal L}_{\cal Q}$ by equations.}  \la{II}

  Given symplectic manifolds ${\cal S}_1, {\cal S}_2$ and Lagrangians ${\rm L}_1 \subset {\cal S}_1$ and ${\cal L} \subset 
  {\cal S}_1\times {\cal S}_2$, we get a Lagrangian ${\rm L}_2 \subset {\cal S}_2$ by acting by the Lagrangian 
   correspondence ${\cal L}$ on the Lagrangian ${\rm L}_1$ in the Weinstein category of symplectic spaces: 
   $$
   {\rm L}_2:= {\cal L}\circ {\rm L}_1 \subset {\cal S}_2.
   $$
  We apply this  to  Lagrangian  correspondence $\Lambda_{\cal Q}$ in (\ref{LLLL}). It acts on the product of   basic Lagrangians  
  $$
  \prod_{q\in {\rm Sing}({\Bbb L})}{\cal L}_{{\rm cube}}\subset  \prod_{q\in {\rm Sing}({\Bbb L})}{\rm Loc}_1({\rm T}_q)
    $$
  parametrised by the  singular points $q$ 
  of the ${\cal Q}-$diagram ${\cal Q}$, and provides  a Lagrangian in ${\rm Loc}_1(\Upsilon_{\cal Q})$: 
  \be \la{LLLLL}
  \Lambda_{\cal Q}\circ\prod_{q\in {\rm Sing}({\Bbb L}) }  {\cal L}_{{\rm cube}} \subset {\rm Loc}_1(\Upsilon_{\cal Q}).
  \ee
     It is described by the equations (\ref{E2})-(\ref{E1}) parametrised by the singular points of ${\Bbb L}$. 
     The resulting Lagrangian  can obtained explicitly by    the elimination of variables.

 \subsection{Proof of Theorem \ref{T1.14}} \la{Sec7.2}

 \subsubsection{Set-up.} Let ${\cal F}$ be an admissible sheaf on a threefold $M$, with the microlocal support at a ${\cal Q}-$diagram ${\cal Q}$, which  
   vanishes on mixed domains. We give an explicit description of ${\cal F}$ in that case. 
   Assuming that ${\cal F}$ is of microlocal rank one, we get a cluster description of the corresponding moduli space, 
   as a subvariety in a non-commutative torus defined by the equations (\ref{E2})-(\ref{E1}). This completes the proof of  Theorem \ref{T1.14}. \\
   
\begin{figure}[ht]
\centerline{\epsfbox{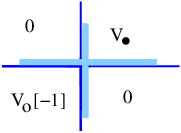}}
\caption{An admissible sheaf  on a surface near an intersection point,  vanishing on mixed domains, with the  microlocal support at the shaded areas and at zero sections over  quadrangles I $\&$ III. The fibers   at  puntcured axes  vanish. Objects $V,W$ are in the degree $0$. The structure  map $h: V_\circ[-1]\to V_\bullet[-1]$ is an isomorphism.}
\label{ncls102a+}
\end{figure}

   Let us first describe ${\cal F}$ on the complement to  quadruple intersection points. In this case the description is reduced to the 2d picture studied in  \cite[Section 10.1]{GKo}. 
   The restriction of ${\cal F}$ to the interior of a $\bullet-$domain $D_\bullet$  is a  constant sheaf  $V_{D_\bullet}$ in the degree $0$. The restriction of ${\cal F}$ to the interior of a $\circ-$domain $D_\circ$  is  the constant sheaf  $V_{D_\circ}[-1]$ in the degree $1$. 
Near an intersection of two surfaces from ${\cal Q}$, it's  transversal section is  described  on Figure \ref{ncls102a+}.

  Let us investigate  now  ${\cal F}$ near a quadruple intersection point $q$.  
 Denote by $B$ a small ball around $q$. The restriction of ${\cal F}$ to the punctured ball 
 $B-q$ is determined   by its restriction to the sphere $\partial B$. The latter is  an admissible dg-sheaf  
 ${\cal F}_{\rm C}$ with the microlocal support at  zig-zag strands of the induced bipartite ribbon graph,  identified with   the graph $\Gamma_{\rm C}$ on the cube ${\rm C}$. 
 Here we use a natural topologically identification $\partial B = {\rm C}$.  The dg-sheaf ${\cal F}_{\rm C}$ 
   vanishes on   mixed domains. 
   The moduli space of such rank one admissible dg-sheaves for any bipartite ribbon graph   is  
 a non-commutative Poisson torus   \cite{GKo}. 
\vskip 2mm

The restriction of the microlocal support of  ${\cal F}$ to $T^*_qM- \{0\}$ is a conical sheaf. 
Our   condition on  ${\cal F}$ near the point $q$ is this:

\begin{itemize}

\item 
{\it The restriction of   microlocal support of  ${\cal F}$  to $T^*_qM- \{0\}$ is concentrated on the four directed conormals  $\eta_1, ..., \eta_4$  to the cooriented surfaces  of ${\cal Q}$ intersecting at     the 
quadruple intersection point $q$.}
  
\end{itemize}
 
 The fiber ${\rm ml}({\cal F})_\eta$ of the microlocal support of ${\cal F}$ at a non-zero covector $\eta\in T^*_qM$ is described as follows. Take a small ball ${\cal U}$ containing $q$. Denote by ${\cal U}_\eta$ its little shift in the direction   $\eta$. Then 
 $$
{\rm ml}({\cal F})_\eta:=  {\rm Cone}\Bigl({\cal F}_{\cal U} \stackrel{{\rm res}}{\lra} {\cal F}_{{\cal U}\cap {\cal U}_\eta} \Bigr).
 $$
A priori  we do not know the fiber ${\cal F}_q={\cal F}_{\cal U}$. However by the condition $\bullet$, the cone   vanishes for any $\eta$ other than  the conormals  to   cooriented surfaces intersecting at  $q$. So  complexes ${\cal F}_{{\cal U}\cap {\cal U}_\eta}$, due to  the isomorphism with ${\cal F}_{q}$,  form a trivial local system on 
the punctured conormal sphere 
$$
S_q^2- \{\eta_1, ..., \eta_4\}, \qquad S_q^2:= (T^*_qM-\{0\})/\R_{>0}.
$$

Let us calculate  the local system $ {\cal L}$ on $S^2_q - \{\eta_1, ..., \eta_4\}$ whose fiber at a covector $\eta$ is given  by 
\be
{\cal L}_\eta:= {\rm R}\Gamma ({\cal F}_{{\cal U}\cap {\cal U}_\eta}).
\ee
 Then we impose the condition that the local system $  {\cal L}$ is trivial, and define   ${\cal F}_q$ as  its  fiber  at   generic $\eta$:
\be\la{FAQ}
{\cal F}_q:= {\cal L}_\eta.
\ee
Then  there is a unique df-sheaf ${\cal F}$, whose restriction to $B-q$ was described above,   and the fiber at $q$ given by (\ref{FAQ}). Let us implement this plan. \\

From now on we identify topologically the conormal sphere  $S_q^2$ and the cube ${\rm C}$ with the bipartite graph $\Gamma_{\rm C}$. 
The covectors $\eta_i$ point to the $\circ-$vertices of the cube, see Figure \ref{AB8}.
We need to calculate the monodromy of the local system $ {\cal L}$   when the covector $\eta$ goes around a cube $\circ-$vertex. 

Pictures below depict  zig-zag strands for the graph $\Gamma_{\rm C}$, drown on a relevant part of the cube surface.  

 \subsubsection{Warm up: the local system around $\bullet-$vertex is trivial.} 
\begin{figure}[ht] 
\centerline{\epsfbox{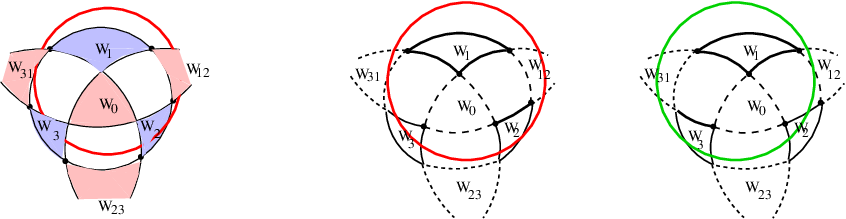}}
\caption{The local system ${\cal L} $ is trivial on the nose.}
\label{ncl2+}
\end{figure}

Consider the local system $  {\cal L}$ near a $\bullet-$vertex $v_\bullet$ of the cube ${\rm C}$  
 on the ray spanned by the covector $-\eta_1$. 
  The seven colored triangles  on Figure \ref{ncl2+} are  intersections of  triangular cones centered at  $q$, 
  containing all vertices of the cube except    $v_\bullet$.  
  
 When the covector $\eta$ is near  $-\eta_1$,  the spherical projection of the domain  ${\cal U}\cap {\cal U}_\eta$ lies inside of a red disc - that is the disc inside of the red circle   on  Figure \ref{ncl2+}.  The disc contains the central red triangle and one of the  blue triangles. Triangles opposite to them do not intersect the  disc. 
 Other four triangles intersect the disc partially. 
 The middle picture on Figure \ref{ncl2+} tells that central triangle  has sides cooriented inside of the triangle, by our convention   on Figure \ref{ncls102a+}. So it is  a $\bullet-$triangle. 
  The  domain on the left of Figure \ref{ncl2+} contains  four red $\bullet-$triangles,  three blue $\circ-$triangles, and three   rectangles. There are six outside crossing points, depicted by $\bullet$ on the left of Figure \ref{ncl2+}.\footnote{The only other crossing points are the three vertices of the red central triangle.} The red disc contains three consecutive       crossing points.   
    
  We claim that   the local system $  {\cal L}$ near the covector $-\eta_1$    is trivial on the nose. 
 Indeed,  the domain ${\cal U}\cap {\cal U}_\eta$ for a   covector $\eta$ near $-\eta_1$ is   the red disc     on Figure \ref{ncl2+}.   Let us calculate the parallel transport of the fiber $ {\cal L}_\eta$  of  $ {\cal L}$. Moving the   $\eta$ around   the covector $-\eta_1$ amounts to moving the  red disc around. There are six different positions of the disc, specified by the consequtive triples of the crossing  points inside the disc. The first one is   the red disc. The next   is   the green disc on the right of Figure \ref{ncl2+},  and so on. \vskip 1mm
 
Figures \ref{ncls102a+} - \ref{ncl2+} tell that the fibers inside of the red triangles  are    $W_0, W_{12}, W_{23}, W_{31}$  in  
 the degree $0$, and inside the  blue ones are $W_1[-1], W_2[-1], W_3[-1]$   in the degree $1$. 
 The fibers at the punctured arcs   are   zero. The fibers at the solid arcs are isomorphic to the fibers inside of the nearby  blue   triangle.  
 For each crossing  separating two colored domains there is an isomorphism between the objects at the domains:
 \be \la{45}
 W_i \stackrel{\sim}{\lra} W_0, \ \ i=1,2,3; \qquad 
 W_{j} \stackrel{\sim}{\lra} W_{ij}, ~~\ W_{i} \stackrel{\sim}{\lra} W_{ij}, \ \ i \not = j.
 \ee
  
  To calculate the complex  ${\rm R}\Gamma({\cal F}_{{\cal U}\cap {\cal U}_\eta})$  we use the following general recipe: 
 \begin{itemize}
 
 \item Given a domain ${\cal U}$   intersecting   strata by contractible subsets, to calculate   $\rm R\Gamma{\cal F}_{{\cal U}}$ we  restrict the data defining the dg-sheaf ${\cal F}$ to ${\cal U}$. Then $\rm R\Gamma{\cal F}_{{\cal U}}$ is the  complex provided by the restricted dg-sheaf data. We need only the restrictions to the strata ${\cal V}$ which lie entirely in ${\cal U}$.  
 \end{itemize}
 
Therefore   ${\rm R}\Gamma({\cal F}_{{\cal U}\cap {\cal U}_\eta})$ for the    red 
 disc on Figure \ref{ncl2+} is   the following complex\footnote{The central open triangle carries the object $W_0$ in the degree $0$. The  sheaf on the closed triangle is  $j_{!}W_0$, where $j$ is the embedding of the open triangle to the closed one. So we get the  cohomology  $W_0[-2]$. Equivalently,    for a Cech cover of the triangle by   open discs ${\cal U}_1, {\cal U}_2, {\cal U}_3$   centered at the vertices,  only the section over ${\cal U}_1 \cap {\cal U}_2\cap {\cal U}_3$ contributes, so we get $W_0[-2]$.} in  degrees $[1,2]$:
  $$
 W_1\oplus W_1\oplus W_3 \lra W_0.
 $$ 
 The similar complex for the   green disc delivers the same result. Since the answer is invariant under the cyclic shift, the other four domains give the same.  
 So the parallel transport from one domain to the other  is   trivial, and the local system ${\cal L}$ is trivial, as expected, since   ${\cal F}$ has no microlocal support  nearby $-\eta_1$.

\subsubsection{The local system ${\cal L}$ near a $\circ-$vertex.}   
\begin{figure}[ht]
\centerline{\epsfbox{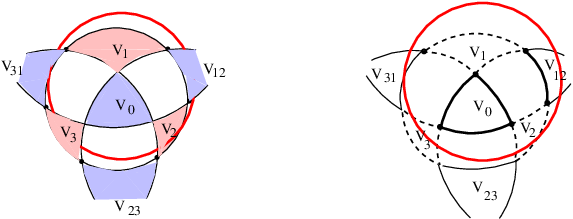}}
\caption{Calculating the $\rm R\Gamma $ for the restriction   of a dg-sheaf ${\cal F}$ to the  red disc.}
\label{ncl1+}
\end{figure}

Let us    calculate  the ${\rm R}\Gamma$ for the restriction of   ${\cal F}$ to the red disc   on  Figure \ref{ncl1+}. 
Figures \ref{ncls102a+} and  \ref{ncl1+} tell that the fibers at   the points of  blue domains  are  the  objects $V_0[-1], V_{12}[-1], V_{23}[-1], V_{31}[-1]$  in  
 the degree $1$, and for the  red  ones are $V_1, V_2, V_3$   in the degree $0$. 
 For each crossing point 
 there is an isomorphism between the objects at the colored domains:
 \be \la{45}
\psi_{0, j}: V_0 \stackrel{\sim}{\lra} V_i, \ \ i=1,2,3; \qquad 
 \varphi_{ij}:  V_{ij} \stackrel{\sim}{\lra} V_j, ~~   \varphi_{ji}: V_{ij} \stackrel{\sim}{\lra} V_i, \ \ i \not = j.
 \ee

Let us introduce shorthands like:
 $$
  V_{0+12+23}:=  V_0 \oplus V_{12} \oplus V_{23} , \quad V_{0+23+31}:=  V_0 \oplus V_{23} \oplus V_{31}, \quad  V_{0+ 31+12}:= V_0 \oplus V_{31} \oplus V_{12}, \quad \mbox{etc.}
  $$ 
The recipe $\bullet$   to the dg-sheaf  for the red disc ${\cal U}$  on Figure \ref{ncl1+} gives the  following complex   in   degrees $[1,2]$: 
 \be \la{sc}
  V_{0+12+31}\to V_1.
    \ee 
Indeed, the domain with the object $V_3$ does not contribute. 
The punctured arc on the boundary of the domain $V_2$, which lies inside of   the red disc,  contributes $0$.  The unique punctured triangle inside of the red disc contributes   $V_0[-2]$.   The blue triangles     contribute  
  $V_{31}[-1]$ and  $V_{12}[-1]$: the first is the section over the vertex, and the second over the unique segment inside of the red disc.

 There   are   quasiisomorphisms given by the   quotients by  acyclic subcomplexes 
  $  V_{13}   \stackrel{ }{\to} V_1$ and  $  V_{12}   \stackrel{ }{\to} V_1$: 
  \be
\begin{split}
&  ( V_{0 + {12} +{13}} \lra V_1 ) \lra   V_{0+12}.\\
&  ( V_{0 +12+13} \lra V_1 ) \lra   V_{0+13}.\\
\end{split}
  \ee 
    \begin{figure}[ht]
\centerline{\epsfbox{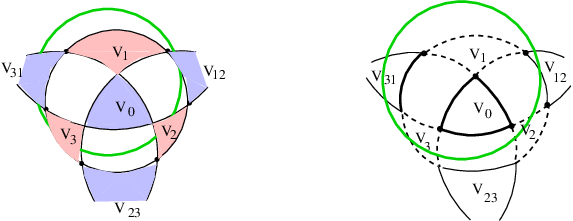}}
\caption{$\rm R\Gamma({\cal F} )$ for   the  green disc   is   identified on the nose with $\rm R\Gamma({\cal F} )$ for the red disc on Figure \ref{ncl1+}.}
\label{ncl1+++}
\end{figure} 
Rotating the red disc, so that it contains  the next triple of crossings,  we get   the green disc on Figure \ref{ncl1+++}. It produces the same on the nose complex (\ref{sc}). 
 So the rotation amounts to the tautological isomorphism.  \vskip 2mm

\begin{figure}[ht]
\centerline{\epsfbox{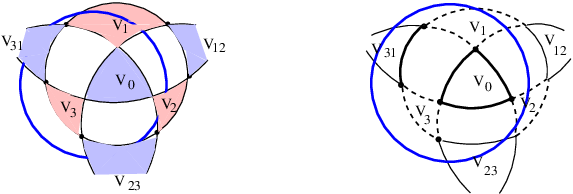}}
\caption{Calculating $R\Gamma({\cal F} )$ for the restriction   of a dg-sheaf ${\cal F}$ to the  blue disc.}
\label{ncl1++}
\end{figure}

Rotating  similarly the green disc we get the blue disc on Figure \ref{ncl1++}, and the   complex  degrees $[1,2]$:
  \be
  V_{0+31 +23 }\to V_3.
    \ee 
  There   are   quasiisomorphisms given by the   quotients by  acyclic subcomplexes 
  $  V_{13}   \stackrel{ }{\lra} V_3$ and  $  V_{23}   \stackrel{ }{\lra} V_3$: 
  \be
\begin{split}
&  ( V_{0 +23+13} \lra V_3 ) \lra   V_{0+23}.\\
&  ( V_{0 +23+13} \lra V_3  ) \lra   V_{0+13}.\\
\end{split}
  \ee    

\vskip 2mm

We repeat the disc rotations  till we come back to the  red disc on Figure \ref{ncl1+}. 
Then  we get a diagram of quasiisomorphisms:    
\be \la{MG3}
\begin{gathered}
    \xymatrix{
    &     
    V_{0+12+23}\to V_2 \ar[ld] \ar[rd]&  &  \ar[rd]  
    V_{0+23+31}\to V_3 \ar[ld]  &&  \ar[rd] 
      V_{0+31+12}\to V_1     \ar[ld]     &  \\
V_{0+12}      &&         V_{0+23}      &&    V_{0+31}      &&      V_{0+12}  \\
        }
\end{gathered}
 \ee
  
Let us calculate the composition of the  quasiisomorphisms between the bottom objects in (\ref{MG3}):
 \be \la{threemaps} 
  V_{0} \oplus V_{12}      \lra           V_{0} \oplus V_{23}      \lra    V_{0} \oplus V_{31}      \lra        V_{0} \oplus V_{12}.
\ee

Each of the  maps in (\ref{threemaps}) is given by the $2 \times 2$ lower triangular matrix.   Let us calculate their product. 
\begin{enumerate}

\item  The  map $V_0 \lra V_{0}$ is the identity map  since it is the composition of the  identical maps: 
$$
V_0 \lra V_0 \lra V_0\lra V_0.
$$ 

\item The   map $V_{12} \lra V_{12}$ amounts to the composition of the  maps 
$$
V_{12}  \lra V_{23}   \lra V_{31} \lra V_{12} .
$$
The latter is  the composition of the  six isomorphisms  $V_{ij} \lra V_i$ and their inverces  in (\ref{45}):

\be \la{MG3bb}
\begin{gathered}
    \xymatrix{
    &    V_{12} & \ar[l]      V_2   \ar[r] & V_{23} & \ar[l]  V_3 \ar[r] & V_{ 31}         &  \ar[l] V_{1} \ar[r]  
          &   
V_{ 12}.                    \\
        }
\end{gathered}
 \ee

\item  The   map $V_0 \lra V_{12}$ is given by the sum of the following three terms, defined below:
\be \la{Sumr}
(\ref{MG3c}) + (\ref{MG3b}) + (\ref{MG3e}).
 \ee
 \end{enumerate}

Each of them is   obtained by composing  the following "up and down" quasiisomorphisms:
  
  \be \la{MG3c}
\begin{gathered}
    \xymatrix{
    &     
    V_{2}    \ar[rd] &  &  
    V_{ 3 }   \ar[rd]   &&  
      V_{1}    \ar[rd]    &  &  \\
V_{0 }   \ar[ur]    &  {-}&        V_{23}    \ar[ur]     & - &     V_{31}    \ar[ur]   & - &      V_{12}  \\ }
\end{gathered}
 \ee
 
  \be \la{MG3b}
\begin{gathered}
    \xymatrix{
    &     
    V_{3}    \ar[rd] &  &  
    V_{1}   \ar[rd]    
     \\
V_{0 }   \ar[ur]    &  {-}&        V_{31}    \ar[ur]     & - &     V_{12}        \\ }
\end{gathered}
 \ee
 
  \be \la{MG3e}
\begin{gathered}
    \xymatrix{
    &      
    V_{1}    \ar[rd] &  &       
     \\
V_{0 }   \ar[ur]    &  {-}&        V_{12}        \\ }
\end{gathered}
 \ee
Let us elaborate the isomorphism $V_0\lra V_{23}$ from the left triangle in (\ref{MG3c}). 
 \be \la{MG3c1}
\begin{gathered}
    \xymatrix{
    &      
   V_{0 +12 +23}  \lra V_{2}  \ar[rd] &  &       
     \\
V_{0 +12}   \ar[ur]    &  {-}&     V_0\oplus    V_{23}        \\ }
\end{gathered}
\ee

Let $(a,b) \in V_{0+12}$. It lifts to a cycle 
 $(a,b,b') \to 0$ in the complex $V_{0+12+23}  \lra V_{2}$, where 
  $
 \psi_{0,2}(a) + \varphi_{12}(b) + \varphi_{32}(b')=0.
  $ 
 Then it projects to the element 
 $$
b'= -(\varphi_{32}^{-1}\circ \psi_{0,2}(a) +\varphi_{32}^{-1}\circ\varphi_{12}(b))  \in   V_{23}.
 $$
 We stress the $-$ sign in this formula by the sign $-$ in (\ref{MG3c}) - (\ref{MG3e}). 
\subsubsection{The meaning of conditions 2) and 3).}  
 
An admissible dg-sheaf near a quadruple intersection point  which vanishes on  mixed domains  induces a similar dg-sheaf on the surface of the cube ${\rm C}$. Since the latter vanishes on mixed domains, it is described by a rank one local system ${\cal L}_{\rm C}$ on the   graph $\Gamma_{\rm C}$. 
  
  White vertices $w$ of the graph $\Gamma_{\rm C}$ 
 match  oriented diagonals $bw$ of the cube, and hence the conormals $\eta_w \in \{\eta_1, ..., \eta_4\}$, and the zig-zag loops on $\Gamma_{\rm C}$.  We denote by $\gamma_w$ the zig-zag loop corresponding to $w$. 
 
\begin{figure}[ht]
\centerline{\epsfbox{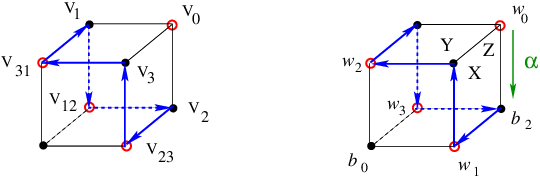}}
\caption{The zig-zag loop $\gamma_{w_0}$ on the bipartite graph $\Gamma_{\rm C}$. The monodromies $X_{w_0},Y_{w_0},Z_{w_0}$ along the loops based at the vertex $w_0$, around the faces sharing $w_0$. The   faces are  ordered clockwise.}
\label{3d10}
\end{figure}

 The composition (\ref{MG3bb}) is   the monodromy of   ${\cal L}_{\rm C}$  along the zig-zag loop $\gamma_{w_0}$ on $\Gamma_{\rm C}$, see   Figure \ref{3d10}. 
 
 Denote by $X_w,Y_w,Z_w$ the   monodromies around the faces sharing a $\circ-$vertex $w$,    over the loops based at  $w$. The  order of the faces  and the loop orientations are compatible with the cyclic order  induced by the ribbon structure of  $\Gamma_{\rm C}$. It is the clockwise order on Figure \ref{3d10},  looking at the cube from the  outside.

 Denote by $\alpha_2  \alpha_1$ the composition of a path  $\alpha_1$ followed by $ \alpha_2$. Let $\alpha$ be the edge path  $w_0 \to b_2$  on   Figure \ref{3d10}.   Then  the   monodromy $  {\rm Mon}_{\gamma, b_2}$   along the zig-zag loop $\gamma$ starting at the vertex $b_2$ is:
 \be \la{MM}
 \alpha^{-1}  {\rm Mon}_{\gamma, b_2}
 \alpha =   Z_{w_0}Y_{w_0}X_{w_0}.  \ee
 
 There are obvious equalities, where on the right stands a map $V_0 \lra V_0$ given by the composition: 
 \be \la{VVV}
 \begin{split}
 &-X_{w_0}  = V_0 \lra V_2 \lra V_{23} \lra V_3 \lra V_0,\\
 &-Y_{w_0}  = V_0 \lra V_3 \lra V_{31} \lra V_1 \lra V_0,\\
 &-Z_{w_0}  = V_0 \lra V_1 \lra V_{12} \lra V_2 \lra V_0.\\
  \end{split}
 \ee

The  $-$ sign  on the left appears since we consider   {\it twisted} dg-sheaves. So making a circle, we rotate the tangent vector by $2\pi$, thus getting an extra $-1$ factor, while on the right stands the plain monodromy.
 
 \bl  For each white vertex $w$ of the graph $\Gamma_{\rm C}$, the condition that   composition (\ref{MG3bb}) is the identity map for   the conormal $\eta_w$ just means that the   local system ${\cal L}_{\rm C}$ on the graph $\Gamma_{\rm C}$ has 
 trivial monodromy along the zig-zag loop $\gamma_w$. Equivalently,  we have the monomial relation 
 \be \la{ZYX}
 Z_wY_wX_w=-1.
 \ee
It is equivalent to the  monomial relation (\ref{MR92}).  
  \el

 \begin{proof} The first and second claims follows from (\ref{MM}) and (\ref{VVV}). The last claim ...
 \end{proof}

\bl For each white vertex $w$ of  the bipartite graph $\Gamma_{\rm C}$, the condition that      the sum (\ref{Sumr}) is zero for the conormal  $\eta_w$  is equivalent  
to  the following relation: 
\be \la{ZY}
  Z_wY_w +Z_w=1.
 \ee 
 \el

 \begin{proof} Composing each of the maps (\ref{MG3c}) -  (\ref{MG3e}) with the map $V_{12}\lra V_0$ we get   (\ref{ZY}) by (\ref{VVV}) and (\ref{ZYX}).  
 \end{proof}

\subsubsection{Conclusion.} Let $q$ be a quadruple intersection point. Let ${\cal F}^\circ$ be an admissible dg-sheaf on the punctured ball $B-q$ vanishing on  mixed domains. 
Then the condition that the composition (\ref{threemaps}) is the identity map is equivalent to the following two:

(i) There is an admissible dg-sheaf ${\cal F}$ on $U$ whose restriction to $U-q$ is given by ${\cal F}^\circ$, and 

(ii) The microlocal support  of ${\cal F}$ at $q$ is supported at the conormals to   surfaces  of ${\cal Q}$ intersecting at $q$.

\subsection{The basic non-commutative cluster Lagrangian in $\mathcal{A}-$coordinates} \la{SECT7}

In Section \ref{SECT7} we give an alternative description of the basic non-commutative cluster Lagrangian. 
It is in fact a description using the non-commutative ${\cal A}-$coordinates introduced in \cite{GKo} and reviewed in Section \ref{SSEECC5}. 
However   Section \ref{SECT7}  is essentially self-contained, and does not require the reader  to know the definition of non-commutative ${\cal A}-$coordinates.  
In fact it provides the crucial example of the latter. 

The description of the basic non-commutative cluster Lagrangian in Section \ref{SECT7}  relates it to  the non-commutative two by two move in ${\cal A}-$coordinates 
introduced in  \cite{GKo}, and uncovers its surprising $\A_4-$symmetry.

\subsubsection{ An  ${\rm A}_4-$invariant Lagrangian  subvariety ${\Bbb L}$ in  a non-commutative  torus   ${\cal A}_{\rm cube}$.} 
\begin{figure}[ht]
\centerline{\epsfbox{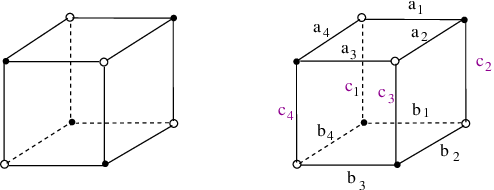}}
\caption{The bipartite ribbon graph $\Gamma_{\rm cube}$ on the cube. On the right: the ${\cal A}-$coordinates on the graph.}
\label{AB11}
\end{figure}

The   
    1-skeleton of the cube $\rm C$  is a bipartite graph $\Gamma_{\rm cube}$, see Figure \ref{AB11}.  It provides   a non-commutative torus 
 ${\cal A}_{\rm cube}$, with the algebra of functions generated  by the variables $(a_i, b_i, c_i)$, $i \in \Z/4\Z$, subject to the  
 monomial relations, corresponding to the vertices of the cube: 
\be \la{S14}
 \begin{split}
 &a_4a_1c_1 =-1, \ \  b_2b_1c_2 =-1, \ \ a_2a_3c_3 =-1,  \ \    b_4b_3c_4 =-1;\\
  &  b_1b_4c_1 =-1,  \ \   a_1a_2c_2 =-1, \ \  b_3b_2c_3 =-1, \ \  a_3a_4c_4 =-1.   \\   
  \end{split}  \ee
Here the first four relations are assigned to  the $\circ-$vertices, and the last four to the $\bullet-$vertices. 
 
Each  monomial   is the   product of the variables at the edges sharing a vertex, in the order compatible with the cyclic order  at the vertex, shown as counterclockwise on the pictures.

 Recall the  system (\ref{S1}) of 12 equations, illustrated on Figure \ref{AB3}:
  \be \la{S1}
 \begin{split}
& a_4b_2+(c_3c_1)^{-1}=1, \ \ \ \  b_2a_4+ (c_4c_2)^{-1}=1, \\
 & a_2b_4+(c_1c_3)^{-1}=1, \ \ \ \  b_4a_2+(c_2c_4)^{-1} =1, \\
&\\
 & a_1b_3+(b_2a_4)^{-1}=1, \ \ \ \  b_3a_1+ (a_2b_4)^{-1}=1, \\
 & a_3b_1+(b_4a_2)^{-1}=1, \ \ \ \  b_1a_3+ (a_4b_2)^{-1}=1, \\
&\\
 & c_1c_3+(b_3a_1)^{-1}=1, \ \ \ \  c_3c_1+(b_1a_3)^{-1} =1, \\
  & c_2c_4+(a_3b_1)^{-1}=1, \ \ \ \  c_4c_2+ (a_1b_3)^{-1}=1. \\
  \end{split}
  \ee

 \bt  \la{S111}
 Equations (\ref{S1})   define  an  ${\rm A}_4-$invariant Lagrangian  subvariety ${\Bbb L}$ in  the non-commutative  torus   ${\cal A}_{\rm cube}$ given  by   equations  (\ref{S14}), with 
   $a_i, b_i, c_i \in R^*$. 
        \et
        
  Theorem \ref{S111} is proved in  Section \ref{S111*}.   
 \begin{figure}[ht]
\centerline{\epsfbox{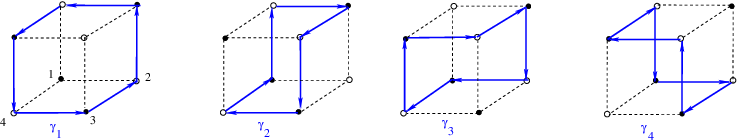}}
\caption{Four zig-zag hexagons $\gamma_1, ..., \gamma_4$ on the cube.}
\label{AB1}
\end{figure}       

We prove in  Lemma \ref{UVIN} below that  equations (\ref{S14})+(\ref{S1}) are equivalent to the equations (\ref{E1})+(\ref{E2}).

\vskip 2mm
The  system of equations   (\ref{S1})  coincides with the one considered in \cite[Section 6]{GKo}. 
It was shown there that  the latter describes 
the  non-commutative cluster Lagrangian provided by the graph of the two by two move on 
non-commutative cluster ${\cal A}-$varieties, described   by formulas (\ref{FRM1}).  
Our description makes evident the  hidden $\A_4-$symmetry of this Lagrangian, which was totally obscure in (\ref{FRM1}). \vskip 1mm

   Using   monomial relations (\ref{S14}), equations (\ref{S1}) can   be written in a polynomial form.\footnote{For example, since by (\ref{S14}) we have $a_4a_1=-c_1^{-1}$ and $a_2a_3=-c_3^{-1}$, we get $c_1^{-1}c_3^{-1} = a_4a_1a_2a_3$. So the equation $a_4b_2+(c_3c_1)^{-1}=1$ from (\ref{S1}) is equivalent to $a_4b_2+a_4a_1a_2a_3=1$ from (\ref{S1a}). Note that the signs in monomial equations (\ref{S14}) can be altered without breacking the above arguments.}
  \be \la{S1a}
 \begin{split}
& a_4b_2+ a_4a_1a_2a_3 =1, \ \ \   b_2a_4+ a_1a_2a_3a_4 =1. \\
 & a_2b_4+ a_2a_3a_4a_1=1, \ \ \   b_4a_2+a_3a_4a_1a_2 =1. \\
&\\
 & a_1b_3+a_1c_1b_1c_2=1, \ \ \ \  b_3a_1+ c_1b_1c_2a_1=1. \\
 & a_3b_1+c_2a_1c_1b_1=1, \ \ \ \  b_1a_3+ b_1c_2a_1c_1=1. \\
&\\
 & c_1c_3+a_2c_2b_2c_3=1, \ \ \ \  c_3c_1+c_3a_2c_2b_2 =1. \\
  & c_2c_4+ c_2b_2c_3a_2=1, \ \ \ \  c_4c_2+ b_2c_3a_2c_2=1. \\
  \end{split}
  \ee

 \subsubsection{Another perspective to equations   (\ref{S1}).}    
\begin{figure}[ht]
\centerline{\epsfbox{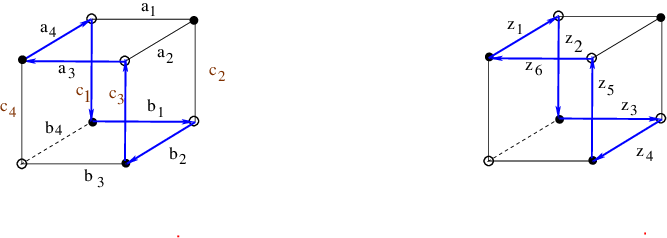}}
\caption{An oriented zig-zag hexagon on the cube.}
\label{AB3}
\end{figure}

There are four zig-zag hexagons on the bipartite graph $\Gamma_{\rm cube}$ on Figure \ref{AB1}. 
They correspond to the principal diagonals of the cube. Namely, the plane perpendicular    to such a diagonal and passing through its center cuts a hexagon on the cube, which is    
the zig-zag loop assigned to this diagonal.   
 Each   zig-zag hexagon gives rise to three equations, described as follows. 
Zig-zags  are oriented so that we turn right at   $\bullet-$vertices, and turn left at  $\circ-$vertices.\footnote{On the pictures  we look to a face from the outside of the cube.}. 
Pick a zig-zag $\gamma$, and a $\bullet-$vertex  on it. Then the edges of $\gamma$   are ordered starting at the edge exiting the $\bullet-$vertex  
and following the orientation of $\gamma$. 
Denote by $z_1, z_2,z_3, z_4, z_5, z_6$ the    elements of $R$ at the  edges, see Figure \ref{AB3}.

\bl  \la{LE3.13} The elements  $z_1, z_2,z_3, z_4, z_5, z_6\in R^\times$ at the   edges of $\gamma$ satisfy  the following equations: 
\be \la{eqq}
z_1z_4+ (z_5z_2)^{-1}=1, \ \ \ \ z_5z_2+ (z_3z_6)^{-1}=1, \ \ \ \ z_3z_6+ (z_1z_4)^{-1}=1.  \ \ \ \ 
\ee
 \be \la{MR92}
z_3z_6 \cdot z_1z_4 \cdot z_5z_2=-1.
\ee 
\el
\noindent
The  relations  are invariant under the cyclic shift by $2$, and so independent  on the choice of a $\bullet-$vertex.  

\begin{proof}  For example, for the  zig-zag $\gamma$ on Figure \ref{AB3},  the   equations (\ref{eqq}) are 
\be \la{RZZ}
 a_4b_2+(c_3c_1)^{-1}=1, \ \ \ \ c_3c_1+ (b_1a_3)^{-1}=1, \ \ \ \ b_1a_3+ ( a_4b_2)^{-1}=1.  \ \ \ \ 
\ee
Using monomial relations (\ref{S14}),  we have 
$
c_1b_1a_3a_4b_2c_3 =   (b_3c_4b_4)^{-1}=(-1)^{-1} =-1. 
$
Equivalently,  
 $b_1a_3a_4b_2c_3c_1 = -1$. This is   equivalent to (\ref{MR92}).  
 \end{proof}
 
\bp \la{P3.11} The monomial equation (\ref{MR92}) follows from the monomial equations $(\ref{S14})$ at $p$. 

 The twelve equations (\ref{eqq}) assigned to the   point $p$ and   four surfaces $H_i$ intersecting at $p$ are equivalent. 
\ep

\begin{proof}  The three equations (\ref{eqq}) are equivalent thanks to the following observation.

 \bl \la{LE3.12} Let us assume that $xyz=-1$. Then the following three equations are equivalent:
$$
x+y^{-1}=1  \ \    \longleftrightarrow  \ \ y+z^{-1}=1  \ \ \longleftrightarrow \ \  z +x^{-1}=1.
$$
\el

\begin{proof} 
Multiplying the first   from the right by $y$ and using $xy=-z^{-1}$ we get the second.  And so on. \end{proof}

By Lemma \ref{LE3.12},  any    equation  in (\ref{S1}) is equivalent to one of the first four: 
 \be \la{S1*}
 \begin{split}
& a_4b_2+(c_3c_1)^{-1}=1, \ \ \ \  b_2a_4+ (c_4c_2)^{-1}=1. \\
 & a_2b_4+(c_1c_3)^{-1}=1, \ \ \ \  b_4a_2+(c_2c_4)^{-1} =1. \\
  \end{split}
  \ee

\bl Each of the   equations (\ref{S1*}) is obtained by  conjugation of one of them. 
  \el
  
  \begin{proof} 
Conjugating the first equation by $b_2$ we get the second equation.   Conjugating the third equation by $b_4$ we get the fourth equation.   
Conjugating the first equation by $c_3$ we get the third equation. \end{proof}

The proof of Proposition \ref{P3.11} follows immediately from these lemmas.  
  \end{proof}

\subsubsection{The rescaling action.} 
Let us label principal diagonals of the cube by the set $\Z/4\Z$, so that the $i-$th diagonal contains the vertex sharing the edges with the coordinates $a_{i-1}, a_{i}$ on Figure \ref{AB3}. 
Denote by $\gamma_i$ the zig-zag  hexagon perpendicular to the $i-$th  diagonal. 
Then an edge $E$ carries the zig-zag assigned to the diagonals generating the plane parallel to the edge $E$. 

Let us rescale each section $s_\gamma$ trivializing the line bundle on the zig-zag $\gamma$ by $s_\gamma \lms \lambda_\gamma s_\gamma$.  

Rescaling the sections $s_{\gamma_i}$ by $\lambda_i$  provides an action of the group $(R^\times)^4$ on the ${\cal A}-$coordinates: 
\be \la{S14a}
\begin{split}
&a_1 \lms \lambda^{-1}_1 a_1 \lambda_2, \ \ a_2 \lms \lambda^{-1}_2 a_2 \lambda_3, \ \ a_3 \lms \lambda^{-1}_3 a_3 \lambda_4, \ \ a_4 \lms \lambda^{-1}_4 a_4 \lambda_1, \\
&b_1 \lms \lambda_4^{-1} b_1 \lambda_3, \ \ b_2 \lms \lambda_1^{-1} b_2 \lambda_4, \ \ b_3 \lms \lambda_2^{-1} b_3 \lambda_1, \ \ b_4 \lms \lambda_3^{-1} b_4 \lambda_2, \\
&c_1 \lms \lambda^{-1}_2 c_1 \lambda_4, \ \ c_2 \lms \lambda^{-1}_3 c_2 \lambda_1, \ \ c_3 \lms \lambda^{-1}_4 c_3 \lambda_2, \ \ c_4 \lms \lambda^{-1}_1 c_4 \lambda_3. \\
\end{split}
\ee

The diagonal subgroup $R^\times$ acts by the conjugation. 

\bp  The action (\ref{S14a}) of the group $(R^\times)^4$ on  ${\cal A}-$coordinates    has the following properties: 

i) Conjugates    equations (\ref{S14}) - (\ref{S1}).

ii) Preserves   the  non-commutative $2-$form $\Omega$. 
\ep

\begin{proof} i) The action of the group $(R^\times)^4$ conjugates   monomial equations (\ref{S14}):
$$
a_ia_{i+1}c_2 \lms \lambda_i^{-1}(a_ia_{i+1}c_2)\lambda_i; \ \ \ \ b_{i+1}b_{i}c_{i+1} \lms \lambda_{i}^{-1}(b_{i+1}b_{i}c_{i+1})\lambda_{i}, \ \ \ \ \forall i \in \Z/4\Z.
$$

Each monomial $z_az_{a+3}$ in   relation (\ref{eqq}) assigned to a zig-zag $\gamma_i$ is conjugated by $\lambda_i^{-1}$. 

So each of the three relations   assigned to a zig-zag $\gamma_i$ is conjugated by $\lambda_i^{-1}$, e.g.
$$
a_4b_2 + (c_3c_1)^{-1}  \lms \lambda_4^{-1}(a_4b_2 + (c_3c_1)^{-1} )\lambda_4.
$$

ii)  Set $\{a,b\}':= da db b^{-1}a^{-1}$. We have 
  $$
\{a_i, a_{i+1}\}' \lms \lambda_i^{-1} \{a_i, a_{i+1}\}'\lambda_i ,\ \ \ \  \{b_{i+1}, b_i\}'   \lms \lambda_{i}^{-1} \{b_{i+1}, b_i\}'    \lambda_{i}. 
 $$
The  2-form  $\{a,b\}$ is the projection of $\{a,b\}'$ to the coinvariants of the cyclic shift. 
 \end{proof}

\subsubsection{Comparing with mutation formulas for the ${\cal A}-$coordinates from \cite{GKo}.}  One can rewrite equations   (\ref{S1}),  using   monomial relations  (\ref{S14})  to eliminate $c_i$'s,  as a birational transformation:  
   \be \la{S3}
 \begin{split}
&(a_1, a_2, a_3, a_4) \lms (b_1, b_2, b_3, b_4), \\
  \end{split}
  \ee                 
   \be \la{S4}
 \begin{split}
& a_4b_2= 1- a_4a_1a_2a_3, \ \ \ \  b_2a_4 = 1-a_1a_2a_3a_4, \\
 & a_2b_4=1- a_2a_3a_4a_1, \ \ \ \  b_4a_2=1-a_3a_4a_1a_2. \\
  \end{split}
  \ee     
Monomial relations (\ref{S14}) imply  
 $
 a_ia_{i+1}= b_{i+1}b_i, \ \ i \in \Z/4\Z.
 $ 
  However  the ${\rm A}_4-$symmetry is  broken.   \vskip 2mm
 
On the other hand, formulas (\ref{FRM1}) for a two by two move, borrowed from \cite{GKo},  
look as follows:    
  \be \la{S24}
 \begin{split}
& a_4\overline b^{-1}_2= 1-a_4a_1a_2a_3, \ \ \ \  \overline b^{-1}_2a_4 = 1-a_1a_2a_3a_4, \\
 & a_2\overline b^{-1}_4=1- a_2a_3a_4a_1, \ \ \ \  \overline b^{-1}_4a_2=1-a_3a_4a_1a_2. \\
  \end{split}
  \ee    
  
\begin{figure}[ht]
\centerline{\epsfbox{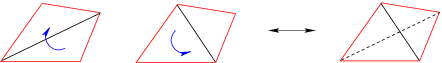}}
\caption{Cutting a tetrahedron into two rectangles  induces the opposite orientations of   rectangles.}
\label{AB4}
\end{figure}

  Formulas (\ref{S24}) transform to formulas (\ref{S4})  via the substitution 
  $$
  \overline b_i:= b_i^{-1}, \ \ i \in \Z/4\Z.
  $$ 
Here is the  explanation. A flip of a triangulation of the rectangle $ABCD$   does not change the rectangle orientation. On the other hand, let us cut a tetrahedron into two 
 rectangles, see Figure \ref{AB4}. We can identify the rectangles  matching the corresponding edges. However then   their orientations    inherited from an orientation of the tetrahedron are opposite to each other. 
 Changing the orientation of  a bipartite graph  means changing the cyclic orders of the edges at vertices to the opposite one. The    ${\cal A}-$coordinates on the edges   change as follows:
\be
a_{E} \lra \overline a_E:= a_E^{-1}. 
\ee

Summarising, we get the following result. 

\bt  Mutation formulas (\ref{S4}) for the non-commutative cluster ${\cal A}-$transformation   corresponding to a two by two move are equivalent to 
the  system of equations (\ref{S1}), where the variables $(a_i, b_i, c_i)$ are related by   monomial relations (\ref{S14}), corresponding to the cube vertices. 
  
        \et

\subsubsection{Geometry of the two  by two move.} 
\begin{figure}[ht]
\centerline{\epsfbox{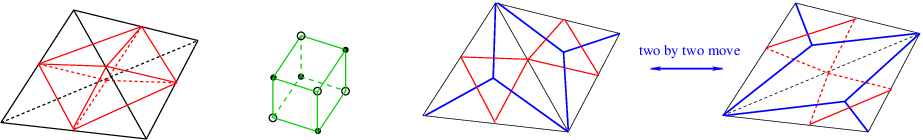}}
\caption{On the left: a black tetrahedron $\&$ the inscribed red octahedron.  In the middle: the green cube, inscribed into the red octahedron, with a bipartite graph. On the right: a    two by two move, blue.} 
\label{AB12}
\end{figure}

Take a tetrahedron ${\rm T}$. Consider the octahedron ${\rm O}$ given by the convex hull of the centers of the edges of ${\rm T}$.  
 Let  ${\rm C}$ be the cube   given by the   
 convex hull of the centers of the faces of the octahedron ${\rm O}$. It carries a  graph $\Gamma_{\rm C}$  
  given by the 1-sceleton of the cube ${\rm C}$.  
 The graph $\Gamma_{\rm C}$ is a bipartite graph: its $\bullet-$vertices   correspond to the  triangles of the octahedra ${\rm O}$ which lie on the faces of the tetrahedra, see Figure \ref{AB12}. 
 
 The bipartite graph $\Gamma_{\rm C}$   on Figure \ref{AB12} is identified with the one on Figure \ref{AB11}. 
The coordinates on the edges of $\Gamma_{\rm C}$ are   identified with coordinates on the edges of the octahedra. 
 We identify them with   coordinates on the blue edges in the two by two move, see Figure \ref{AB12}. So there are bijections: 
$$
 \mbox{$\{$green edges of the cube ${\rm C}$$\}$ = $\{$red edges of the octahedron ${\rm O}$$\}$ = $\{$blue edges of the two by two move$\}$.}
$$
  By the duality between the cube and the octahedron, the ${\cal A}-$variables are assigned to the edges of the octahedron ${\rm O}$, and the monomial 
 relations correspond to the triangular faces of the octahedron.   
 Then  the ${\cal A}-$variables   are assigned to the edges of the blue cube, given by the union of the two blue graphs on the right of Figure \ref{AB12}. The blue cube has two kind of vertices: 
 the centers of the faces of the tetrahedra, and the   vertices of the tetrahedra. So we get  the bipartite graph on the $1-$sceleton of the cube. 
 
\subsubsection{Proof of Theorem \ref{S111}.}  \la{S111*} We already know that ${\Bbb L}$ is isotropic for the canonical   2-form $\Omega$.   On the other hand, take a generic quadruple 
        $({\cal A}, {\cal B}, {\cal C}, {\cal D})$ of decorated flags in a 2-dimensional $R-$vector space. 
        Then a triangulation ${\rm A}{\rm C}$     of the rectangle ${\rm A}{\rm B}{\rm C} {\rm D}$ gives rise to a collection of cluster ${\cal A}-$coordinates 
        on the corresponding bipartite graph $\Gamma_2$. A flip of this triangulation provides a collection of cluster ${\cal A}-$coordinates 
        for the corresponding mutated bipartite graph $\Gamma'_2$.  Note that this is also a corollory of the very general claim in Theorem \ref{MTHGK}.

\subsubsection{The commutative case.}  The commutative torus ${\cal A}_{\rm cube}$   is  5-dimensional.  
Indeed, there are $12$  ${\cal A}-$variables  at 
the  edges   of the cube. There are $8$ monomial relations $r_v$ assigned to the vertices $v$ of the cube, see (\ref{S14}). 
There is a  single relation between the relations:
$$
\prod_{\bullet-\mbox{vertices $b$} }r_b = \prod_{\circ-\mbox{vertices $w$}}r_w.
 $$
   So we get a torus of dimension 
 $12-8+1=5$.   It  describes moduli space of  2-dimensional local systems on the 4-punctured sphere with  a choice of an invariant vector near each puncture, considered modulo 
simultaneous rescaling of all of them.

 The  group $(R^\times)^4$  acts by rescaling of  trivializations of the four zig-zags.  
 The diagonal subgroup $R_{\rm diag}$ acts trivially.  So we get a torus of dimension 
 $12-8+1-3=2$.

  \section{Boundaries at infinity of singular Lagrangians} \la{SECTION9}
  
\subsubsection{The Lagrangian $ {\Bbb L}_{\cal H}\subset T^*X$.}   Let $X$ be a  manifold with boundary $Y$, which carries  a finite collection ${\cal H}$ of smooth cooriented hypersurfaces $\{H_i\}$   with disjoint Legendrians. 
 We assume that the hypersurfaces $H_i$ intersect the boundary transversally by smooth cooriented hypersurfaces $\partial H_i \subset Y$ with disjoint Legendrians. 
 Consider the Lagrangian in $T^*X$ given by the zero section $X^\circ$ and the conormal bundles to  the oriented hypersurfaces:
 \be \la{L1}
 {\Bbb L} = {\Bbb L}_{\cal H}:= X^\circ \cup \coprod_iT^*_{H_i}X.
\ee

  
\subsubsection{The boundary Lagrangian $\partial( {\Bbb L}_{\cal H})$.}     Our goal is define the boundary  of  the Lagrangian $ {\Bbb L} =  {\Bbb L}_{\cal H}$, which is  a Lagrangian  of the dimension one less  in a symplectic space ${\cal S}$ introduced below:
\be
\partial {\Bbb L} \subset {\cal S}.
\ee
It is instructive to give the answer first in the case when the boundary $Y$ is empty. 

\begin{figure}[ht]
\centerline{\epsfbox{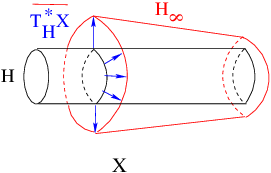}}
\caption{The hypersurface $H_\infty$ is the boundary at infinity   of $T^*_HX$. The Lagrangian $\partial_\infty  {\Bbb L}$ is   the zero section   of $T^*{H_\infty}$.  } 
\label{3dLL}
\end{figure}

Take a sufficiently large compact $C \subset  {\Bbb L}$, so that $ {\Bbb L}-C$ is a union of components  $\R_{[0, \infty)}\times H_i$:
 $$
  { {\Bbb L}-C} = \coprod_i \Bigl({\R_{[0, \infty)}\times H_i}\Bigr). 
  $$
Let  $H_{i, \infty}$ be  the copy of $H_i$, perceived as $\infty \times H_i$.  
Then  $\partial {\Bbb L}$ is  the boundary at infinity $\partial_\infty {\Bbb L}$, defined as  the zero section Lagrangian in  the symplectic space ${\cal S}:=  \coprod  T^*H_{i, \infty}$, see Figure \ref{3dLL}:
\be \la{89}
\begin{split}
 \partial {\Bbb L}:= \partial_{\infty} {\Bbb L} &:= \coprod  H_{i, \infty} ~~ \subset 
~~ \coprod  T^*H_{i, \infty}.\\
\end{split}
\ee

Now let $H$ be a single cooriented hypersurface in $X$ with boundary $\partial H\subset Y$, and $L= T^*_HX$.  

\begin{figure}[ht]
\centerline{\epsfbox{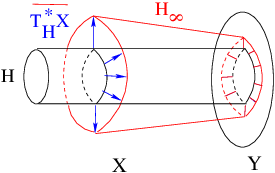}}
\caption{ The Lagrangian $\partial  {\Bbb L}$ is obtained by gluing  the boundary end  $\partial (H_\infty)$ of  $\partial_\infty  {\Bbb L}= H_\infty$ with the boundary at infinity  $(\partial H)_\infty$ of  $T^*_{\partial H}Y$.  } 
\label{3dL}
\end{figure}

Let us define  the following diagram: 

  \be \la{MGIG}
\begin{gathered}
    \xymatrix{
        T^*Y   & T^*H\      
     \\
    T^*_{\partial H}Y \ar[u]_{\alpha_2}^{\cup}&\ar[u]_{\beta_2}^{\cup}H\\
\partial H\times [0,1]  \ar[u]_{\alpha_1}^{\cup} \ar[r]^{=}&     \ar[u]_{\beta_1}^{\cup}  \partial H \times [0,1]   \\ }
\end{gathered}
\ee  
   The map  $\alpha_1$ is the embedding of a collar neighborhood $\partial H\times [0,1] $ of the zero section to $T^*_{\partial H}Y$. The 
   $\alpha_2:T^*_{\partial H}Y \hra T^*Y$ is the natural map.  
 The map $\beta_1$ identifies  $\partial H\times [0,1]$ with a collar neighborhood of  $\partial H$ in $H$. The  $\beta_2$ is  the zero section $H\hra T^*H$.  
 Composing  the vertical maps, we get the diagram
   \be \la{87}
\begin{gathered}
    \xymatrix{
        T^*Y   && T^*H\      
     \\
&\ar[ul]^{\alpha_2\circ \alpha_1}\partial H\times [0,1]   \ar[ur]_{\beta_2\circ \beta_1}&      \\ }
\end{gathered}
\ee    
 We define  the symplectic variety ${\cal S}$    by gluing  symplectic varieties $T^*Y$ and $T^*H$ over  $\partial H \times [0,1]$,  embedded to  $T^*Y$  by the map $\alpha_2\circ \alpha_1$, and to 
 $T^*H$  by  the map $\beta_2\circ \beta_1$, see diagram  (\ref{87}):     \be \la{E11}
 {\cal S}\stackrel{}{:= }T^*Y ~\ast_{\partial H \times [0,1]} ~ T^*H.
 \ee  
 The Lagrangian $\partial  {\Bbb L}$ is obtained by a similar gluing, using the bottom half of  diagram (\ref{MGIG}), see Figure \ref{3dL}: 
  \be \la{E12}
 {\partial  {\Bbb L}}:= \Bigl(Y^\circ \cup T^*_{\partial H}Y\Bigr)  ~\ast_{\partial H \times [0,1]} ~ H.
 \ee 
 The  embedding 
 $
 \partial  {\Bbb L}\subset {\cal S}
 $ 
  is induced from the embeddings $Y^\circ \cup T^*_{\partial H}Y \hra T^*Y$ and $H \hra T^*H$.

  Applying this construction to the union of hypersurfaces $H_i$ we get the symplectic space ${\cal S}$ and the Lagrangian $ {\Bbb L} \subset {\cal S}$ in general. 

\bd Given a finite collection of smooth hypersurfaces $\{H_i\}$ with disjoint Legendrians  in a manifold $X$ with the boundary $Y$, the symplectic variety ${\cal S}$ is defined by gluing $T^*Y$  and $\coprod_i T^*H_i$:
\la{E111}
\be
 {\cal S}:= T^*Y ~\ast_{\coprod_i \partial H_i \times [0,1]} ~ \coprod_i T^*H_i.
 \ee  
 Its Lagrangian $\partial  {\Bbb L}$ is defined by 
\be
 {\partial {\Bbb L}}:=  {\Bbb L}_Y ~\ast_{\coprod_i \partial H_i \times [0,1]} ~ \coprod_i H_i, 
 \ee
where 
$
  {\Bbb L}_Y:= Y^\circ \cup \cup_i T^*_{\partial H_i}Y\subset T^*Y. 
$
 \ed

\subsubsection{The boundary at infinity of Lagrangians.} 
Here is a more general picture. 
Let ${\Bbb L}$ be a possibly singular non-compact closed Lagrangian subset in a symplectic space $M$, which contains a compact  $C$, such that 
   ${\Bbb L}-C$ has the following shape. There exists  a possibly singular Lagrangian ${\Bbb L}'$ of the dimension one less than ${\Bbb L}$ in a symplectic space $U_\varepsilon(L')$, which can be retracted to ${\Bbb L}'$, and a neighborhood $M_\varepsilon$ in $M$, 
   such that 
$$
({\Bbb L}-C \subset M_\varepsilon)= \Bigl({\Bbb L}' \subset U_\varepsilon(L') \Bigr) \times \Bigl(  [0, \infty)\subset ([0, \infty) \times (-\varepsilon, \varepsilon))\Bigr), \ \ \ \ \ \ {\rm dim}({\Bbb L}') = {\rm dim}({\Bbb L})-1.
$$
The second factor is the Lagrangian $[0, \infty)$ in a plane  $(p, q)$,  where $0 \leq p < \infty$ and   $-\varepsilon <q< \varepsilon$, with the symplectic form $dp\wedge dq$, so that
 $$
{\Bbb L}-C = {\Bbb L}'\times[0, \infty).
$$
  Then the boundary at infinity $\partial_\infty{\Bbb L} $ is  the Lagrangian ${\Bbb L}' $ in the symplectic space $U_\varepsilon({\Bbb L}')$:
$$
\partial_\infty {\Bbb L} :={\Bbb L}'\subset U_\varepsilon({\Bbb L}').
$$

In particular, consider a smooth hypersurface $H$  without the boundary in a manifold $X$. It gives rise to 
the  Lagrangian ${\Bbb L}$ given by the union of  $T^*_HX$ and the zero section $X^\circ$. {The Lagrangian ${\Bbb L}$ has singularity  
at  $T^*_HX \cap X^\circ$ of type (a single vertex graph with 3 edges) $\times \R^{n-2}$, where $n = {\rm dim}(X)$.}Then applying the general definition above we arrive at the definition (\ref{89}). 

\vskip 1mm

One can define Lagrangians at infinity with corners. Let us spell the case of codimension one corners, which is relevant to our story. 
In this case we assume that the complement ${\Bbb L}-C$ has two ends. 
Namely, there exists  
a possibly singular Lagrangian ${\Bbb L}''$ of the dimension two less than ${\Bbb L}$ in a symplectic space $U_\varepsilon(L'')$, which can be retracted to ${\Bbb L}'$, and a neighborhood $M_\varepsilon$ in $M$, 
   such that $$
{\Bbb L}-C = \Bigl({\Bbb L}'' \subset U_\varepsilon({\Bbb L}'') \Bigr) \times \Bigl( [0, \infty)^2 \subset ([0, \infty) \times (-\varepsilon, \varepsilon))^2\Bigr), \ \ \ \ \ \ {\rm dim}({\Bbb L}'') = {\rm dim}({\Bbb L})-2.
$$
The second factor is the subspace  of  the symplectic space  $\R^4$ with  coordinates $(p_1, q_1, p_n, q_n)$ with the symplectic form 
$dp_1\wedge dq_1 + dp_n \wedge dq_n$, given by  $0 \leq p_n, q_1 < \infty$ and   $-\varepsilon <q_n, p_1< \varepsilon$. 
It contains the Lagrangian 
$[0, \infty)^2$ with coordinates $(p_n, q_1)$. We compactify it to a square $[0, \infty]^2$. The two sides of the square containing the corner $(\infty, \infty)$ form an angle 
$([0, \infty]\times \{\infty\}) \cup (\{\infty\} \times [0, \infty])$. 
The corner Lagrangian $\partial_{(2)}{\Bbb L}$ is 
the product of the Lagrangian ${\Bbb L}'' $  and the angle:
$$
\partial_{(2)}({\Bbb L}) := {\Bbb L}''\times \Bigl(([0, \infty]\times \{\infty\}) \cup (\{\infty\} \times [0, \infty])\Bigr).
$$
It contains the corner $ {\Bbb L}''\times \{\infty\}^2$. The Lagrangian $\partial_{(2)}({\Bbb L}) $ lies in the symplectic space 
$$
\partial_{(2)}{\cal S} = U_\varepsilon({\Bbb L}'') \times \Bigl(([0, \infty]\times \{\infty\}) \cup (\{\infty\} \times [0, \infty])\Bigr)\times (-\varepsilon, \varepsilon)^2.
$$
       \begin{figure}[ht]
\centerline{\epsfbox{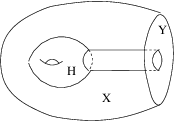}}
\caption{  } 
\label{3dL1}
\end{figure}    
\subsubsection{ The sheaf of categories ${\cal C}_{\Bbb L}$ and its restriction  to the boundary Lagrangian.} 
 Recall   the singular Lagrangian ${\Bbb L}$ in (\ref{L1}), and the  sheaf of categories  ${\cal C}_ {\Bbb L}$ on  $ {\Bbb L}$,  
discussed in Section \ref{1.1}.   Recall  the zero section $H^\circ \subset T^*H$. For open sets  $\R_{(a, \infty]}\times H_i$ there are canonical equivalences
 \be \la{caneq}
 \begin{split}
 &{\cal C}_{\R_{(a, \infty]}\times H_i}  \stackrel{\sim}{\lra} {\cal C}_{ H_i}; \\ 
 & {\cal C}_{ {\Bbb L}-C}  = \prod_i  {\cal C}_{\R_{(a, \infty]}\times H_i} \stackrel{\sim}{\lra} \prod_i{\cal C}_{H^\circ_{i, \infty}}.\\
  \end{split}
  \ee
  So the restriction   functor $
   {\rm Res}_{ {\Bbb L}-C}: {\cal C}_ {\Bbb L} \lra  {\cal C}_{ {\Bbb L}-C}
   $
    provides, using (\ref{caneq}),  the restriction at infinity functor  
 $$
 {\rm Res}_\infty: {\cal C}_ {\Bbb L}  \stackrel{}{\lra} {\cal C}_{{\cal H}^\circ_\infty}:=  \prod_i{\cal C}_{ H^\circ_{i, \infty}}.
 $$
 
 The transversal intersection of the collection of hypersurfaces ${\cal H}$ with the boundary $\partial X$ induces a collection $\partial {\cal H}$ on the boundary. There is the Lagrangian ${\Bbb L}_{\partial {\cal H}}\subset T^*\partial X$, see (\ref{L1}). 
The restriction to the boundary $\partial X$ of $X$ provides the functor
 $$
  {\rm Res}_{\partial X}: {\cal C}_ {\Bbb L}  \stackrel{}{\lra}   {\cal C}_{{\Bbb L}_{\partial  {\cal H}}}. 
  $$
 There are also the restriction functors, where  the bottom line functors are provided by the top line:
 \be
 \begin{split}
 &{\cal C}_{H^\circ_{i}} \lra    {\cal C}_{\partial H^\circ_{i}}, \ \ \ \ \ \    {\cal C}_{{\Bbb L}_{\partial  {\cal H}}} \lra   {\cal C}_{\partial H^\circ_i}.\\
 &{\cal C}_{{\cal H}^\circ_\infty} \lra    {\cal C}_{\partial {\cal H}^\circ_{\infty}}, \ \ \ \ \ \    {\cal C}_{{\Bbb L}_{\partial  {\cal H}}} \lra   {\cal C}_{\partial {\cal H}^\circ_\infty}.\\
 \end{split}
 \ee
  The functors ${\rm Res}_\infty$ and ${\rm Res}_{\partial X}$ gives rise to  the restriction functor to the fibered product of categories:
   \be \la{RESF1}
  {\rm Res}: {\cal C}_ {\Bbb L}  \stackrel{}{\lra}  {\cal C}_{{\cal H}^\circ_\infty} \times_{{\cal C}_{\partial {\cal H}_\infty^\circ}} {\cal C}_{{\Bbb L}_{\partial  {\cal H}}}.
  \ee
 
 The following can  be reduced from the results of Brav-Dyckerhoff \cite{BrD1}-\cite{BrD2}. 
 \bt \la{Th9.2}
 The image of the restriction functor (\ref{RESF1}) is derived Lagrangian. 
 \et

\subsubsection{The case of  threefolds.} Consider  a threefold $M$  with boundary $\partial M$, equipped with a ${\cal Q}-$collection of smooth cooriented surfaces $\{S_i\}$. 
Recall the smooth closed surface $\Upsilon_{\cal Q}$. 

  \bl   
The boundary  $\partial {\Bbb L}^\circ$ of the Lagrangian ${\Bbb L}^\circ$ is the smooth closed surface $\Upsilon_{\cal Q}$:
   $$
  \partial {\Bbb L}^\circ=    \Upsilon_{\cal Q}.
    $$\el
   
  The restriction to the boundary  provides  projections
  \be
  \begin{split}
 & \pi_1: {\rm Loc}_1(\Sigma) \lra {\rm Loc}_1(\partial \Sigma).\\
 & \pi_2: {\rm Loc}_1(S_i) \lra {\rm Loc}_1(\partial S_i).\\
  \end{split}
  \ee
      
   Recall the   restriction functor:
  \be \la{RESF2}
\begin{split}
&  {\rm Res}_{\partial  {\Bbb L}^\circ}: {\cal C}^\circ_ {\Bbb L}  \stackrel{}{\lra}  \prod_i {\rm Loc}_1(S_i) \times_{{\rm Loc}_1(\partial \Sigma)} {\cal C}^\circ_{{\cal T}}. \\
\end{split}
  \ee 
   The category ${\cal C}^\circ_{{\cal T}}$ is equivalent to the groupoid $ {\rm Loc}_1(\Sigma)$. 
 So we arrive at the restriction functor
  \be \la{RESF3}
 \begin{split}
{\rm Res}_{\partial  {\Bbb L}^\circ}: {\cal C}^\circ_{{\Bbb L}} \ \  \stackrel{}{\lra}  \ \ 
  &{\rm Loc}_1({\Upsilon}_\Sigma) \ \ = \ \ \prod_i {\rm Loc}_1(S_i) \times_{{\rm Loc}_1(\partial \Sigma)} {\rm Loc}_1(\Sigma). \\
\end{split}  \ee   

  \bt \la{THE4.4}
 The image of the restriction functor (\ref{RESF3}) is Lagrangian. 
 
   In particular,  if the surfaces $S_i$ are discs,  we get a Lagrangian subvariety 
     \be \la{RESF5}
  {\rm Res}_ {\partial {\Bbb L}^\circ}({\cal C}^\circ_{{\Bbb L}})  \subset  {\rm Loc}_1(\Sigma). 
  \ee   \et
 
 \begin{proof} The first claim   follows from Theorem \ref{Th9.2}.
  If  surfaces $S_i$ are discs, then 
   ${\rm Loc}_1(S_i)$ are  points, and  the   restriction functor (\ref{RESF3}) boils down to
   $
  {\rm Res}_ {\partial {\Bbb L}^\circ}: {\cal C}^\circ_{\Bbb L}  \lra {\rm Loc}_1(\Sigma). 
 $ 
 \end{proof}

 \section{Spectral covers  for  alternating and ${\cal Q}-$diagrams in threefolds} \la{SECT3aa}

 Let $M$ be a threefold. Let ${\cal Q}$ be a ${\cal Q}-$diagram in $M$. Recall the Lagrangian ${\Bbb L}_{\cal Q}^\circ$, and its subset  ${\Bbb L}_\times^\circ$ obtained by cutting out  
  singular points of ${\Bbb L}_{\cal Q}^\circ$, 
 which correspond to  quadruple intersection points 
 of  ${\cal Q}$. 
 Recall the threefold 
 $$
 M_\times := M-\{\mbox{quadruple intersection points}\}.
  $$
   In Section \ref{SECT3aa}, under  mild assumptions on ${\cal Q}$, we  construct a smooth threefold with boundary  $\Sigma_{\cal Q}$.    The threefold $\Sigma_{\cal Q}$ is homeomorphic to 
 ${\Bbb L}_{\times}^\circ$, and  comes with  a finite projection  
  $$
  \pi_{\cal Q}: \Sigma_{\cal Q} \lra M_\times.
  $$
  To achieve this, we construct  the spectral cover for an alternating diagram ${\cal H}$ in  $M$. The construction  uses {\it zig-zag surfaces} assigned to alternating diagrams. 
  We apply   this construction to the alternating diagram obtained  by resolving quadruple intersection points of ${\cal Q}$, and cutting out the obtained $\bullet-$tetrahedra.

 \subsubsection{A singular surface $\bS_{\cal H}$ for an alternating  diagram ${\cal H}$.}     
  
 Let ${\cal H}$ be an alternating diagram  in a threefold $M$. Its singular points are the isolated {triple intersection points}, and   the {\it edges} $E$, 
 given by   components of the singular locus of ${\cal H}$  minus  triple intersection points.

 Let us assume   that the $\bullet-$ and $\circ-$domains are polyhedrons. Let us    
 define  a singular surface 
$$
\bS_{\cal H}\subset X.
$$ 
Let $D$ be  a \underline{colored} polyhedron. Pick an internal point $c_{D }\in D $. The cone  over the $1-$skeleton of   $D $ with the vertex at
   $c_{D}$  
is the {\it singular surface} $\bS_{D } \subset D$.    
 
\bd \la{SSH}  The singular surface $\bS_{\cal H} $ is the union of    surfaces $\bS_{D}$ in all colored polyhedrons:$$
\bS_{\cal H}:= \cup_{D_\bullet} \bS_{D_\bullet} \cup_{D_\circ} \bS_{D_\circ}. 
$$ 
\ed

  Given an edge  $E$ of   ${\cal H}$, there are exactly two colored polyhedrons $D_{E, \bullet}$ and $D_{E, \circ}$ sharing the edge. 
Their colors are different. In each of the  two polyhedrons  $D_{E, *}$ there is a triangle $t_{E, \ast}$   with the base  $E$ and the vertex at   $c_{D_{E,*}}$. 
Denote by $r_E$ the rectangle given by their   union:
$$
r_E = t_{E, \bullet} \cup t_{E, \circ}.
$$

{\it The singular surfaces $\bS_{\cal H}$ are 2d  analogs of bipartite ribbon graphs.}  

{\it The points $c_D$ and rectangles $r_E$ are  the analogs of  vertices and edges of a bipartite graph}.


\subsubsection{Cooriented zig-zag surfaces  for an alternating diagram ${\cal H}$.} 
\begin{figure}[ht]
\centerline{\epsfbox{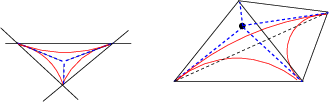}}
\caption{Smoothified strands/surfaces $\gamma_{D, i}$ (red)  in    a triangle/simplex.} 
\label{3d4}
\end{figure}

For each  face $f$ of  a colored  polyhedron $D$, consider the cone $C_{D,f}$  over  $f$ with the vertex $c_D$.  
The polyhedral surface $\gamma_{D, f}$ is  the boundary of   $C_{D,f}$. It is the cone over    $\partial f$. 
A coorientation  of   $\gamma_{D,f}$  is  
a coorientation of the surface obtained by smoothifying  its corners,  see Figure \ref{3d4}. 
We coorient the  surface $\gamma_{D, f}$ inside of  the cone $C_{D,f}$ if   $D$ is $\circ-$polyhedron, and outside  if $D$ is $\bullet-$polyhedron. \vskip 2mm

\begin{figure}[ht]
\centerline{\epsfbox{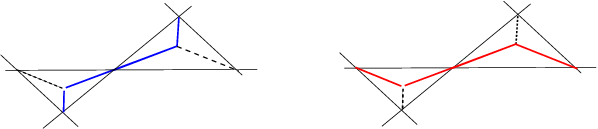}}
\caption{The two zig-zag strands, colored in blue and red, sharing a common edge.} 
\label{3d5}
\end{figure}

For each of the two triangles $t_{E, \ast}$  sharing an $E$ of ${\cal H}$, where $\ast = \circ, \bullet$, there are   two  surfaces $\gamma_{D_*, f_1}$ and $\gamma_{D_*, f_2}$ containing the cone. Here $f_1, f_2$ are the two faces of $ D_{E, *}$ sharing   $E$.   Each of the two surfaces $\gamma_{D_\circ, f}$ can be  connected to just one of the two   surfaces $\gamma_{D_\bullet, f}$ 
so that they form a cooriented disc, see Figure \ref{3d5} for the 2d analog.    
This way each disc extends to the unique  cooriented polygonal surface  $\gamma$, which we call a {\it zig-zag surface}.  \vskip 1mm

{\it Zig-zag surfaces are analogs of   zig-zag strands of a bipartite graph $\Gamma$}.

\subsubsection{The spectral cover.} Assume  that each zig-zag surface $\gamma$ cuts $M$ into two domains. For just one of them, denoted by ${\Bbb D}_\gamma$,      the "inside the domain"  boundary coorientation   match the coorientation of   $\gamma$. Take a copy    ${\cal D}_\gamma$   of the domain ${\Bbb D}_\gamma$.     
Then the disjoint union    of all   domains ${\cal D}_\gamma$ projects tautologicaly onto $M$, so that  each domain $D_\gamma$ maps to its copy ${\Bbb D}_\gamma$ in $M$:
$$
\coprod_{\gamma} {\cal D}_\gamma \lra M.
$$
We glue the  domains ${\cal D}_{\gamma}$ and  ${\cal D}_{\gamma'}$   if their projections   share a rectangle $r_E$, getting the {\it spectral space} $\overline \Sigma_{\cal H}$:
\be \la{SIG}
 \overline  \Sigma_{\cal H}:=  \frac{\coprod_{\gamma} {\cal D}_\gamma}{\mbox{gluing  domains ${\cal D}_{\gamma}$ and  ${\cal D}_{\gamma'}$ along the rectangle $r_E$} }.
\ee
    Then there is the canonical map, which identifies the domain ${\cal D}_\gamma$ to the one ${\Bbb D}_\gamma$:
   \be \la{22Ma}
    \overline \Sigma_{\cal H} \lra M.
   \ee

\bd 
An alternating diagram is called {\rm ideal} if there are special points 
$\{s_i\}$ 
on the boundary  of $M$ such that for each zig-zag surface $\gamma$, the domain ${\Bbb D}_\gamma$ is a (semi-)ball containing just one special point.
\ed

    If ${\rm dim}M=2$, we recover the ideal webs  \cite{G}. 
 Consider 
the {\it spectral threefold}
$$
\Sigma_{\cal H}:= \overline \Sigma_{\cal H}- \{\mbox{singular points}\}.
$$
It  generalises the spectral surface $\Sigma_{\Gamma}$ assigned to a bipartite ribbon graph $\Gamma$ in \cite{GKe}.  Note that the spectral surface is smooth, while the  space 
$\overline \Sigma_{\cal H}$ is singular. 
  The map (\ref{22Ma})  induces the {\it spectral cover} map
   \be \la{22M}
   \pi_{\cal H}: \Sigma_{\cal H} \lra M_\times.
   \ee
     It generalises the spectral cover construction from \cite{G}. 
     
        \bt   \la{THM} 

 The map (\ref{22M}) is  ramified only over  singular $\bullet-$edges and $\bullet-$vertices    of $\bS_{{\cal H}}$.

 For  an ideal alternating diagram $\cal H$,    the  spectral threefold $\Sigma_{\cal H} $ is homeomorphic to the   Lagrangian ${\Bbb L}^\circ_{{\cal H}}$. 
\et

 \begin{proof}

 1. The map $\pi$ is evidently unramified outside of the gluing points, and 
at the gluing points  in the open rectangles $r^\circ_E$. It is unramified at the $\circ-$vertices and $\circ-$edges in  the $\circ-$domains. 
Indeed, the domains ${\Bbb D}_\gamma$ sharing a $\circ-$vertex or a $\circ-$edge cover  each nearby point just once, see Figure \ref{3d7}: 
 near a  $\circ-$point  the domains ${\Bbb D}_\gamma$ are   the cones $C_{D_\circ, f}$, which cover the polyhedron $D_\circ$. 

The map $\pi$ is ramified over the $\bullet-$vertices and $\bullet-$edges. Indeed, near a $\bullet-$vertex, the domains ${\Bbb D}_\gamma$ are the complements to the cones $C_{D_\bullet, f}$. 
Therefore the  ramification index at a $\bullet-$vertex $c_{D_\bullet}$  is equal to the $\#$(faces   of ${D_\bullet}$) - 1. Similarly, the  
 ramification index at the $\bullet-$edge   over  a vertex $v$ of 
$ D_\bullet$ is equal to the $\#$(faces   of  $ D_\bullet$ meeting at   $v$)    $-1$.

  \begin{figure}[ht]
\centerline{\epsfbox{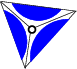}}
\caption{The polyhedron $D_\circ$  is the union of    blue domains ${\Bbb D}_\gamma \cap D_\circ$} 
\label{3d7}
\end{figure} 

  The domains ${\cal D}_\gamma$ are smooth. So the  singularities of   $\Sigma_{\cal H}$ are contained in the set of   glued points.  
Clearly the points where we glue two domains over an open rectangle $r^\circ_E$ are smooth. It is easy to see that the points on the open parts of the   
segments connecting the point  $c_D$ with   vertices of a  polyhedron  $D$ are   smooth. A simple local analysis, using the description of  domains ${\Bbb D}_\gamma$ near the  point $c_D$ as the cones $C_{D, f}$ or their complements, shows that  the points $c_D$ are smooth as well. So  $\Sigma_{\cal H}$  can be singular only at the points corresponding to the triple intersection points of ${\cal H}$. The rest follows from Lemma \ref{L3.8}. \vskip 2mm

2. Since the alternating diagram ${\cal H}$ is ideal, the conormal bundle  $T^*_\gamma M$ to the cooriented zig zag surface $\gamma$ is homeomorphic to the punctured domain ${\Bbb D}_\gamma^\circ$. Let us isotope
 zig-zag surfaces $\gamma\lra \gamma'$ by  pushing the cone boundaries  $\partial C_{D,f}$ towards the corresponding face $f$. Then $\gamma'$ is just the    the corresponding surface $H_i$ of ${\cal H}$. Adding back the cones $C_{D,f}$ we recover the domains ${\Bbb D}_\gamma^\circ$. But since 
 $
 D = \cup_f C_{D,f},
 $
 where the union is over all faces $f$ of $D$, adding all these cones is equivalent to gluing the zero sections of $T^*D$ for all colored domains $D$. 
 
 \end{proof}

Since  the map $\pi_{\cal H}$  is a cover map in   codimension $\leq 2$,  the number of  points in the generic fiber is the same. It is called the {\it degree} of the map $\pi_{\cal H}$.  
It   is an invariant of a diagram, called {\it the rank}.

\subsubsection{Singularities of the Lagrangian ${\Bbb L}_{\cal H}^\circ$ assigned to an alternating diagram in $M$}
 
 \bl \la{LagH} Let ${\cal H}$ be an alternating diagram in a threefold $M$. Then  the Lagrangian ${\Bbb L}^\circ_{\cal H}$ is  singular.  Its singularities $t_\circ$ are  the zero covectors  $0 \in T_t^*M$ at the triple intersection points $t$ of ${\cal H}$. A small ball around   $p_\circ$ is homeomorphic to a cone over a pair of pants.  
\el

\begin{figure}[ht]
\centerline{\epsfbox{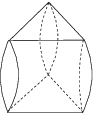}}
\caption{For an alternating diagram ${\cal H}$ in a threefold, the Lagrangian ${\Bbb L}^\circ_{\cal H}$ near a singular point is   a cone over a pair of pants.}
\label{3d6}
\end{figure}

\begin{proof} Given a triple intersection point $t$, take 
  a little sphere in $T^*X$ containing $0 \in T^*_tX$.  Its intersection   with ${\Bbb L}^\circ_{\cal H}$ is obtained by gluing two triangles $\tau_\bullet$ and $\tau_\circ$,  provided by the intersection with the zero sections over the $\bullet$ and $\circ$ domains sharing the point $p$, with three discs provided by the conormal bundles to the 
surfaces intersecting at $t$. We glue the triangles and discs as shown on Figure \ref{3d6}. Each disc can be viewed as a rectangle, and we glue two its opposite sides to   two triangles, leaving the other two sides untouched. We get a sphere with three holes. Each hole is a bigon whose sides are unglued sides of the rectangles. 
\end{proof}

  \bl \la{L3.8} The point  $t'\in \overline \Sigma_{\cal H}$ corresponding to a   triple intersection point  $t\in M$  is singular. A small ball around $p'$ is homeomorphic to a cone over a pair of pants. 
Therefore $ \overline \Sigma_{\cal H}$ is a space with boundary.
\el

\begin{proof}   A neighborhood of a point $t'$ is homeomorphic to a neighborhood of the Lagrangian ${\Bbb L}^\circ_{\cal H}$ near the point $0 \in T^*_tM$. The singularities of the Lagrangian ${\Bbb L}^\circ_{\cal H}$ were described in Lemma \ref{LagH}. \end{proof}

\subsubsection{Spectral threefolds  associated with ${\cal Q}-$digrams.} \la{Sect4.0.5}


   Given a ${\cal Q}-$diagram ${\cal Q}$ in a threefold $M$, let us resolve its quadruple intersection points, getting an alternating diagram ${\cal H}$. 
   Then there is the spectral space $\Sigma_{\cal H}$ and the spectral cover map (\ref{22M}). 
   
   Cutting out little balls $B_\bullet$ surrounding the $\bullet-$tetrahedra in ${\cal H}$ arising from the quadruple intersection points in ${\cal Q}$, we get a threefold with boundary $M_\times$. Let $p:T^*M \lra M$ be the canonical projection.  Set 
   $$
  \Sigma_{\cal Q} :=  \Sigma_{\cal H} \cap  p^{-1}(  M_\times).
   $$ 
   We call    $\Sigma_{\cal Q}$ the {\it spectral threefold  assigned to ${\cal Q}$}.   
      
   \bl \la{L4.6} The  spectral threefold $\Sigma_{\cal Q}$  is a smooth threefold with boundary. 
   
   If ${\cal Q}$ is an ideal diagram, the spectral threefold   $\Sigma_{\cal Q}$ is homeomorphic to the Lagrangian ${\Bbb L}^\circ_\times$:
    $$
    \Sigma_{\cal Q} = {\Bbb L}^\circ_\times.
        $$
   \el
   
  \begin{proof}  The first claim follows from Lemma  \ref{L1.20}. The second  follows from Theorem \ref{THM}. \end{proof}
   
   If ${\cal Q}$ is ideal, the projection $T^*M \to M$    induces a projection of the spectral threefold  to $M$:
    \be \la{22MQ}
   \pi_{\cal Q}: \Sigma_{\cal Q} \lra M_\times.
   \ee
 It is  called the {\it spectral cover} map.

\section{Cluster    Lagrangians in  non-commutative character varieties} \la{Sec7.3}

   \subsection{Two flavors of  cluster varieties from ${\cal Q}-$diagrams  of discs}  \la{SECT6}

\subsubsection{Cluster varieties from alternating  diagrams ${\cal T}$ on a surface}    Given such a ${\cal T}$, there is the stack 
$  
 {\cal A}_{[{\cal T}]} $ with 
Zariski open substack:
      \be \la{ZOPa}
     {\cal A}_{\cal T}^{\circ}\subset {\cal A}_{[{\cal T}]}.
     \ee
     
Next, there is   the substack ${\cal U}_{[ {\cal T}]}\subset {\cal X}_{[{\cal T}]}$  of sheaves with {\it trivial} rank one microlocalization at the conormal bundles to the strands of ${\cal T}$. It has a  substack ${\cal U}^\circ_{{\cal T}}\subset {\cal U}_{[{\cal T}]}$ 
of the sheaves \underline{vanishing on  mixed domains}.  
  These  substacks  are non-commutative tori. 
 
 The  collections of the tori $\{{\cal A}^\circ_{{\cal T}}\}$ and $\{{\cal U}^\circ_{{\cal T}}\}$, assigned to  alternating diagrams  in the given admissible deformation class,  
provide non-commutative cluster varieties structures of two types  \cite{GKo}: 
 \vskip 1mm
 

   The  stack ${\cal A}_{[{{\cal T}}]}$ is  a non-commutative cluster ${\cal A}-$variety,  equipped with the canonical 2-form $\Omega_{\cal A}$.  
   
  The  stack ${\cal U}_{[ {\cal T}]}$  is a 
     non-commutative cluster symplectic variety with the symplectic form $\Omega$. 
     
      \vskip 2mm
     
Recall the groupoid ${\rm Loc}_1({\bf \Sigma}_{\cal T})$ of flat line bundles on the compactified spectral surface ${\bf \Sigma}_{\cal T}$. 
  One has
 \be \la{MDIS} 
 \begin{split}
& {\cal U}^{\circ}_{ {\cal T}}  \stackrel{}{=}    {\rm Loc}_1({\bf \Sigma}_{\cal T}).
\end{split}
\ee
Let  ${\rm Loc}^{\rm triv}_1({\bf \Sigma}_{\cal T})$ be the groupoid of flat line bundles on  ${\bf \Sigma}_{\cal T}$, trivialized at the former punctures. Then
    \be \la{MDIS*} 
{\cal A}_{\cal T}^\circ  \stackrel{\sim}{=}    {\rm Loc}^{\rm triv}_1({\bf \Sigma}_{\cal T}).
\ee

 The   torus ${\cal A}^\circ_{{\cal T}}$ carries  a  non-commutative 
 2-form $\Omega_{\cal A}$, see  (\ref{AF}). The   2-form $\Omega_{\cal A}$  descends under the natural map ${\cal A}^\circ_{{{\cal T}}} \lra {\cal U}^\circ_{{{\cal T}}}$, which forgets  trivializations,  to 
torus (\ref{MDIS}), providing a symplectic form $\Omega$ there.   

\subsubsection{Cluster varieties from ${\cal Q}-$diagrams of discs in a threefold}

 Let us now turn to ${\cal Q}-$diagrams in a threefold $M$. 

\bd Given  a ${\cal Q}-$diagram  of cooriented \underline{discs} ${\cal Q}$    in a threefold $M$,  we consider: 
    \vskip 1mm
    
 
 1)   The stack ${\cal A}_{[{\cal Q}]}$ 
parametrising  ${\cal F}\in {\cal X}_{[{\cal Q}]}$ $+$  \underline{trivializations}  of their microlocal support line bundles.

    2)  Its substack, specified by the extra condition of  \underline{vanishing on mixed domains}:
        \be \la{ZOP}
      {\cal A}_{\cal Q}^{\circ}\subset {\cal A}_{[{\cal Q}]}.
     \ee
   \ed
  The image of embedding  (\ref{ZOP})  is  Zariski open in the target. 
   Forgetting  trivializations   we get projections
 \be \la{CPRO}
 {\cal A}^\circ_{\cal Q} \lra  {\cal X}^{\circ}_{\cal Q},  \ \ \ \   {\cal A}_{[{\cal Q}]}\lra {\cal X}_{[{\cal Q}]}.
 \ee
They are quotients by the action of the   {group} 
$
{\cal R}:= {R^\times}^{\{\mbox{surfaces of  ${\cal Q}$}\}},  
$ 
acting  by rescaling  trivializations.  
  \vskip 1mm

Assume now that a ${\cal Q}-$diagram of cooriented discs ${\cal Q}$ intersects the boundary $\partial M$ transversally by an alternating diagram of loops ${\cal T}$. 
Then the restriction functor provides functors
\be
\begin{split}
&{\rm Res}: {\cal X}_{[{\cal Q}]} \lra {\cal U}_{[ {\cal T}]}.\\
&{\rm Res}: {\cal X}^\circ_{{\cal Q}} \lra {\cal U}^\circ_{{\cal T}}.\\
\end{split}
 \ee 
 Indeed, since any local system on the disc is trivial, the image of the restriction functor is a sheaf with trivial microlocalization at the punctured conormal bundle to any zig-zag loop. 
 So it lies in ${\cal U}_{[ {\cal T}]}$. 
 We set
  \be \la{14} \begin{split}
& {\cal L}_{[{\cal Q}]} := {\rm Im} ~{\rm Res}({\cal X}_{{[\cal Q} ]}) \subset    {\cal U}_{[{\cal T}]}.    \\
& {\cal L}^{\circ}_{\cal Q} := {\rm Im} ~{\rm Res}({\cal X}^\circ_{\cal Q}) \subset    {\cal U}^\circ_{\cal T}\stackrel{(\ref{MDIS})}{=}  {\rm Loc}_1({\bf \Sigma}_{\cal T}).     \\
\end{split}
\ee 
There are restriction functors for  the ${\cal A}-$stacks:
  \be \la{14} \begin{split}
& {\cal L}^{\cal A}_{[{\cal Q}]} := {\rm Im} ~{\rm Res}({\cal A}_{{[\cal Q} ]}) \subset    {\cal A}_{[{\cal T}]}.    \\
& {\cal L}^{{\cal A}, \circ}_{{\cal Q}} := {\rm Im} ~{\rm Res}({\cal A}^\circ_{{\cal Q} }) \subset    {\cal A}^\circ_{{\cal T}}\stackrel{}{=}  {\rm Loc}_1^{\rm triv}({\bf \Sigma}_{\cal T}).     \\
\end{split}
\ee

 \subsection{${\cal A}-$coordinates from ${\cal Q}-$diagrams of  discs}  \la{SECTAC}

 We start from the definition of the graph of singularities ${\cal S}_{\cal Q}$ of any ${\cal Q}-$diagram.

\bd   \la{GS} The  intersection points of a  ${\cal Q}-$diagram of surfaces ${\cal Q}$  form the {\em graph of singularities}   ${\cal S}_{\cal Q}$. Its vertices are the  quadruple  intersection points, and the edges are segments of  the  double intersections. 
\ed

For example, the  singularity graph for the basic ${\cal Q}-$diagram in the cube ${\rm C}$  is given by the twelve edges connecting the center   with  midpoints of the cube edges. 
  We orient the edges out of the center.  
The edges of the singularity graph match the cube edges. So we can think that the 
${\cal A}-$coordinates for the basic ${\cal Q}-$diagram are assigned to the (unoriented) edges of the cube.  
So
\vskip 1mm

$\{${\it The edges of the singularity graph of the basic diagram}$\} \ \ \stackrel{1:1}{\Longleftrightarrow} \ \ \{${\it The cube edges}$\}$.

\vskip 2mm

\begin{figure}[ht]
\centerline{\epsfbox{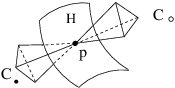}}
\caption{The two triangular cones  $C_\circ, C_\bullet$ with the vertex at the   quadruple point $p$ on the surface $H$.} 
\label{AB10}
\end{figure}

Any quadruple intersection point $p$ of a ${\cal Q}-$diagram is locally diffeomorphic to  the basic ${\cal Q}-$diagram. So there are $6$ singular lines passing through the point $p$, providing  $12$ edges 
 sharing $p$. 
 For each pair $(p, H)$ where $p\in H$ is a quadruple intersection point on a surface $H$ of  a ${\cal Q}-$diagram  there are two  triangular cones $C_\bullet, C_\circ$ with the vertex $p$, see Figure \ref{AB10}. 
 They are given by the other three surfaces intersecting at $p$. 
On Figure \ref{AB8},  these are the cones over the two triangles on the cube   containing the vertices of the  diagonal perpendicular to  $H$.

\begin{figure}[ht]
\centerline{\epsfbox{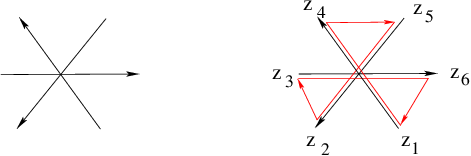}}
\caption{A triple crossing diagram and the relation $z_3z_6 \cdot z_1z_4 \cdot z_5z_2=-1$.}
\label{ncc1}
\end{figure}

We assign      to   the \underline{oriented} edges $E$ of the 
singularity graph ${\cal S}_{\cal Q}$ the  {\it ${\cal A}-$variables} $ a_E\in R^\times $ subject to the following relation. Let  $\overline E$ be the oriented edge obtained  by reversing the orientation of  $E$. Then
\be \la{INV}
a_{\overline E} = a_{E}^{-1}.
\ee
Let us  introduce a system of equations ${\cal R}_{\cal Q}$ on the ${\cal A}-$variables $\{a_E\}$. 
  Let $a,b,c$ be the ${\cal A}-$variables on the  edges of the cone $C_\ast$,  oriented out of   $p$,   ordered counterclockwise, looking from the point $p$.
For each pair $(H, p)$, where $p$ is the quadruple point  on the surface $H$ of the diagram,   we have the  triple crossing diagram $H \cap {\cal Q}$ on   $H$ with a crossing point $p$.  Let  $z_1, z_2,z_3, z_4, z_5, z_6\in R^\times$ be the ${\cal A}-$variables  
on the  edges of the singularity graph   ${\cal S}_{\cal Q}$   which lie in $H$,  share the vertex $p$, and oriented out of  $p$, see Figure \ref{ncc1}.

\bd  \la{Def3.8} Let  ${\cal Q}$ be a ${\cal Q}-$diagram  of cooriented discs in an oriented threefold. Consider the following   set of equations ${\cal R}_{\cal Q}$  on   ${\cal A}-$variables   at the  oriented edges $E$ of the singularity graph ${\cal S}_{\cal Q}$, satisfying    (\ref{INV}): 

For each pair $(H, p)$, where $p\in H$ is the  quadruple intersection point  of ${\cal Q}$, we have:
 \vskip 1mm

1) For each of  the two cones,  $C_\bullet$ and $C_\circ$,  the monomial relation
\be \la{M1}
abc=-1.
\ee

2)  
\be \la{eqqQ}
z_1z_4+ (z_5z_2)^{-1}=1, \ \ \ \ z_5z_2+ (z_3z_6)^{-1}=1, \ \ \ \ z_3z_6+ (z_1z_4)^{-1}=1.  \ \ \ \ 
\ee
 \be \la{MR92Q}
z_3z_6 \cdot z_1z_4 \cdot z_5z_2=-1.
\ee 
\ed

 Each of the four surfaces $H_i$ containing the quadruple intersection point $p$ provides  three equations (\ref{eqqQ}). So we get twelve equations at  $p$. 
By Proposition \ref{P3.11}, the monomial equation (\ref{MR92Q}) follows from the ones $(\ref{M1})$, and 
 the twelve equations (\ref{eqqQ}) at $p$  are equivalent.

\bt \la{MTHGK} Let ${\cal Q}$ be a ${\cal Q}-$diagram of cooriented discs  in an oriented threefold $M$ intersecting  the boundary $\partial M$ by  an alternating diagram ${\cal T}$. 
Then:

\begin{enumerate}

 \item The ${\cal L}^{\cal A}_{[{\cal Q}]}$   is  a non-commutative  cluster Lagrangian in the cluster  ${\cal A}-$variety  ${\cal A}_{[{\cal T}]}$.  
 
 \item The cluster atlas for the Lagrangian ${\cal L}^{\cal A}_{[{\cal Q}]}$ is given by the Lagrangians  ${\cal L}^{{\cal A}, {\circ}}_{{\cal Q}}\subset {\rm Loc}_1^{\rm triv}({\bf \Sigma}_{\cal T})$, where 
 diagrams  ${\cal Q}$ are in the  same admissible deformation class.   
 Each of them is    described   by the systems of equations ${\cal R}_{\cal Q}$.  
 
 \item  Similarly, 
 ${\cal L}_{[{\cal Q}]}$   is  a non-commutative  cluster Lagrangian in the symplectic stack $  {\cal U}_{[{\cal T}]}$.

\item   In the commutative case, ${\cal L}_{[{\cal Q}]} $  is a {\it $K_2-$Lagrangian}  in the $K_2-$cluster symplectic stack ${\cal U}_{[{\cal T}]}$. \end{enumerate}
 \et


\paragraph{Remark.}

The substacks ${\cal L}_{[{\cal Q}]}$ and  ${\cal L}^{\cal A}_{[{\cal Q}]}$    can be defined for arbitrary ${\cal Q}-$diagrams. They carry ${\cal A}-$coordinates, defined below. However they are not  Lagrangian 
if the surfaces of the diagram are not discs. 
\vskip 2mm

        \begin{proof}   (1) Follows from Theorem \ref{THE4.4}.  
        \vskip 1mm

        (2)  
  Given an object ${\cal F}$ of ${\cal A}^\circ_{\cal Q}$, let us define  elements $a_E \in R^\times$ at the oriented edges $E$ of the singularity graph ${\cal S}_{\cal Q}$, called  the 
{\it ${\cal A}-$coordinates}. 
   Any edge $E$ is contained in    just two surfaces:
$$
E \subset H_1  \cap H_2.
$$
The  transversal disc $\pi_E$ to   $E$  contains two 
cooriented lines $\gamma_i := \pi_E \cap H_i$. 
The orientations of the edge $E$ and the   threefold $M$ provide  an orientation of the disc $\pi_E$.  The lines are ordered, say $h_1, h_2$, by the orientation of $\pi_E$.
  Consider the unique $\bullet-$sector $b_{h_1, h_2} \subset \pi_E$, and a  path 
$\gamma$ transporting the conormal to $h_1$ to the conormal  to $h_2$ inside the sector. Denote by $\gamma(s_{H_1})$ the parallel transform of the section $s_{H_1}$   along the path $\gamma$. We define an element $a_E \in R^\times$ by writing 
$$
\gamma(s_{H_1}) = a_E s_{H_2 }.
$$

 The claim that the 2-form $\Omega$ vanishes on ${\cal L}^{{\cal A}, \circ}_{\cal Q}$ boils down to the calculation in  \cite[Section 6]{GKo}.  
The subspace $
{\cal L}^\circ_{\cal Q}  \subset {\rm Loc}_1^{\rm tr}(\Gamma) \stackrel{(\ref{MDIS}) }{=} {\cal U}_{{\cal T}}
$  is Lagrangian by (1). Therefore ${\cal L}^{{\cal A}, \circ}_{\cal Q}$ maximal isotropic. 

  So   for each quadruple intersection point $p$ there are
  $8$ monomial equations 
  (\ref{S14}), and equations (\ref{S1}).
  
  \bl \la{UVIN} Equations (\ref{S14})+(\ref{S1}) are equivalent to the equations (\ref{E1})+(\ref{E2}).
\el

\begin{proof} Given an edge $E$ of the graph $\Gamma_{\rm cube}$, there are two zig-zags $\gamma_{E-}$ and $\gamma_{E+}$ containing the midpoint $e$ of  $E$. 
So there are two trivializations $s_+$ and $s_-$ 
at the fiber ${\rm L}_e$ of a flat line bundle $\rm L$ over $e$. Then 
$
s_+= a_{\rm E} s_-$, where $a_{\rm E}\in R^\times.
$
The $a_{\rm E}$ is the ${\cal A}-$variable assigned to the edge $E$. The claim follows. 
\end{proof}
        (3)    The  projection 
$
    {\cal A}_{[{\cal Q}]}\lra {\cal U}_{[{\cal Q}]} 
 $ 
  is  the  quotient by the action of the    {group}  $
 {R^\times}^{\{\mbox{surfaces of  ${\cal Q}$}\}} 
$   by rescaling  trivializations. 
  Then (3) follows  from (1) and (2). 
        \vskip 1mm
        
         (4) This is obvious from the general shape $x+y=1$ of equations (\ref{eqqQ}). See also Lemma \ref{K2}. \end{proof}


 \subsection{${\cal Q}-$diagrams of  discs  from   ideal triangulations of threefolds} \la{sec5.3}
 
In Section \ref{sec5.3} we consider  threefolds $M$,  possibly with boundary,  glued from simplexes. 
Some faces of the simplexes may not belong to $M$.   For example,  take a compact threefold  and remove  a finite non-empty 
 collection of points,  referred to as punctures. Hyperbolic threefolds with cusps are in this category.

Given a collection of  {\it special points}  on the boundary of  $M$,   an {\it ideal triangulation} ${\tau}$ of $M$ is a triangulation with the vertices at these points. 
 Given an integer $m \geq 2$ and  an ideal triangulation $\tau$ of a   threefold $M$, we  define an ideal ${\cal Q}-$diagram ${\cal Q}_{\tau, m}$ in $M$. It provides a singular surface
 $$
 \bS_{{\tau}, m}\subset M. 
 $$
The construction from Section \ref{Sect4.0.5}   assigns to the ${\cal Q}-$diagram ${\cal Q}_{{\tau},m}$ 
the  $m:1$ spectral cover:
   \be \la{mpi}
\pi: \Sigma^\times_{{\tau}, m}\lra M^\times, 
\ee
as well as  its completion over the singular points $\pi: \Sigma_{{\tau}, m}\lra M$.
  
\subsubsection{Hypersimplicial decompositions.} Given a pair of integers $p,q \geq 0$, the hypersimplex $\Delta^{p,q}$ is a polyhedron of the dimension   $p+q+1$,   defined as an integral hyperplane section of the unit cube: 
  $$
  \Delta^{p,q}:= \{x\in \R^{p+q+2}~~ |~~ x_0+ ... + x_{p+q+1}=p+1, ~~ 0 \leq x_i \leq 1\}.
  $$  
So $\Delta^{0,0}$ is a segment; $\Delta^{0,1}$ and $\Delta^{1,0}$ are triangles; $\Delta^{2,0}$ and $\Delta^{0,2}$ are tetrahedra, and $\Delta^{1,1}$ is an octahedron. \vskip 2mm
 
 The boundary of the hypersimplex $\Delta^{p,q}$ is the union of $p+q+2$ hypersimplices of type  $\Delta^{p-1,q}$ and $p+q+2$ hypersimplices of type  $\Delta^{p,q-1}$. 
 It is obtained by setting one of the coordinates $x_i$ to $0$ or $1$.\footnote{Note that $\Delta^{-1,q}$ and $\Delta^{p,-1}$ are the empty sets.}\vskip 2mm
 
{\it Examples}. 1. The boundary of the simplex $\Delta^{n,0}$ consists of $n+1$ simplices of type $\Delta^{n-1,0}$.
  
 2.  The boundary of the octahedron $\Delta^{1,1}$ consists of the four $\Delta^{1,0}-$triangles,   called    
   $\circ-$triangles, and four $\Delta^{0,1}-$triangles,   called  $\bullet-$triangles. 
   Each boundary  triangle shares its edges with  triangles of different type.      
For the octahedron  on Figure (\ref{3d2}) the $\bullet-$triangles   are   on the faces of the tetrahedron.  \\
 
 Given  integers $d\geq 1$,  $m\geq 2$, take the  hyperplane section of the positive cone in $\R^{d+1}$:
 \be \la{104}
 {\rm T}_{d}:= 
 \{x\in \R^{d+1}~~ |~~ x_0+ \ldots +x_{d}=m, ~~ x_i \geq 0\}.
 \ee 
 The {\it hypersimplicial decomposition}  ${\cal H}_{{\rm T}_{d}, m}$  
 of the simplex ${\rm T}_{d}$ in (\ref{104})  is given by the $d+1$ families of $m-1$  hyperplanes parallel to one of the faces, defined  by the equations 
  $x_i =1, ..., m-1$. It cuts the simplex ${\rm T}_{d, m}$ into hypersimplices of type $\Delta^{p,q}$ where $p+q=d-1$. \vskip 2mm

{\it Examples}. 1. If $d=2$ we get an $m-$decomposition of a triangle into looking up and down triangles. 
   
\begin{figure}[ht]
\centerline{\epsfbox{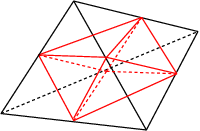}}
\caption{The hypersimplicial decomposition ${\cal H}_{{\rm T}, 2}$ of  a tetrahedron   into four $\Delta^{2,0}-$tetrahedra   and an octahedron.} 
\label{AB12a}
\end{figure}

 2.  The $d=3$ case. We consider the  tetrahedron    \be \la{TETm}
 {\rm T}:= \{x\in \R^4~~ |~~ x_0+x_1+x_2+x_3=m, ~~ x_i \geq 0\}, ~~~~~~m \geq 2.
 \ee 
  The {\it hypersimplicial decomposition}  ${\cal H}_{{\rm T}, m}$  cuts the tetrahedron ${\rm T}$ into    $\Delta^{2,0}-$tetrahedra    
  (looking up),   octahedra $\Delta^{1,1}$, and  $\Delta^{0,2}-$tetrahedra   (looking down), see Figure \ref{AB12a}. \vskip 2mm
  
 \bd   
  Given an integer $m \geq 2$ and  an ideal triangulation ${\tau}$ of a manifold $X$, the  {\it hypersimplicial decomposition} 
  ${\cal H}_{{\tau},m}$ of $X$ is   the union of   hypersimplicial decompositions ${\cal H}_{{\rm T},  m}$ of each simplex ${\rm T}$ of ${\tau}$:   
  \be \la{QTmh}
X = {\cal H}_{{\tau},m} = \cup_{{\rm T} \in {\tau}} {\cal H}_{{\rm T}, m}.
\ee
 \ed

\subsubsection{The singular hypersurface $\bS_{{\tau}, m}\subset X$.}
 
\begin{figure}[ht]
\centerline{\epsfbox{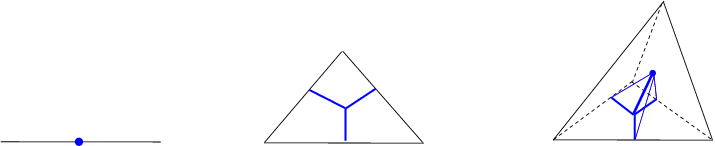}}
\caption{The point $\bS^{0,0}$. The singular graph $\bS^{1,0}$.   A quarter of the  singular surface $\bS^{2,0}$.} 
\label{3d1}
\end{figure}

\begin{figure}[ht]
\centerline{\epsfbox{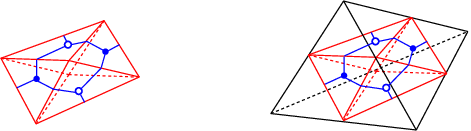}}
\caption{The  surface $\bS^{1,1}$ in the  octahedron is the cone over the   bipartite graph  on its boundary.} 
\label{3d3}
\end{figure}

 We define inductively a singular polyhedral hypersurface  
 $$
 \bS^{p,q}\subset \Delta^{p,q}.
 $$
The point $\bS^{0,0}$ is the center of the segment $\Delta^{0,0}$. 
\bd  \la{INDD} 
 The hypersurface $\bS^{p,q}$ is the cone over   hypersurfaces 
$\bS^{p-1,q}$ and  $\bS^{p,q-1}$ in the   boundary hypersimplices,  centered at the center of  $\Delta^{p,q}$,  see Figures \ref{3d1},  \ref{3d3}. \vskip 2mm
\ed
 
 \bd Given an integer $m \geq 2$ and  an ideal triangulation ${\tau}$ of $X$,  the singular hypersurface 
 $$
\bS_{{\tau}, m} \subset X 
 $$
is the union of  singular hypersurfaces $\bS^{p,q} $  in the hypersimplices 
of   hypersimplicial decomposition (\ref{QTmh}).
\ed

 The singular hypersurfaces $\bS^{p,q} $ fit together since  they induce on each of the faces of the hypersimplicial decomposition the  same singular hypersurface.  
 
When $X$ is a surface $S$, the surface graph $\bS_{{\tau}, m}$  is  the bipartite graph $\Gamma_m\subset S$ assigned to an integer $m \geq 2$ and a triangulation of a surface in \cite{G}. 
\vskip 1mm

{\it The singular hypersurfaces $\bS_{{\tau}, m}\subset X$ are higher dimensional analogs of bipartite graphs}. \vskip 1mm

Let us elaborate the case when $X$ is a threefold $M$. 


\subsubsection{Hypersimplicial decompositions of a threefold.}  
 
The {\it hypersimplicial decomposition ${\cal H}_{\tau, m}$}    of a threefold $M$ consists of  the octahedra ${\rm O}_p$,  
    the   tetrahedra $t_\circ$, which are $\Delta^{2,0}-$hypersimplices,  
 and the   tetrahedra $t_\bullet$, which are $\Delta^{0,2}-$hypersimplices:
\be
M = \cup_p{\rm O}_p \cup_{\bullet} t_{\bullet} \cup_{\circ} t_{\circ}.
\ee
 The  tetrahedra   are glued to  the  octahedra. The colors of    triangles are consistent:  a $\bullet-$triangle of a
tetrahedron is glued to a $\bullet-$triangle of the octaherdron, and the same for the $\circ-$triangles. 
The octahedrons ${\rm O}_p$ correspond  to the quadruple intersection points $p$ of the ${\cal Q}-$diagram ${\cal Q}_{{\rm T},m}$.

\subsubsection{The singular surface $\bS_{{\tau}, m}\subset M$.}  Singular surfaces $\bS^{2,0}$ / $\bS^{0,2}$ are  obtained by connecting the center of the   tetrahedron with the graphs  on  its faces, see Figure \ref{3d1}. 
The edges of the singular locus of surfaces $\bS^{0,2}$ 
/ $\bS^{2,0}$ 
are called  the singular $\bullet-$edges / $\circ-$edges. The centers of the  tetrahedra $\bS^{0,2}$ / $\bS^{2,0}$   are  called the singular $\bullet-$vertices / $\circ-$vertices.

The  surface $\bS^{1,1}\subset \Delta^{1,1}$ is the cone over the bipartite graph on the boundary of the octahedron. Its singular locus 
is given by the four $\bullet-$edges connecting the center with the $\bullet-$vertices of the boundary triangles, and four similar $\circ-$edges. 
Given a triangulation $\tau$ of $M$, the singular surface $\bS_{{\tau}, m}\subset M$ 
  is the union of   surfaces $\bS^{2,0}, \bS^{0,2}, \bS^{1,1}$ in  the  hypersimplicial decomposition ${\cal H}_{{\tau}, m}$. 
Its singular locus consists of the $\bullet-$edges and $\circ-$edges. The endpoints of  singular edges are   centers of the hypersimplices of the hypersimplicial decomposition.

  \subsubsection{Zig-zag surfaces for the  singular surface $\bS_{{\tau}, m} \subset   M$.}     
 Each   singular surface  $\bS^{p,q}\subset \Delta^{p,q}$,   $p+q=2$, carries four cooriented {\it zig-zag surfaces}. Their coorientations    consistent with   coorientations of zig-zag strands on  bicolored  graphs  induced on the boundary of $\Delta^{p,q}$.   Precisely, they look as follows.

\begin{figure}[ht]
\centerline{\epsfbox{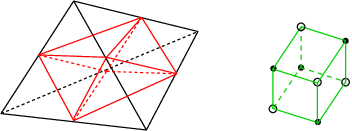}}
\caption{The green cube is inscribed into the red octahedron. 
  Its  $\bullet-$vertices are the centers of the octahedron $\bullet-$triangles - the ones on the faces of the black tetrahedron. The $\circ-$vertices are the centers of the octahedron $\circ-$triangles.} 
\label{3d2}
\end{figure}

Zig-zag surfaces for    $\bS^{2,0} $ /  $\bS^{0,2} $  are isotopic to hemisphers around    vertices of  the  tetrahedra $\Delta^{2,0}$ / $ \Delta^{0,2}$. They are isotopic to a part of   $\bS^{2,0}$ /  $\bS^{0,2}$ surronding the vertex, and 
  coorientated  
  out/towards  the center of  $\bS^{2,0}$ / $ \bS^{0,2}$. 
Zig-zag surfaces for   $\bS^{1,1}$   described as follows. 
  In  the   octahedron  $\Delta^{1,1}$  on Figure \ref{3d2},   take 
 the four conical surfaces in the green cube ${\rm C}_p$      given by the centered at $p$ cones    over  the cooriented\footnote{cooriented towards the $\bullet$ vertex of the cube which is not on the zig-zag strand.}  zig-zag hexagons $\gamma_1, ..., \gamma_4$  on Figure \ref{AB1}, and  expand them to the  octahedron. \vskip 2mm
 
Concatenating 
   zig-zag surfaces in the hypersimplices we get   zig-zag surfaces   for the  hypersiplicial decomposition  ${\cal H}_{{\tau}, m}$ of $M$.    
 On the other hand, let us consider  the ${\cal Q}-$diagram ${\cal Q}_{{\tau},m}$.

\subsubsection{The ${\cal Q}-$diagram   of surfaces from an ideal triangulation of a threefold $M$.}  \la{Sect5.3.6}

\bd The  ${\cal Q}-$diagram ${\cal Q}_{{\rm T},m}$   is given by the following $4m$ cooriented  planes in the tetrahedron   (\ref{TETm}):  
 \be \la{Q4} 
 x_i=  k + \frac{1}{2}, ~~ i\in \{0,1,2,3\}, ~~ k = 0, ..., m-1.
 \ee
 \ed

Each plane  is parallel to one of  the faces      and cooriented towards  this face. 
 The  quadruple intersections $q$  are   the centers of  octahedrons $\Delta^{1,1}$. They have  positive half-integral coordinates,  summing to $m$: 
  $$
  q=(q_0, q_1, q_2, q_3); ~~~q_i \in \frac{1}{2}+\Z_{\geq 0}, ~~q_0+q_1+q_2+q_3=m.
   $$
  The ${\cal Q}-$diagram  ${\cal Q}_{{\rm T},m}$   cuts  a standard triple crossing diagram   on each   face.    So we can concatenate them   for the tetrahedra ${\rm T}$ of an ideal triangulation ${\tau}$ of    $M$,   getting the  ${\cal Q}-$diagram   
 \be \la{SSQ}
    {\cal Q}_{{\tau},m} := \cup_{{\rm T} \in {\tau}}{\cal Q}_{{\rm T},m}\subset M.
 \ee  
It gives rise to the singular surface $\bS^{\cal Q}_{{\tau}, m}$.

\subsubsection{The two kinds  of  singular surfaces in threefolds coincide.} Given a triangulation ${\tau}$ of   $M$, and an integer $m \geq 2$, there are    two  collections of singular surfaces in $M$:

 \begin{enumerate}
 \item  The singular surface $\bS^{\cal H}_{{\tau}, m}$, assigned  to the $m-$hypersimplicial decomposition ${\cal H}_{{\tau}, m}$.
  
 \item  The singular surface $\bS^{\cal Q}_{ {\tau}, m}$ assigned to the ${\cal Q}-$diagram ${\cal Q}_{{\tau}, m}$    in (\ref{SSQ}).  
 \end{enumerate}
 
  \bp
 The singular surfaces $\bS^{\cal H}_{{\tau}, m}$  and  $\bS^{\cal Q}_{ {\tau}, m}$ coincide, matching 
   their  zig-zag surfaces.  
 So the related two spectral covers of $M$ coincide. 
 \ep
 
\begin{figure}[ht]
\centerline{\epsfbox{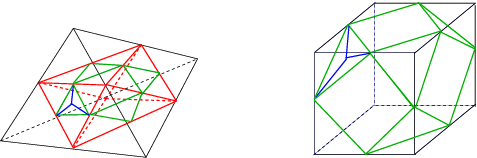}}
\caption{On the left: the red octahedron of the $m=2$ hypersimplicial decomposition of the tetrahedron. 
The  ${\cal Q}-$diagram with the quadruple intersection in the center   cuts   eight green triangles on the octaherdon. 
On the right: the intersection of the ${\cal Q}-$diagram with a cube have eight green triangles, matching the  green triangles on the left. The domains inside of the green triangles are colored. } 
\label{3d2b}
\end{figure}

 \begin{proof} It is sufficient to prove this for a  tetrahedron ${\rm T}$. We claim that $\bullet$ / $\circ-$domains are identified with the $\Delta^{2,0}$ / $\Delta^{0,2}$ simplices of the hypersimplicial decomposition of ${\rm T}$.  Figure \ref{3d2b} shows that the intersections of  singular surfaces $\bS^{\cal H}_{{\rm T}, m}$  and  $\bS^{\cal Q}_{ {\rm T}, m}$  with the 
 octahedron $\Delta^{1,1}$ of the hypersimplicial decomposition  coincide on the nose. Indeed, $\bS^{\cal Q}_{ {\rm T}, m}\cap \Delta^{1,1}$ is the cone over the eight 
 blue trivalent graphs with a single vertex sitting in the eight green triangles. These green triangles are the intersections of eight colored - that is $\bullet$ or $\circ$ -  triangular cones centered at the quadruple intersection point $q$ with the boundary of $\Delta^{1,1}$. Therefore ,by the very definition, $\bS^{\cal H}_{ {\rm T}, m}\cap \Delta^{1,1}$  is the cone over the same blue graphs. So we get
$$
 \bS^{\cal Q}_{ {\rm T}, m}\cap \Delta^{1,1} = \bS^{\cal H}_{ {\rm T}, m}\cap \Delta^{1,1}.
$$
The claim that the intersections of  singular surfaces $\bS^{\cal H}_{{\rm T}, m}$  and  $\bS^{\cal Q}_{ {\rm T}, m}$  with each of the 
 tetrahedra $\Delta^{2,0}$ / $\Delta^{0,2}$ of the hypersimplicial decomposition  also coincide on the nose  is clear from Figure \ref{3d2b} as well.  
 
By the very definitions, this identification match the corresponding zig-zag surfaces. \end{proof}

 Each   zig-zag surface of $\bS_{{\rm T}, m}$    bounds a ball/halfball containing   a single puncture $s$ of  $\partial M$, 
  and coorientated  out of the puncture.     
So these zig-zag surfaces form an ideal ${\cal Q}-$diagram. 
   The  spectral space $\Sigma_{{\rm T}, m}$ is obtained by gluing the copies of these domains. 
   There is    
  a  spectral map
 \be \la{mpi1}
\pi: \Sigma_{{\rm T}, m}\lra M.
\ee

\bc 
 The spectral space $\Sigma_{{\tau}, m}$ is homeomorphic to the  Lagrangian ${\Bbb L}^\circ_{{\tau}, m}$ associated with ${\cal Q}_{{\tau},m}$. 
   The spectral map $\pi$ in (\ref{mpi1}) is an $m:1$  cover map, ramified over the   singular $\bullet-$edges and $\bullet-$vertices    of $\bS_{{\tau}, m}$. 
Its ramification degree  is  $3$ at the singular $\bullet-$vertices   
  of $\bS_{{\tau}, m}$, and $2$ over the singular $\bullet-$edges. 
 \ec 

\begin{proof} This is a special case of Theorem \ref{THM}. The number of domains ${\cal D}_{s,i}$ over a point closed to a puncture of $\partial M$ is $m$. Therefore the degree of the map $\pi$ is also $m$. 
\end{proof}

\subsection{Cluster    Lagrangians in  non-commutative character varieties} \la{sec7.3}

\subsubsection{Character varieties on surfaces of various flavors.}    The stack ${\cal U}_m(S)$ parametrising local systems of 
$m$-dimensional $R-$vector spaces on a surface  
 $S$ with  unipotent monodromies around the punctures 
 and    an invariant flag near each puncture, \cite{FG1},  \cite{GKo}.

The stack ${\cal A}_{m}(S)$ parametrises  twisted $m-$dimensional $R-$local systems on $S$, with an invariant decorated flag near each puncture. Forgetting  decorations we get a surjective map
$
p: {\cal A}_{m}(S) \lra {\cal U}_{m}(S).
$

To fit these stacks  into the general framework, recall  a class of    bipartite ribbon graphs $\Gamma$ on $S$, called   {\it ${\rm GL}_m-$graphs}  \cite{G}. 
The collection ${\cal T}$  of  zig-zag loops on $\Gamma$ has the following features,  characterizing them:
 \vskip 2mm

i)  For each puncture $s$ on $S$ there are exactly $m$ zig-zag loops of the collection ${\cal T}$ surrounding $s$.

ii) Each zig-zag loop of the collection ${\cal T}$ surrounds just one puncture.

iii) The rank of the spectral cover is equal to $m$, see \cite{G} for details. \vskip 2mm

\bl \la{BE}
One has 
$$
 {\cal U}_{m}(S) = {\cal U}_{[{\cal T}]}, \quad {\cal A}_{m} (S)= {\cal A}_{[{\cal T}]}.
$$
 \el 
 
 Lemma  \ref{BE} explains   the role of framing in the definition of  the stacks 
  ${\cal U}_m(S)$  and ${\cal A}_m(S)$. 
 
 \begin{proof}   
 
   Let $S$ be a surface with the finite set of punctures $\rm P$. Denote by ${\cal C}_m$ the collection  given by  $m$ concentric circles near each puncture. 
 Then by the very definition we have a canonical  equivalence 
  $$
{\cal U}_{{\cal C}_m} = {\cal U}_m(S), \ \ \ \  {\cal A}_{{\cal C}_m} = {\cal A}_m(S).  
 $$ One can deform  ${\cal C}_m$ to a 
  collection of zig-zag loops  of a  ${\rm GL}_m-$graphs graph   on $S$  \cite{G}. So  the invariance of these stacks under admissible deformations proves the claim. 
  \end{proof}

    The stack ${\cal U}_m(S)$  
   carries a   non-commutative cluster symplectic structure, invariant under the mapping class group ${\rm Mod}_S$ of $S$. 
 Its cluster tori  include the ones assigned to 
 the ${\rm GL}_m-$graphs  on   $S$.

\subsubsection{Cluster Lagrangians in moduli spaces of non-commutative  local systems on surfaces.}  
  The stack ${\rm Loc}_m(S)$ of $m-$dimensional local systems on a closed surface $S$ is symplectic. 
  A threefold $M$ with the boundary $S$ gives rise to a
   Lagrangian,   given by the     local systems on $S$ which  extend to $M$:
   $$
   {\rm Loc}_m(M) \subset {\rm Loc}_m(S). 
   $$
   This   is generalized to  punctured surfaces as follows. Let $S$ be a surface with punctures,   $S'$ is 
   the  surface  compactifying $S$, and $M$ is a threefold with   the boundary $S'$ and a set ${\rm P}$ {\it special boundary points}    matching the punctures on $S$. 
  Let  ${\cal L}_m({M, {\rm P}})$ be the stack  of  non-commutative  $m-$dimensional 
framed local systems on $S$, with unipotent monodromies around the punctures, which can be extended to $M$:
\be \la{LLM}
 {\cal L}_m({M, {\rm P}})\subset {\cal U}_m(S).
 \ee

\bt \la{T4.1} The substack (\ref{LLM}) is a  {\it non-commutative \underline{cluster} Lagrangian}.    
Each ideal triangulation ${\cal T}$ of $M$ provides a system  of equations  ${\cal R}_{{\cal Q}_{{\cal T},m}}$ on the variables $\{a_E\}$, describing an open part of   $ {\cal L}_m({M, {\rm P}})$.
\et

In the commutative case we recover the main result of \cite{DGG}. \vskip 1mm
 
 For example, let $\rm D$ is a ball and $S= S^2-\{p_1, ..., p_n\}$. 
   Let  ${\cal F}_m$ be the variety of flags in an $m-$dimensional vector space. Then
   $$
    {\cal L}_m(\rm D, P)  =   (\underbrace{{\cal F}_m \times \ldots \times {\cal F}_m}_{\mbox{$n$ copies}})/{\rm GL}_m.
    $$
    Indeed, any local system on the ball is trivial. So a framing is   the same   as  $n$ flags in a fiber over a point,  obtained by parallel transporting the flags to the point, 
 considered modulo the action of   ${\rm GL}_m(R)$.  
 
 Note  that     ${\rm Loc}_m(\rm D)$ is just a point.   So the framing is  essential.

  \subsubsection{Cluster description of the Lagrangian $ {\cal L}_m({M, {\rm P}})$.}  
     
We  combine Lemma \ref{BE} with its 3d analog.  
        Namely, we extend  concentric circles of the collection ${\cal C}_m$ to  concentric hemi-spheres 
 in $M$, getting a collection of cooriented surfaces ${\cal S}_m\subset M$. We get 
  the  Lagrangian
   $$
\ \ \ \ \ \ {\cal U}_{[{\cal S}_m]}  = {\cal L}_{\cal Q} = {\cal L}_m(M, {\rm P}).
$$
 Now Theorem \ref{T4.1} follows from Theorem \ref{MTHGK} and  the following crucial result. 

   \bt One can deform 
  ${\cal S}_m$  to a ${\cal Q}-$diagram ${\cal Q}$   in $M$,  so that 
 $
 \partial {\cal Q} $ is  the collection of zig-zag loops of a ${\rm GL}_m-$graph $\Gamma$ on $S$. 
 \et

\begin{proof}

 By the construction, each surface in the diagram ${\cal Q}_{{\tau}, m}$ contains inside just one puncture $p$ of $S$, and  cooriented out of the point $p$. 
So one can admissibly deform the surfaces towards the unique punctures $p$ they contain, we get a new diagram of surfaces ${\cal Q}'_{{\tau}, m}$
 such that each puncture is surrounded by exactly $m$ surfaces. 
 The deformation does not change the related moduli spaces.   So we get an equivalence of stacks
 \be
{\cal L}_{{\cal Q}_{\tau, m}} \stackrel{\sim}{\lra}   {\cal L}_m({M, {\rm P}}).
 \ee
 
Its intersection with the cluster symplectic  torus $ {\rm Loc}_1({\bf \Sigma}_{{\tau}, m})$ in ${\cal U}_{m, S}$
  contains a subtorus defined by the system of equations ${\cal R}_{{\cal Q}_{{\tau},m}}$. It parametrises sheaves vanishing on the mixed domains for the 
  alternating diagram associated with the graph $\Gamma_{{\tau}, m}$. 
  \end{proof}

 \bt Let $M$ be a threefold  with boundary $S'$.  Let $\Gamma$ be a ${\rm GL}_m-$bipartite graph on  $S$.   Then there exists a  ${\cal Q}-$diagram ${\cal Q}_\Gamma$ in $M$  inducing on the boundary  the collection
    of the zig-zag strands of $\Gamma$. 
 \et
 
 \begin{proof}

 Take an ideal triangulation ${\tau}$ of $M$. 
 It gives rise to the  ${\cal Q}-$diagram ${\cal Q}_{{\tau}, m}$, see Section \ref{Sect5.3.6}.
 The  ideal triangulation ${\tau}$   induces an ideal triangulation  ${\tau}_S$ of $S$. Let $\Gamma_{{\tau}, m}$ be the standard ${\rm GL}_m-$bipartite graph  
 on $S$ inscribed into this triangulation \cite{G}.  The diagram ${\cal Q}_{{\tau}, m}$ 
 induces on $S$ the collection  of zig-zag strands   of the bipartite graph $\Gamma_{{\tau}, m}$. 
 Now let $\Gamma$ be any  ${\rm GL}_m-$bipartite graph on  $S$. 
 Then there exists a sequence of two by two moves of the graph $\Gamma$ which transforms it to the graph $\Gamma_{{\tau}, m}$ \cite{G}.  
 
 \bl \la{SM}
 Given any ${\cal Q}-$diagram ${\cal Q}$  in $M$ which induces on the boundary the collection ${\cal Z}_\Gamma$ of zig-zag strands of a bipartite graph $\Gamma$ on $S$, 
  and given a two by two move $\mu: \Gamma \lra \Gamma'$, there is a canonical  ${\cal Q}-$diagram ${\cal Q}'$   in $M$ which induces on the boundary 
  the collection ${\cal Z}_{\Gamma'}$ of zig-zag strands of a $\Gamma'$. 
  \el
  
   \begin{proof} Take the standard tetrahedron ${\rm T}$ with the standard collection ${\cal Z}_{\rm T}$ of zig-zag strands on it. A two by two move $\mu_\alpha: \Gamma \lra \Gamma'$ is performed inside of a certain locus $\alpha$ on $S$, shown on Figure \ref{Amut}. The effect of the two by two move  $\mu_\alpha$  on the collection of zig-zag strands ${\cal Z}_\Gamma$ on $S$   can be described by gluing to $S$ the  
   tetrahedra ${\rm T}$  along the two faces of the tetrahedra, matching the zig-zag strands on these two faces with the ones inside of the   locus $\alpha$ on $S$. Then the 
   zig-zag strands on the other two faces of ${\rm T}$ induce the collection of zig-zag strands of $\Gamma'$ inside the locus $\alpha$, see Figure \ref{AB12}. 
   
   On the other hand, there is a canonical ${\cal Q}-$diagram ${\cal Q}_{\rm T}$  in ${\rm T}$ which induces the above collection of  zig-zag strands on the boundary $\partial {\rm T}$. 
   Therefore the effect of the two by two move $\mu_\alpha$ on the ${\cal Q}-$diagram ${\cal Q}$ amounts to gluing to it the  ${\cal Q}-$diagram ${\cal Q}_{\rm T}$. The obtained   
   ${\cal Q}-$diagram ${\cal Q}'$   is the one we need.  \end{proof}

        Given a ${\rm GL}_m-$graph  $\Gamma$ on $S$, there is a sequence of two by two moves which transforms it to the   ${\rm GL}_m-$graph    $\Gamma_{{\rm T}, m}$.
 Applying Lemma \ref{SM} to this sequence of moves we get  the  ${\cal Q}-$diagram ${\cal Q}_\Gamma$.\end{proof}

  \section{Appendices}

 \subsection{A:  Non-commutative cluster ${\cal A}-$varieties} \la{SSEECC5} 
 
 In Section \ref{SSEECC5} we recall briefly some  constructions from \cite{GKo} which we use in the paper.

 A {\it twisted} flat vector bundle over a surface $S$ is a  vector bundle over the punctured tangent bundle $TS- \{\mbox{\rm zero section}\}$  which, for each $p\in S$,  has the monodromy $-{\rm Id}$ 
 around a loop in $T'_pS:= T_pS-\{0\}$ generating $H_1(T'_pS, \Z)$. 
Take a    twisted flat line bundle ${\cal L}$ on an oriented    surface $S$ with boundary $\partial  S$.  Pick a  vector  field $v$ tangent to the boundary   of  $S$, following the boundary orientation. 
A trivialization  of ${\cal L}$ on  the boundary $\partial  S$  is given by  a non-zero section   of the restriction of ${\cal L}$ to $v$.  
\vskip 2mm

 Let   $\Gamma$ be a bipartite ribbon graph, and  $S_\Gamma$ the associated surface.  The conjugate bipartite ribbon graph $\Gamma^*$ 
 provides the spectral surface 
  $
  \Sigma_\Gamma:= S_{\Gamma^*}.
  $ 
Then 
  $ 
 \mbox{  
$\{$zig-zags on $S_\Gamma$$\}$ =  $\{$boundary components of $\Sigma_\Gamma\}$}.
$  Consider the groupoid:
\be
\begin{split}
{\cal A}^\circ_\Gamma:= &\{\mbox{Twisted   flat line bundles on  the spectral 
  surface ${\Sigma_\Gamma }$,   trivialized
at   the boundary}\}\\
  =&\{\mbox{Twisted   flat line bundles on  the  
  surface ${S_\Gamma }$,   trivialized
at the zig-zag strands}\}.\\
\end{split}
\ee

Let us give  a coordinate description of the non-commutative tori 
${\cal A}^\circ_\Gamma$ via the ${\cal A}-$coordinates.

\bd \la{DDD} Let $\Gamma$ be a ribbon graph. The {\em ${\cal A}-$coordinates}\footnote{Calling the elements $\Delta_{E}$ "coordinates" is an abuse of terminology since they are not independent.  } on $\Gamma$ are   the  elements $\{\Delta_E\in R^\times\}$ assigned to 
 the  \underline{oriented} edges $E$ of $\Gamma$,    satisfying  the following monomial relations:

  \begin{itemize}
 \item Let $E$ be an oriented edge, and   $\overline E$ the same edge  with  the opposite orientation. Then
    \be \la{MEQi}
  \Delta_{E}\Delta_{\overline E}=-1.
  \ee  
 \item Let   $E_1, ..., E_n$ be  the edges    incident to any   vertex $v$ of $\Gamma$, oriented out of the vertex $v$, whose order 
is compatible with  their cyclic order. Then     
   \be \la{MEQ}
  \Delta_{E_1}\Delta_{E_2} \ldots \Delta_{E_n}=-1.
  \ee 
\end{itemize}
 \ed

 We   use the name {\it ${\cal A}-$decorated ribbon graph} for a ribbon graph    with a collection of ${\cal A}-$coordinates.   \\

\begin{figure}[ht]
\centerline{\epsfbox{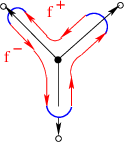}}
\caption{}
\label{ncls00}
\end{figure}

Given a   
  ribbon graph $\Gamma$, the  ${\cal A}-$coordinates on   $\Gamma$ parametrise 
 the isomorphism classes of twisted   flat line bundles on  the  
  surface ${S_\Gamma }$,   trivialized
at   the boundary. 
Indeed, an   edge $E$ of $\Gamma$ determines  two  faces $f_E^+$ and $f_E^-$  of $S_\Gamma$.   
  There is a  canonical homotopy class of path    $p_E$, connecting  $f_E^+$ and $f_E^-$, 
see Figure \ref{ncls00}.  
\noindent  The parallel transport along the   $p_E$   acts on the trivializations $s_E^\pm$ at the boundary components   $f_E^\pm$,  defining   elements $\Delta_{E}\in R^\times$:
 \be \la{para}
 {\rm par}_{p_E}: s^+_{E} \lra \Delta_{E} s^-_{E}, \ \ \ \ \Delta_{E}\in R^\times. 
 \ee 
 Going around any edge, or any vertex along the cyclic order at this vertex, as   on Figure \ref{ncls00},   we rotate the tangent vector by $2\pi$,    getting both relations  (\ref{MEQi}) and (\ref{MEQ}). 
  The $-$ sign   results from the fact that going around the circle amounts to the monodromy $-1$. Conversely, a collection of the ${\cal A}-$coordinates on $\Gamma$ determines a twisted flat line bundle on $S_\Gamma$.

We apply this convention to the bipartite ribbon graph $\Gamma^*$, since   the twisted flat line bundles lives on $S_{\Gamma^*}$.  
Note that  the bipartite ribbon graphs  $\Gamma$ and $\Gamma^*$  are identical as  graphs, but differs by the cyclic order of the edges in the $\bullet-$vertices. 
We  depict the cyclic order on   $\Gamma$ as the  counterclockwise, and 
   use below   the convention  
   for the graph $\Gamma^*$  in  terms of the graph $\Gamma$ as follows: 
 \be \la{CONVa} 
\begin{split}
&\mbox{ \it The counterclockwise product of the elements on the edges sharing   a  $\circ-$vertex is equal to $-1$.}\\
&\mbox{ \it The counterclockwise product of the elements on the edges sharing   a $\bullet-$vertex $v$ is   $(-1)^{{\rm val}(v)-1}$.}\\
\end{split}
\ee 
So describing twisted  flat line bundles on the  surface    $S_{\Gamma^*}$ in the terms of   $\Gamma$ we use convention (\ref{CONVa}).

\subsubsection{Two by two moves of ${\cal A}-$coordinates on bipartite ribbon graphs.}  
\begin{figure}[ht]
\centerline{\epsfbox{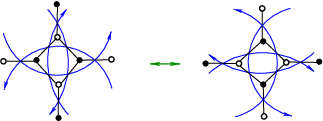}}
\caption{A two by two move of bipartite graphs.}
\label{ncls125}
\end{figure}

\begin{figure}[ht]
\centerline{\epsfbox{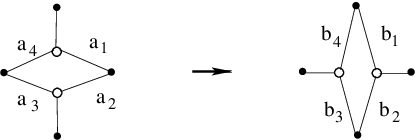}}
\caption{The  two by two move of ${\cal A}-$decorated bipartite ribbon graphs is described by formulas (\ref{FRM1}).} 
\label{Amut}
\end{figure}

Recall the two by two move $\Gamma \to \Gamma'$ of bipartite graphs on Figure \ref{ncls125}.  
Let us define  the corresponding  birational isomorphism  $$
{\cal A}^\circ_\Gamma \lra {\cal A}^\circ_{\Gamma'}.
$$
Consider  coordinates $\{a_1, a_2, a_3, a_4\}$ and $\{b_1, b_2, b_3, b_4\}$ at the internal edges on the graphs on Figure \ref{Amut}. 

Set 
\be
A_i = a_ia_{i+1}a_{i+2}a_{i+3}, ~~~~B_i = b_ib_{i+1}b_{i+2}b_{i+3}, \ \  \ \ i \in \Z/4\Z.
\ee

Consider the following transformation of    the coordinates  illustrated on Figure \ref{Amut}. 
\be \la{FRM1}
 \begin{split}
 &  b_1  =    (1 -A_3^{-1} ) a_3 \ = \ \ \ 
 a_3  (1 - A_4^{-1} );  \\
  &   b_2 =  (1- A_4)^{-1}a_4 = \ \ \ 
 a_4(1-A_1)^{-1};\\
&\mbox{and the cyclic shift $i \lms i+2$ of these two equations.}\\
  \end{split}
 \ee

\begin{figure}[ht]
\centerline{\epsfbox{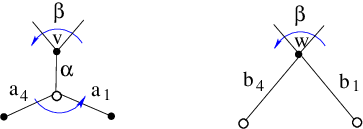}}
\caption{Using conventions (\ref{MEQ}),  $a_4a_1 \alpha=-1$,  $\alpha  \beta= -1$, $b_4b_1\beta= -1$.  
So $  a_4a_1  = -(b_4b_1)^{-1}$.}
\label{AB}
\end{figure}

 Note that   for any $i \in \Z/4\Z$ we have 
 $
\quad a_{i}A_{i+1}  =   A_{i}a_{i}, \quad A_i = B_{i+2}^{-1},\quad 
 $
 and 
\be \la{RAB}
\begin{split}
&b_{i}b_{i+1}= -(a_{i } a_{i+1})^{-1}.   \\ 
\end{split}
\ee
For example,   $b_1b_2a_1a_2= a_3(1-A_4^{-1})(1-A_4)^{-1}a_4a_1a_2  =-1$. 
Alternatively, Figure \ref{AB} shows that  formula  (\ref{RAB}) is consistent with  equation (\ref{MEQ}).
 
 \bp \la{CONSS} \cite{GKo}
 The square of the two by two transformation (\ref{FRM1}) is the identity transformation. 
 The coordinate transformation (\ref{FRM1}) for a two by two move satisfies the pentagon relation.   
 \ep

Let us define the non-comutative cluster variety ${\cal A} = {\cal A}_{|\Gamma|}$ assigned to a bipartite ribbon graph $\Gamma$. Its clusters are parametrised by the bipartite ribbon graphs 
which can be obtained from $\Gamma$ by two by two moves, and  by  shrinking    two-valent vertices 
We assign to each cluster the non-commutative torus ${\cal A}_\Gamma$, and glue them using the 
coordinate transformations (\ref{FRM1}) for the two by two moves. Thanks to Lemma \ref{CONSS}, this definition makes sense.

\subsubsection{ The  canonical 2-form on a  non-commutative cluster ${\cal A}-$variety.}

 The cyclic envelope of the tensor algebra ${\rm T}(V)$ of a graded vector space $V$  is a graded vector space spanned by the elements 
\be
\begin{split}
& v_1   v_2   \ldots   v_n := (v_1 \otimes v_2 \otimes \ldots \otimes v_n)_{\cal C},\\
&(v_1 \otimes v_2 \otimes \ldots \otimes v_n)_{\cal C}  = (-1)^{|v_1| \cdot (|v_2| + \ldots + |v_n|)}(v_2 \otimes \ldots \otimes v_n\otimes v_1)_{\cal C}.\\
\end{split}
\ee
We define the noncommutative analog of $d\log (a_1) \wedge \ldots \wedge d\log (a_n)$ by setting
  \be
  \{a_1, \ldots,  a_n\}:= da_1   \ldots da_n a_n^{-1} \ldots a_1^{-1}. 
\ee

For  a  3-valent  ribbon graph $\Gamma_v$ with a  vertex  $v$,     decorated cyclically by $a, b, c$ with $abc=\pm 1$, 
 we set 
 \be
 \Omega_{v}:= \{a, b\}= da ~db ~b^{-1}a^{-1}.
  \ee
The condition $abc=\pm 1$ implies that it is  invariant under the  
   cyclic shift  $(a,b,c) \lms (b,c,a)$.   Note that 
 \be \la{TRC}
 3 d \{a,b\} = -(a^{-1}da)^3 -  (b^{-1}db)^3 - (c^{-1}dc)^3 .
\ee 
   Given a   ribbon graph $\Gamma_v$ with a single vertex $v$ of valency $>3$, we expand this vertex by adding two-valent vertices of the opposite color, 
   producing a  bipartite ribbon graph $\Gamma'_v$   with trivalent vertices of the original color, and two-valent vertices of the opposite color, and set 
    \be
 \Omega_{v}:=   \sum_{x \in \Gamma'_v} \Omega_{x}.
 \ee   
The  sum is over trivalent  vertices  of the   graph $\Gamma'_v$. 
It does not depend on the choice of   $\Gamma'_v$.

  Given an arbitrary ${\cal A}-$decorated  {bipartite} ribbon graph  $\Gamma$,  we  introduce the  2-form  
     \be \la{AF}
 \Omega_\Gamma:=  \sum_{w \in \Gamma} \Omega_{w} - \sum_{b \in \Gamma} \Omega_{b}.
 \ee   
 Here the first sum is over all $\circ-$vertices   of $\Gamma$, and the second over all $\bullet-$vertices.   
 
 The 2-form $\Omega_\Gamma$  is invariant under the two by two moves. 
This just means that   using   notation (\ref{FRM1}),  
 \be
 \{a_4, a_1\}-\{a_1, a_2\} +  \{a_2, a_3\}-\{a_3, a_4\}  \ \stackrel{}{= }\  \{b_1, b_2\}-\{b_2, b_3\} + \{b_3, b_4\}-\{b_4, b_1\}.   
 \ee
 \vskip 2mm

 We assign to each   edge ${E}$ of an  ${\cal A}-$decorated ribbon graph $\Gamma$, decorated by an element $a_E$  the  3-form 
  \be
 \omega_{E}:= (a_{E}^{-1}da_{E})^3. 
  \ee
Then it follows from (\ref{TRC}) for an  ${\cal A}-$decorated bipartite ribbon graph $\Gamma$, we have:
 \be
d\Omega_\Gamma = \sum_{E: \ \mbox{external edges of $\Gamma$}} \omega_E.
\ee 

\bt The 2-form $\Omega$ is preserved by the two by two move given by formulas (\ref{FRM1})

\et

 \begin{proof} This follows from   \cite[Theorem 4.12]{GKo}. Note however that formulas (\ref{FRM1}) where written in \cite{GKo} with different sign:
 \be \la{FRM2}
 \begin{split}
 &  b_1  =    (1 +A_3^{-1} ) a_3;\ \ \ \    b_2 =  (1+A_4)^{-1}a_4; \\
 &  b_3  =   (1 +A_1^{-1} )a_1; \ \ \ \  b_4 =  (1+A_2)^{-1}a_2.\\
  \end{split}
 \ee 
 However the 2-forms do not change if we replace $a_i$ by $-a_i$, and similarly for $b_i$. So changing $a_1\lms -a_1, b_3 \lms -b_3$, and keeping the rest six variables intact we transform 
 (\ref{FRM2}) to (\ref{FRM1}).\end{proof}

 \vskip 2mm
 {\it Comments on Proposition \ref{CONSS}.} The issue  is that we use formulas (\ref{FRM1}) rather than (\ref{FRM2}).  Yet the first claim follows by the same trick with the sign change as above. 
 
 \subsection{B: Admissible dg-sheaves} \la{SEC7.1}

  We use the dg-enhamncement of the derived category of complexes of sheaves, constructible with respect to a stratification on $X$ from \cite[Section 10.1]{GKo}.
  Its objects are called dg-sheaves.  
  
  Let ${\cal H}$ be a collection   of    cooriented hypersurfaces in a manifold $X$. 
A complex of  constructible sheaves on $X$ is ${\cal H}$-supported    
 if its microlocal  support is contained in  the union of the zero section of $T^*X$ and the conormal bundles to the 
  cooriented hypersurfaces in  ${\cal H}$. 
 
  Below we assume that the hypersurfaces in  ${\cal H}$ have normal crossing intersections. 
Among   ${\cal H}$-{supported} dg-sheaves, we distinguish the  subcategory of {\it ${\cal H}$-admissible dg-sheaves} introduced in \cite[Section 10.1]{GKo}.  
Below we spell the definition in the important for us cases when ${\rm dim}(X) = 2,3$.

An ${\cal H}$-admissible  dg-sheaf  on a threefold $M$ is given by the following data:

\begin{itemize}
  
  \item A   complex of local systems    ${\cal F}^\bullet_{\cal D}$ on each domain  ${\cal D}$ of  $M-{\cal H}$, concentrated in   degrees $[0, 1]$:

 \item  For each  component $\gamma$  of ${\cal H} - $ $\{$codimension $\geq 2$ strata of ${\cal H}$$\}$, sharing  domains  ${\cal D}^+$, and  ${\cal D}^-$ 
 and cooriented towards ${\cal D}^+$, a  map  of complexes  
   \be \la{CONEGx}
   \begin{split}
&{\varphi}_{\gamma}:  {\cal F}^\bullet_{{\cal D}^-} \to {\cal F}^\bullet_{{\cal D}^+},\\
&H^i{\rm Cone}(\varphi_{\gamma}) =0~~\mbox{if  $i\not =0$}.\\
\end{split}
\ee

 \begin{figure}[ht]
\centerline{\epsfbox{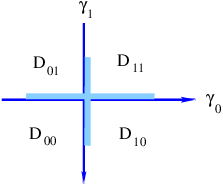}}
\caption{The stratification near a codimension two intersection of  cooriented hyperplanes $\gamma_0, \gamma_1$  has a codinension two strata $\gamma_0\cap \gamma_1$, four  codimension one strata, and four domains  ${\cal D}_{\ast\ast}$.}
\label{ncls102a}
\end{figure}

\item For each codimension two strata, given by intersection of two codimension one strata, a homotopy: 
 $$
 h: {\cal F}^\bullet_{{\cal D}_{00}}  \lra {\cal F}^\bullet_{{\cal D}_{11}}[-1], ~~~~dh  + h d =  \varphi_{2} \circ \varphi_{1} -\varphi_4 \circ \varphi_3.
  $$
 \be \la{rs22a} 
\begin{gathered}
    \xymatrix{ 
        {\cal F}^\bullet_{{\cal D}_{01}}\ar[r]^{\varphi_{2}}    & {\cal F}^\bullet_{{\cal D}_{11}}  \\
         {\cal F}^\bullet_{{\cal D}_{00}}\ar[u]^{\varphi_{1}}  \ar[ru]^{h} \ar[r]_{\varphi_{3}}        & \ar[u]_{\varphi_{4}}  {\cal F}^\bullet_{{\cal D}_{10}}}
\end{gathered}
 \ee
 
 \item Similarly, for each codimension three strata, a higher homotopy
 $$
 h_2: {\cal F}^\bullet_{{\cal D}_{000}}  \lra {\cal F}^\bullet_{{\cal D}_{111}}[-2], ~~~~dh_2  + h_2 d =  \ldots .
  $$
 
\item  Let    ${\rm D}$   be   the signed sum of the differentials on ${\cal F}^\bullet_{{\cal D}_{\ast\ast\ast}}$,  the maps  $\varphi_{\ast}$, homotopies $h$, and higher homotopy $h_2$.  Then there is 
the  following complex, given by the direct sum of the stalks assigned to the vertices of a cube with the differential ${\rm D}$:
 \be \la{CONDg}
  \begin{split}
 & {\cal F}^\bullet_{{\cal D}_{000}} \oplus \Bigl({\cal F}^\bullet_{{\cal D}_{001}} \oplus {\cal F}^\bullet_{{\cal D}_{010}}\oplus  {\cal F}^\bullet_{{\cal D}_{001}}\Bigr)[-1]\oplus \Bigl({\cal F}^\bullet_{{\cal D}_{011}} \oplus {\cal F}^\bullet_{{\cal D}_{101}}\oplus  {\cal F}^\bullet_{{\cal D}_{011}}\Bigr)[-2]  \oplus
 {\cal F}^\bullet_{{\cal D}_{111}}[-3]. \\
 \end{split}
 \ee
  \end{itemize}
  
  We assume that this data satisfies the folloing condition
\be \la{CAC}
\mbox{Complexes (\ref{rs22a}) and (\ref{CONDg}) are acyclic.} 
\ee

When  $X$ is a surface,   there are no higher homotopies $h_2$, and the data  (\ref{CONDg}) is void.   

\bp Conditions (\ref{CAC}) are equivalent to   conditions that the microlocal support of ${\cal F}$ does not contain any vectors in 
the conormal bundles to   codimension two and three strata minus $\cup_iT^*_{H_i}X$, respectively. 
\ep

When $X$ is a surface, this follows from   \cite[Lemma 10.5]{GKo} and the discussion there. 

The case of the threefold,   involving the condition that the complex (\ref{CONDg}) is acyclic, is similar. 

Precisely, 
given a point $q\in X$, denote by $T^+_{q, {\cal H}}X$ the open  cone in $T^*_qX$, given by 
positive linear combinations of  non-zero covectors from $T^*_{H_i}$, where $q \in H_i$.   Then we have the following Lemma. 

\bl
Given an ${\cal H}$-supported   complex of sheaves ${\cal F}$ on $X$, and a codimension three stratum $q$ of ${\cal H}$, the  fiber of the microlocal support of ${\cal F}$  at any  $\eta 
\in T^+_{q, {\cal H}}X$ is quasiisomorphic to   complex  (\ref{CONDg}). 
 \el
 
 One easily sees, just as in Section 10 of loc. cit.,  that this implies the claim. 
 
So  the  restriction of an ${\cal H}$-admissible dg-sheaf    to the conormal bundle to each surface   of ${\cal H}$ is a local system. 
  We say that a dg-sheaf has {\it rank one microlocal support} if these local systems, sitting in the degree zero, 
  are one dimensional. 
 The microlocal support on the zero section of $X$ can be more complicated.


\begin{thebibliography}{BL}
   
   \bibitem[B]{B} Beilinson A.: {\it Higher regulators and values of L-functions},
J. of Soviet Math. 30 (1985), 2036-2070. 
   
   \bibitem[BR]{BR} Berenstein A., Retakh V.: {\it Noncommutative marked surfaces}.     \newblock \href{https://arxiv.org/abs/math/1510.02628}{arXiv:1510.02628}.  
\bibitem[BrD1]{BrD1} Brav C., and Dyckerhoff T.:.{\it Relative Calabi-Yau structures}. 
Compositio Mathematica 155.2 (2019), pp. 372–412. 
\bibitem[BrD2]{BrD2} Brav C., and Dyckerhoff T.:  {\it Relative Calabi-Yau structures II: shifted Lagrangians in the moduli of objects}. Selecta Mathematica 27.4 (2021), p. 63.
\bibitem[D]{D} Deligne P.: {\it Le symbole mod\'er\'e}. Publ. Math\'ematiques de IH\'ES, Volume 73 (1991), pp. 147-181.
\bibitem[DGG]{DGG} Dimofte T., Gabella M., Goncharov A.B.: {\it K-Decompositions and 3d Gauge Theories} 
~\newblock \href{https://arxiv.org/abs/math/1301.0192}{arXiv:1301.0192}.    

\bibitem[FG1]{FG1} Fock V., Goncharov A.: 
{\it Moduli spaces of local systems  and higher Teichm\"uller theory}. Publ. Math. 
IHES, n. 103 (2006) 1-212. ~\newblock \href{https://arxiv.org/abs/math/0311149}{arXive:0311149}. 
\bibitem[FG2]{FG2} Fock V., Goncharov A.: {\it Cluster ensembles, quantization and the dilogarithm.} 
Ann. Sci. L'Ecole Norm. Sup. (2009). ~\newblock \href{https://arxiv.org/abs/math/0311245}{arXive:0311245}.
\bibitem[FN]{FN}  Freed D., Neitzke A.: {\it 3d spectral networks and classical Chern-Simons theory}. 
~\newblock \href{https://arxiv.org/abs/math/2208.07420 }{arXiv:2208.07420}. 


\bibitem[GPS1]{GPS1}  Ganatra S., Pardon J., Shende V.: {\it Covariantly functorial wrapped floer theory on Liouville sectors}, Publ. Math. IHES 131 (2020), 73–200.
\bibitem[GPS2]{GPS2}  Ganatra S., Pardon J., Shende V.:  {\it Microlocal morse theory of wrapped Fukaya categories}, Ann. of Math. (2) 199 (2024), no. 3, 943–1042.
\bibitem[GPS3]{GPS3}  Ganatra S., Pardon J., Shende V.:  {\it Sectorial descent for wrapped Fukaya categories}, J. Amer. Math. Soc. 37 (2024), no. 2, 499–635. 
\bibitem[GSV]{GSV} Gekhtman M.,  Shapiro M.,  Vainshtein A., {\it Cluster algebras and Weil-Petersson forms} ~\newblock \href{https://arxiv.org/abs/math/0309138}{arXiv:0309138}.
\bibitem[G]{G} Goncharov A.B. {\it Ideal webs moduli spaces of local systems, and 3d Calabi-Yau categories}.     In "Algebra, Geometry, and Physics in the $21^{\rm st}$  Century". Kontsevich Festshrift. 
    Prog. in Math., 324, 2017, p. 31-99. ~\newblock \href{https://arxiv.org/abs/math/1607.05228}{arXive:1607.05228}.  
\bibitem[GKe]{GKe} Goncharov A.B., Kenyon R.: {\it Dimers and cluster integrable systems}. Ann. Sci. L'Ecole Norm. Sup., 2013, vol. 46, n 5, 747-813, 
~\newblock \href{https://arxiv.org/abs/math/107.5588}{arXive:1107.5588}.
\bibitem[GKo]{GKo} Goncharov A.B., Kontsevich M.: 
{\it Spectral description of non-commutative local systems on surfaces and non-commutative cluster varieties.} ~\newblock \href{https://arxiv.org/abs/math/107.5588}{arXiv:2108.04168}.  

\bibitem[GKS]{GKS}   Guillermou S., Kashiwara M.,  Schapira P: {\it Sheaf Quantization of Hamiltonian Isotopies and Applications to Nondisplaceability Problems}, Duke Math J. 161 (2012) 201-245.
\bibitem[KS]{KS} Kashiwara M.,  P. Schapira P.: {\it Sheaves on Manifolds}, Grundlehren der Mathematischen Wissenschafte
292, (Springer-Verlag, 1994).
\bibitem[K]{K}  Kontsevich M.: {\it Symplectic geometry of homological algebra}. Arbeitstagung, Bonn, 2009.  
\bibitem[N]{N} Nadler D.: {\it Arboreal singularities} ~\newblock \href{https://arxiv.org/abs/math/1309.4102}{arXive:1309.4122}.
\bibitem[STZ]{STZ} Shende V., Treumann D.,  Williams H.: {\it Legendrian knots and constructible sheaves}. Inventiones mathematicae 207 (3), 1031-1133.
~\newblock \href{https://arxiv.org/abs/math/1402.0490}{arXive:1402.0490}.  \bibitem[STWZ]{STWZ} Shende V., Treumann D.,  Williams H., Zaslow E.: {\it Cluster varieties from Legendrian knots}. Duke Mathematical Journal 168 (15), 2801-2871~\newblock \href{https://arxiv.org/abs/math/1512.08942}{arXive:1512.08942}. 
\bibitem[STW]{STW} Shende V., Treumann D.,  Williams H.: {\it On the combinatorics of exact Lagrangian surfaces}. ~\newblock \href{https://arxiv.org/abs/math/1603.07449}{arXive:1603.07449}.   
\bibitem[Th]{Th} Thurston D.: {\it From domino to hexagons}. ~\newblock \href{https://arxiv.org/abs/math/0405482}{arXive:0405482}.
\end{thebibliography}
\end{document}